\documentclass[11pt,a4paper]{report}

\usepackage{theorem}
\usepackage[dvips]{epsfig}
\usepackage{amsmath}
\usepackage{amssymb}
\usepackage{amscd}
\usepackage[all]{xy} \CompileMatrices  
\usepackage{makeidx} \makeindex
\theoremstyle{plain}
\theorembodyfont{\rmfamily}
\newtheorem{lemma}{Lemma}[chapter]
\newtheorem{defn}[lemma]{Definition}
\newtheorem{rem}[lemma]{Remark}
\newtheorem{prop}[lemma]{Proposition}
\newtheorem{example}[lemma]{Example}
\newtheorem{thm}[lemma]{Theorem}
\newtheorem{cor}[lemma]{Corollary}
\newtheorem{question}[lemma]{Question}

\newcommand{\qua}{\nopagebreak \hspace*{\fill} $\Box$}
\newcommand{\asdim}{\operatorname{asdim}}
\newcommand{\im}{\operatorname{im}}

\newcommand{\diam}{\operatorname{diam}}
\newcommand{\mesh}{\operatorname{mesh}}
\newcommand{\abs}[1]{\lvert#1\rvert}
\newcommand{\norm}[1]{\lVert#1\rVert}
\newcommand{\RRR}{\mathbb{R}}
\newcommand{\ZZZ}{\mathbb{Z}}
\newcommand{\NNN}{\mathbb{N}}

\newcommand{\CCC}{\mathbb{C}}
\newcommand{\cs}{\operatorname{C\!\!S}}
\newcommand{\cncs}{\operatorname{_{cn}\!C\!\!S}}
\newcommand{\compl}{\operatorname{\square\!\!\!\!\!\!\:\angle}}
\newcommand{\vertices}{\operatorname{vert}}
\newcommand{\CAT}{\operatorname{CAT}}
\newcommand{\Star}{\operatorname{Star}}
\newcommand{\supp}{\operatorname{Supp}}

\newcommand{\id}{\operatorname{id}}

\def\kringel{\mathaccent "7017}

\usepackage[hypertex,backref]{hyperref}

\title{Coarse geometry and asymptotic dimension\\[2cm]}
\author{Dissertation\\
  zur Erlangung des Doktorgrades\\ 
  der Mathematisch-Naturwissenschaftlichen Fakult\"aten\\
  der Georg-August-Universit\"at zu G\"ottingen\\[3cm]
  vorgelegt von\\
  \textbf{Bernd Grave}\\ 
  aus\\
  Damme\\[5cm]}
\date{G\"ottingen 2005}

\begin{document}

\maketitle

\ \vfill
D7
\begin{description}
\item[Referent:] Prof. Dr. Thomas Schick
\item[Korreferent:] Prof. Dr. Ulrich Stuhler
\item[Tag der m\"undlichen Pr\"ufung:] Montag, 23. Januar 2006
\end{description}
\thispagestyle{empty}\newpage

\begin{abstract}
We prove that two spaces, whose coarse structures are induced by metrisable compactifications,
are coarsely equivalent if and only if their (Higson) coronas are homeomorphic.

\bigskip

We introduce translation $C^*$-algebras for coarse spaces
which admit a countable, uniformly bounded cover
using pro\-jec\-tion-valued measures. This was already done in \cite{RoeCG}.
Here we give a more complete exposition on the subject
including a result about the independence of the translation $C^*$-algebra from the involved projection-valued measure.
Moreover, we eliminate some mistakes, e.g. we have to be careful about null sets when defining the support of an operator.

\bigskip

We introduce some characterisations of asymptotic dimension in the general setting of coarse spaces
and prove some basic properties such as monotony, a formula for the asymptotic dimension of finite unions
and estimates for the asymptotic dimension of the product of two coarse spaces.

\bigskip

We define coarse cell complexes and prove the obvious conjecture about their asymptotic dimension.

\bigskip

We prove $\asdim(X,\mathcal{E}) = \dim(K\backslash X) + 1$
if the coarse structure $\mathcal{E}$ is induced by a metrisable compactification $K$ of $X$.
As an application, we obtain coarse spaces $X$ and $Y$ such that
$\asdim(X\times Y) < \asdim(X) + \asdim(Y)$.

\bigskip

For $\text{CAT}(\kappa)$-spaces with $\kappa < 0$ having nicely $n$-covered spheres we prove
that the asymptotic dimension is at most $n$.
We apply this result to prove that $\asdim(X) = \dim(X)$ for any complete, simply connected
Riemannian manifold $X$ with bounded, strictly negative sectional curvature.
\end{abstract}

\thispagestyle{empty}\ \newpage 
\thispagestyle{empty}
\vspace*{\fill}
\subsection*{\centering Acknowledgements}
\thispagestyle{empty}
\vspace{6pt}

First of all I wish to thank my advisor, Thomas Schick, for his support,
many enlightening discussions and most especially his patience when I was struggling.

\bigskip

The "Graduiertenkolleg Gruppen und Geometrie" at G\"ottingen University provided the financial support
and mathematical education which benefited me greatly in this endeavor.
And I profited as well from visiting the research group
in Noncommutative Geometry and Operator Algebras at Vanderbilt University.

\bigskip

Special thanks to Paul Mitchener for inspiring discussions about coarse geometry,
his many helpful suggestions and his very special English lessons,
in addition to those taught by Marjory Frauts.

\bigskip

I am also grateful to Elias Kappos, Behnam Norouzizadeh and Moritz Wiet\-haup
for their careful read of this dissertation and their detailed feedback.

\bigskip

Last, but not least, I wish to thank all my friends both in the G\"ottingen math department
and without for their kindness and the part they play in making G\"ottingen
such an attractive place in which to work and live.

\vspace*{\fill}
\vspace*{\fill}
\clearpage
\thispagestyle{empty}\ \newpage 

\pagenumbering{arabic} 

\tableofcontents
\subsubsection*{\hspace{17pt}Bibliography \hspace*{\fill}73}
\subsubsection*{\hspace{17pt}Index\hspace*{\fill}77}
\subsubsection*{\hspace{17pt}Curriculum vitae\hspace*{\fill}79}

\clearpage
\thispagestyle{empty}\ \newpage 

\chapter{Some basics in coarse geometry}

\section[Finitely generated groups from the metric viewpoint]{Finitely generated groups \\ from the metric viewpoint}

In this section we give an example of a setting where coarse geometry naturally arises. 
We define the word metric on each finitely generated group.
This metric depends on the chosen generating set,
but - as we will see soon - the asymptotic behavior does not.
So we can associate a coarse structure to each finitely generated group.

\bigskip

Let $G$ be a group with a finite generating set $A$.
\begin{defn}[word metric] \label{word_metric} \index{word metric}
The distance $d_{(G,A)}(g_1,g_2)$ of $g_1, g_2\in G$ in the word metric associated to the pair
$(G,A)$ is the length of the shortest word in $A$ representing $g_1^{-1}g_2$.
\end{defn}
Alternatively we could have used the Cayley graph of $(G,A)$ to define this metric. Compare
\cite{BH} for more details.

\begin{lemma} \label{wordlemma}
If $B$ is another finite generating set of $G$, there is $\lambda > 0$ such that
$d_{(G,A)} \leq \lambda \cdot d_{(G,B)}$.
\end{lemma}
\textbf{Proof.}
For each $b \in B$ choose a word in $A$ of minimal length representing $b$.
Define $\lambda$ to be the maximum of the length of these words.
Given any word in $B$ representing $g \in G$,
we replace each letter by the chosen word in $A$ and count letters.
\qua

\begin{defn}\label{metricmapproperties}
Let $X,Y$ be pseudometric\footnote{Omitting the condition $d(x_1,x_2)>0$
if $x_1 \neq x_2$ in the definition of a metric space,
we get the definition of a pseudometric space.
In a pseudometric space being of distance zero is an equivalence relation.
Equivalent points are contained in exactly the same open sets.
The map to the quotient can be used to carry over many properties of metric spaces
to pseudometric spaces. From the asymptotic point of view
pseudometric spaces should be as good as metric spaces.}
spaces and $f \colon X \to Y$ a (not necessarily continuous) map.
\begin{itemize}
\item
$f$ is called \emph{coarsely proper}\index{coarsely proper}
if the inverse image of any bounded set is bounded.
\item
$f$ is called \emph{coarsely uniform}\index{coarsely uniform}
if for every $r>0$ there is $s(r)>0$ such that
$d(f(x_1),f(x_2)) \leq s(r)$ for all $x_1,x_2 \in X$ with $d(x_1,x_2) \leq r$.
\item
$f$ is called a \emph{coarse map}\index{coarse map} if it is coarsely proper and coarsely uniform.
\item
Let $S$ be a set. Two maps $f,g \colon S \to X$ are called \emph{close}\index{close maps}
if there is $D>0$ such that $d(f(s),g(s)) < D$ for all $s \in S$.
\item
$f$ is called a \emph{coarse equivalence} if it is a coarse map and if there exists another coarse map
$g \colon Y \to X$ such that $g \circ f$ is close to $id_X$ and $f \circ g$ is close to $id_Y$.
\item
$X$ and $Y$ are called \emph{coarsely equivalent}\index{coarsely equivalent}
if there exists a coarse equivalence from $X$ to $Y$.
\end{itemize}
Compare our coarse maps with asymptotically Lipschitz maps in \cite{Dran}.
\end{defn}

Consider a finitely generated group.
The metric spaces corresponding to different finite generating sets
are obviously coarsely equivalent.
Thus, every functor from the category of metric spaces and arbitrary maps
which is invariant under coarse equivalence gives an invariant of finitely generated groups.
The asymptotic dimension will be an example of such an invariant.

\section{The coarse category}

In this section an axiomatic description of the structure
needed to do coarse geometry will be given.
Compare \cite{Roe2}, \cite{Mitchener} and \cite{RoeCG}.

\begin{defn}[coarse structure] \label{def_coarse_space}\index{coarse structure}
Let X be a set. A collection $\mathcal{E}$ of subsets of $X \times X$
is called a \emph{coarse structure}, and the elements of $\mathcal{E}$ will be called
\emph{entourages}\index{entourage}, if the following axioms are fulfilled:
\begin{enumerate}
\item[(a)] A subset of an entourage is an entourage.
\item[(b)] A finite union of entourages is an entourage.
\item[(c)] The diagonal $\Delta_X := \{(x,x) \mid x \in X \}$ is an entourage.\index{diagonal of a set}
\item[(d)] The inverse $E^{-1}$ of an entourage $E$ is an entourage.
$$E^{-1} := \{(y,x) \in X\times X \mid (x,y)\in E\}$$
\item[(e)] The composition $E_1E_2$ of entourages $E_1$ and $E_2$ is an entourage.
\index{composition of entourages}
$$E_1E_2 :=\{(x,z)\in X\times X \mid \exists_{y\in X} (x,y)\in E_1 \text{ and } (y,z)\in E_2\}$$
\end{enumerate}
The pair $(X,\mathcal{E})$ is called a \emph{coarse space}.\index{coarse space}
Sometimes we will say that $E\subseteq X\times X$ is \emph{controlled} if $E$ is an entourage.\index{controlled}

A coarse space is called \emph{connected}\index{connected coarse structure}
if every point of $X\times X$ is contained in an entourage.
\end{defn}

\begin{rem} For $E_1,E_2,F_1,F_2\subseteq X\times X$ we have
$$(E_1\cup E_2)(F_1\cup F_2) = E_1F_1\cup E_1F_2\cup E_2F_1\cup E_2F_2 \ . $$
\end{rem}

\begin{defn}[bounded sets]\index{bounded}\index{E@$E[...]$}
Let $(X,\mathcal{E})$ be a coarse space, $A\subseteq X$ and $E\in\mathcal{E}$.
We define
$$E[A]:=\{x\in X\mid (x,a)\in E \text{ for some } a\in A\}\ .$$
For a point $x\in X$ we will write $E(x)$ instead of $E[\{x\}]$.\index{E@$E(...)$}
Sets of the form $E(x)$ with $x\in X$ and $E\in \mathcal{E}$ are called \emph{bounded}.  
\end{defn}

\begin{prop}[properties of bounded sets]
Let $(X,\mathcal{E})$ be a coarse space.
\begin{enumerate}
\item Subsets of bounded sets are bounded.
\item If $B\subseteq X$ is bounded, then $B\times B\in\mathcal{E}$.
\item If $B\subseteq X$ is bounded and $E\in\mathcal{E}$, then $E[B]$ is bounded.
\item Let $B_1,B_2\subseteq X$ be bounded sets. The following are equivalent.
\begin{itemize}
\item $B_1\cup B_2$ is bounded.
\item $B_1\times B_2 \in \mathcal{E}$
\item There exists an entourage $E\in\mathcal{E}$ such that $E \cap B_1\!\times\! B_2 \neq \emptyset$.
\end{itemize}
If $(X,\mathcal{E})$ is a connected coarse space, then any finite union of bounded sets is bounded.
\end{enumerate}
\end{prop}
\textbf{Proof.}
Observe that for entourages $E_1, E_2\in\mathcal{E}$ and $A \subseteq X$ we have $E_1E_2[A]=E_1[E_2[A]]$.
Therefore $E[B]$ is bounded if $E\in \mathcal{E}$ and $B\subseteq X$ is bounded.

Let $B_1, B_2\subseteq X$ be bounded sets
and let $E\in\mathcal{E}$ such that $E \cap B_1\!\times\! B_2 \neq \emptyset$.
Take $(b_1,b_2)\in E \cap B_1\!\times\! B_2$ and
observe that $b_1,b_2\in (E\cup\Delta_X)(b_2)$
and $B_1\cup B_2\subseteq (B_1\!\times\! B_1\cup B_2\!\times\! B_2)(E\cup\Delta_X)(b_2)$.
Hence $B_1\cup B_2$ is bounded.
\qua

\begin{defn}\index{generated coarse structure}\index{CS@$\cs$}\index{cnCS@$\cncs$}
Let $X$ be a set and $\mathcal{M}$ a collection of subsets of $X\times X$.
Since any intersection of coarse structures on $X$ is itself a coarse structure,
we can make the following definition.
By $\cs(\mathcal{M})$ we denote the smallest coarse structure containing $\mathcal{M}$,
i.e. the intersection of all coarse structures containing $\mathcal{M}$.
We call $\cs(\mathcal{M})$ the \emph{coarse structure generated by $\mathcal{M}$}.

In the same way, define the connected coarse structure generated by $\mathcal{M}$
and denote it by $\cncs(\mathcal{M})$.
\end{defn}

\begin{rem}\label{explcs}
There is an explicit description of the coarse structure $\cs(\mathcal{M})$
generated by $\mathcal{M}$.
Set $\mathcal{M}_0 := \mathcal{M} \cup \{\Delta_X\}$
and for $n\in\NNN$ define $\mathcal{M}_n$ inductively as follows.
\begin{eqnarray*}
\mathcal{M}_{n+1} & := & \{ M_1 \mid M_2\in\mathcal{M}_n, M_1\subseteq M_2 \}
							  \ \cup \ \{ M^{-1} \mid M\in\mathcal{M}_n \} \\
							 	& \cup & \{ M_1\cup M_2 \mid M_1,M_2\in\mathcal{M}_n \}
							  \ \cup \ \{ M_1M_2 \mid M_1,M_2\in\mathcal{M}_n \}
\end{eqnarray*}
Obviously, $\mathcal{M}_{\infty} := \bigcup_{n\in\NNN} \mathcal{M}_n$ is contained in $\cs(\mathcal{M})$
and since $\mathcal{M}_{\infty}$ is indeed a coarse structure,
we get $\cs(\mathcal{M}) = \mathcal{M}_{\infty}$.
\end{rem}

\begin{rem}
Let $(X,\mathcal{E})$ be a coarse structure and define $\mathcal{E}_{cn} := \cncs(\mathcal{E})$.
For $E\in\mathcal{E}_{cn}$ there are $E'\in\mathcal{E}$, $k\in\NNN$
and sets $A_1,\ldots,A_k,B_1,\ldots,B_k\subseteq X$ which are bounded with respect to $\mathcal{E}$
such that
$$E\ =\ E'\ \cup\ A_1\!\times\! B_1\ \cup\ \cdots\ \cup\ A_k\!\times\! B_k \ .$$

Any set which is bounded with respect to $\mathcal{E}_{cn}$ is a finite disjoint union
of sets which are bounded with respect to $\mathcal{E}$.
\end{rem}

There are notions of compatibility of coarse structure and topology.
\begin{defn}\index{compatibility of coarse structure and topology}\label{CSandTopcomp}
Suppose we are given a topological space $X$.
A coarse structure $\mathcal{E}$ on $X$ is said to be \emph{compatible with the topology}
if
(1) there is a neighborhood of the diagonal $\Delta_X$ which is an entourage and
(2) the closure of any bounded set is compact.
\end{defn}
In \cite{RoeCG} and \cite{HR}, the term proper\index{proper} is being used
for compatibility of coarse structure and topology.
Compatibility of coarse structure and topology in the sense of \cite{Mitchener},
implies compatibility as defined above.

\begin{example}\textbf{(bounded coarse structure)}\index{bounded coarse structure} \\
Let $(X,d)$ be a pseudometric space.
Set $\Delta_r := \{(x,y)\in X\times X\mid d(x,y)<r\}$ and define
$$\mathcal{E}_d := \cs(\{\Delta_r\mid r>0\})
   = \{ E\subseteq X\times X \mid E \subseteq \Delta_r \text{ for some } r>0 \}.$$
It is easy to verify that $(X,\mathcal{E}_d)$ is
a connected coarse space compatible with the topology.
\end{example}

Let $G$ be a finitely generated group. Lemma~\ref{wordlemma} tells us
that word metrics arising from different finite generating systems
induce the same bounded coarse structure $\mathcal{E}_G$.\index{EG@$\mathcal{E}_G$}
Thus, in a natural manner, every finitely generated group $G$ is a coarse space.

\begin{example}\textbf{(continuously controlled coarse structure)}
\index{continuously controlled}\index{compactification} \\
Let $X$ be a Hausdorff space and $\overline{X}$ a compactification of $X$,
i.e. $X$ is a dense and open subset of the compact set $\overline{X}$.
The collection
$$\mathcal{E}_{\overline{X}} := \{E\subseteq X\times X \mid \overline{E}\subseteq X\times X\cup
\Delta_{\overline{X}}\}$$
of all subsets $E\subseteq X\!\times\! X$, whose closure meets the boundary
$(\overline{X}\!\times\!\overline{X})\backslash(X\!\times\! X)$
only in the diagonal, is a connected coarse structure on $X$.
If $\overline{X}$ is metrisable, the coarse structure $\mathcal{E}_{\overline{X}}$
is compatible with the topology.
The proof is not hard. Compare \cite{RoeCG}.
\end{example}

\begin{defn}[close maps]\index{close maps}
Let $(X,\mathcal{E})$ be a coarse space and $S$ a set.
The maps $f \colon S \to X$ and $g \colon S \to X$ are called \emph{close}
if $\{(f(s),g(s)) \mid s\in S)\}$ is an entourage.
\end{defn}
Compare with Definition \ref{metricmapproperties}.

\begin{defn}
Let $X$, $Y$ be coarse spaces and $f\colon X\to Y$ a map. 

\begin{itemize}
\item
We call $f$ \emph{coarsely proper}\index{coarsely proper}
if the inverse image of bounded sets is bounded.
\item
We call $f$ \emph{coarsely uniform}\index{coarsely uniform}
if the image of each entourage under the map
$f\times f \colon X\times X \to Y\times Y$ is an entourage.
\item
We call $f$ a \emph{coarse map}\index{coarse map}
if it is coarsely proper and coarsely uniform.
\item
We call $f$ a \emph{coarse embedding}\index{coarse embedding}
if $f$ is coarsely uniform and
the inverse image of an entourage under $f\times f$ is an entourage.\footnote{
In \cite{RoeCG}, coarse embeddings are also called rough maps.}

\end{itemize}
Note that a coarse embedding is a coarse map.
\end{defn}

We will denote the category of coarse spaces and coarse maps by $\mathcal{C}$.
\index{C@$\mathcal{C}$}
This is the category used in coarse algebraic topology.
Compare \cite{Roe2}, \cite{RoeCG} and \cite{Mitchener}.

For some constructions it turns out to be more convenient
to work in the category of coarse spaces and coarsely uniform maps.
We will denote this category by $\mathcal{D}$.
\index{D@$\mathcal{D}$}
When refering only to connected coarse structures, we will write
$\mathcal{C}_{cn}$\index{Ccn@$\mathcal{C}_{cn}$} 
and $\mathcal{D}_{cn}$\index{Dcn@$\mathcal{D}_{cn}$} respectively.
Observe that there are some forgetful functors between these categories.

A difference between the categories $\mathcal{C}$ and $\mathcal{D}$ is
that in $\mathcal{D}$ products and direct limits exist,
while in $\mathcal{C}$ they do not.
We will deal with products, direct limits, etc. in Section~\ref{CGconstr}.

\begin{defn}\label{defncoarseequiv}
Let $X$, $Y$ be coarse spaces and $f\colon X\to Y$ a map.
\begin{itemize}
\item\index{D-equivalenz@$\mathcal{D}$-equivalence}
$f$ is called a $\mathcal{D}$\emph{-equivalence} if $f$ is coarsly uniform
and there exists a coarsely uniform map $g \colon Y \to X$
such that $g \circ f$ is close to $id_X$ and $f \circ g$ is close to $id_Y$.
\item\index{C-equivalenz@$\mathcal{C}$-equivalence}
$f$ is called a $\mathcal{C}$\emph{-equivalence} if $f$ is coarse
and there exists a coarse map $g \colon Y \to X$
such that $g \circ f$ is close to $id_X$ and $f \circ g$ is close to $id_Y$.
\end{itemize}
\end{defn}

\begin{prop}
Any $\mathcal{C}$-equivalence is also a $\mathcal{D}$-equivalence.
Conversely, any $\mathcal{D}$-equivalence $f$ is a coarse embedding
and hence a $\mathcal{C}$-equivalence.
\end{prop}
\textbf{Proof.}
Obviously, any $\mathcal{C}$-equivalence is a $\mathcal{D}$-equivalence.

Let $\mathcal{E}$ and $\mathcal{E}'$ be the coarse structures for $X$ and $Y$ respectively.
We prove that a $\mathcal{D}$-equivalence $f\colon X\to Y$ is a coarse embedding.
Let $g\colon Y\to X$ be a coarsely uniform map as in Definition~\ref{defncoarseequiv}.
Take $M'\in\mathcal{E}'$ and set $M:=g\times g(M')$.
Define $E:=\{(x,g\circ f(x))\mid x\in X\}\in\mathcal{E}$.
It is sufficient to prove $(f\times f)^{-1}(M')\subseteq EME^{-1}$.
Let $(x_1,x_2)\in(f\times f)^{-1}(M')$.
With $y_i:=f(x_i)$ we have $(y_1,y_2)\in M'$.
Set $\tilde{x}_i:=g(y_i)=g\circ f(x_i)$.
Using $(x_i,\tilde{x}_i)\in E$ and $(\tilde{x}_1,\tilde{x}_2)\in M$, we get $(x_1,x_2)\in EME^{-1}$.
\qua

\begin{defn}\textbf{(coarse equivalence)}\index{coarsely equivalent} \\
We call $f\colon X\to Y$ a \emph{coarse equivalence} if $f$ is a $\mathcal{C}$-equivalence.
A map $g\colon Y\to X$ as in Definition~\ref{defncoarseequiv}
is called a \emph{coarse inverse}\index{coarse inverse} of $f$.
We say that $X$ and $Y$ are \emph{coarsely equivalent}
if there exists a coarse equivalence from $X$ to $Y$.
\end{defn}

\begin{prop}\label{checkubgen}
Let $f\colon X\to Y$ be a map and $\mathcal{M}$ a collection of subsets of $X\times X$.
If $\mathcal{E}$ is a coarse structure on $Y$ and $f(M)\in\mathcal{E}$ for all $M\in\mathcal{M}$,
then $f\colon (X,\cs(\mathcal{M})) \to (Y,\mathcal{E})$ is coarsely uniform.
\end{prop}
\textbf{Proof.}
The proposition follows easily from Remark~\ref{explcs} and the fact that
$f\!\times\! f(E_1\, E_2)\subseteq f\!\times\! f(E_1)\ f\!\times\! f(E_2)$
for all $E_1,E_2\subseteq X\times X$.
\qua

\begin{prop}
Let $(X,\mathcal{E}_X)$ and $(Y,\mathcal{E}_Y)$ be coarse spaces
and suppose that $f_1\colon X\to Y$ and $f_2\colon X\to Y$ are close maps.
If $f_1$ is coarsely proper, coarsely uniform or coarse, the same is true for $f_2$.
If $f_1$ is a coarse embedding or a coarse equivalence, then so is $f_2$.
\end{prop}
\textbf{Proof.}
Define $M := \{(f_1(x),f_2(x)) \mid x\in X\}$ and note that $M\in\mathcal{E}_Y$.

Observe that $f_2\times f_2(E)\subseteq M^{-1} (f_1\times f_1(E)) M$ for $E\in\mathcal{E}_X$.
This implies the statement about coarsely uniform maps.

For a bounded set $B\subseteq Y$ we notice that $f_2^{-1}(B)$ is contained in the bounded set $f_1^{-1}(M[B])$.
Hence, the claim about coarsely proper maps follows.

For $E\in\mathcal{E}_Y$ we have $(f_2\times f_2)^{-1}(E)\subseteq (f_1\times f_1)^{-1}(MEM^{-1})$.
This implies the claim about coarse embeddings.

Suppose that $f_1$ is a coarse equivalence and $g\colon Y\to X$ a coarse inverse of $f_1$.
Observe that $f_1\circ g$ and $f_2\circ g$ are close maps.
Since $g$ is coarsely uniform, we obtain that $g\circ f_1$ and $g\circ f_2$ are close maps.
Since ``being close'' has the same properties as an equivalence relation,
it follows that $g$ is also a coarse inverse of $f_2$.
\qua

\section[Some constructions]{Some constructions in the coarse category} \label{CGconstr}

\subsection*{Pull-back of a coarse structure}
Let$(X,\mathcal{E})$ be a coarse space, $A$ an arbitrary set and $f\colon A\to X$ any map.
\begin{defn}\index{pull-back}
We call $f^*(\mathcal{E}) := \cs\left(\left\{(f\times f)^{-1}(E) \mid E\in\mathcal{E}\right\}\right)$
the \emph{pull-back} of $\mathcal{E}$ by $f$.\footnote{
   Since the collection $\{(f\times f)^{-1}(E) \mid E\in\mathcal{E}\}$
   is not necessarily closed under taking subsets, we need to take the coarse structure generated by this collection.
   Note that this is the only reason for putting $\cs$ here.}
\end{defn}
Clearly, $f^*(\mathcal{E})$ is a coarse structure on $A$ such that $f$ is a coarse embedding.
If $(X,\mathcal{E})$ is connected, the same is true for $(A,f^*(\mathcal{E}))$.

Suppose $f$ is a continuous and proper map between topological spaces.
If $\mathcal{E}$ is compatible with the topology, the same is true for $f^*(\mathcal{E})$.
\begin{defn}\index{restriction}\textbf{(restriction of a coarse structure to a subset)} \\
If $A\subseteq X$ and $f$ is the inclusion map,
we call $f^*(\mathcal{E})$ the \emph{restriction} of $\mathcal{E}$ to $A$ and write $\mathcal{E}|_A$.
\end{defn}

\subsection*{Unions of coarse spaces}
Let $I$ be a set, $(X_i,\mathcal{E}_i)$ a coarse space for every $i\in I$ and $X=\bigcup_{i\in I} X_i$.
Note that the sets $X_i$ for $i\in I$ are not supposed to be disjoint.
\begin{defn}\index{union of coarse spaces}
We call
$\bigvee_{i\in I}\mathcal{E}_i
   := \cs\left(\bigcup_{i\in I}\mathcal{E}_i\right)$
the minimal coarse structure on the union of $(X_i,\mathcal{E}_i)$.
\end{defn}

\begin{rem}
If we are given a coarse structure $\mathcal{E}$ on $X$
such that for all $i\in I$ its restriction to $X_i$ is $\mathcal{E}_i$,
then $\bigvee_{i\in I}\mathcal{E}_i\subseteq\mathcal{E}$.
\end{rem}

\subsection*{Direct limits of coarse spaces}

For information about category theory compare \cite{MacL}.
Let $I$ be a small category and let $F\colon I\to\mathcal{D}$ or $F\colon I\to\mathcal{D}_{cn}$ be a functor.
In particular, we are given a coarse space $(X_i,\mathcal{E}_i)$ for each $i\in Obj(I)$
and a coarsely uniform map $F(h)\colon (X_i,\mathcal{E}_i) \to (X_j,\mathcal{E}_j)$
for each $h\in Mor(i,j)$.

\begin{prop}[existence of direct limits in $\mathcal{D}$ and $\mathcal{D}_{cn}$]
\label{colimex}\index{direct limit}
The categories $\mathcal{D}$ and $\mathcal{D}_{cn}$ are co-complete, i.e. $F$ has a colimit.
More precisely,
$$\underrightarrow{\lim}\ (X_i,\mathcal{E}_i) = \left( \underrightarrow{\lim} X_i\ ,\ \mathcal{E}_{\underrightarrow{\lim}} \right)$$
where $\underrightarrow{\lim} X_i$ is the corresponding colimit in the category of sets
and $\mathcal{E}_{\underrightarrow{\lim}}$ is the intersection of all (connected) coarse structures
such that the maps $f_i\colon X_i\to \underrightarrow{\lim} X_i$
are coarsely uniform for all $i\in Obj(I)$.
\end{prop}
\textbf{Proof.}
Observe that $\mathcal{E}_{\underrightarrow{\lim}} = \cs\left( \bigcup_{i\in Obj(I)} f_i(\mathcal{E}_i) \right)$.
With Proposition~\ref{checkubgen} the claim follows easily from the definitions.
\qua

\bigskip

The analogue of Proposition~\ref{colimex} in $\mathcal{C}$ and $\mathcal{C}_{cn}$
is not true. See Example~\ref{colomR}.
However, coproducts do exist in the categories $\mathcal{C}$ and $\mathcal{C}_{cn}$.

\begin{example} \label{colomR}
For each $i\in\NNN$ take $X_i:=\RRR_+$ with the usual bounded coarse strucure and define
$f_{ji}\colon X_i\to X_j,\ x\mapsto\max\{\ 0,\ x-j+i\ \}$ for $i\leq j$.
Then $\underrightarrow{\lim} X_i$ is just a point.
\end{example}

In some cases we can give a more explicit description for the coarse structure of the direct limit.

\begin{prop}\label{limpor}\index{directed set}
Let $I$ be a directed set\footnote{
   Compare \cite{Spanier} for the definition of a directed set
   and the definition of direct limits in this special case.
}.
If for all $i\in I$ the map $f_i\colon X_i\to \underrightarrow{\lim} X_i$ is injective, then
$$\mathcal{E}_{\underrightarrow{\lim}} = \{D \cup f_j\times f_j(L) \mid 
   j\in I, L\in\mathcal{E}_j, D\subseteq\Delta_{\underrightarrow{\lim} X_i} \}\ .$$
\end{prop}
\textbf{Proof.}
For $j\in I$, $L\in\mathcal{E}_j$ and $D\subseteq\Delta_{\underrightarrow{\lim} X_i}$
the set $D \cup f_j\times f_j(L)$ is in $\mathcal{E}_{\underrightarrow{\lim}}$.

It remains to prove that $\mathcal{E}_c :=
\{D \cup f_j\times f_j(L) \mid j\in I, L\in\mathcal{E}_j, D\subseteq\Delta_{\underrightarrow{\lim} X_i} \}$
is a coarse structure.
In order to conclude that $\mathcal{E}_c$ is closed under composition of elements,
we need the maps $f_i$ to be injective.
\qua

\begin{prop}\label{propcdcoprod}\index{coproduct}
Suppose that the identity morphisms are the only morphisms in the category $I$.
In this case the colimit is just the coproduct
and we denote the coarse structure of the direct limit by $\mathcal{E}_{\coprod}$.
We get
$\mathcal{E}_{\coprod} = \{D\cup E_{i_1}\cup\cdots\cup E_{i_k}
   \mid k\in\NNN, i_1,\ldots,i_k\in Obj(I), E_{i_j}\in\mathcal{E}_{i_j}, D\subseteq\Delta_{\coprod X_i} \}$.
\end{prop}

\subsection*{Infinite disjoint unions of coarse structures}

\begin{example}
$\RRR$ and $\ZZZ$ (both equipped with their usual bounded coarse strucure) are coarsely equivalent,
but the corresponding countable coproducts $\coprod_\NNN \RRR$ and $\coprod_\NNN \ZZZ$ are not.
\end{example}
\textbf{Proof.}
Any map $\coprod_\NNN \RRR \to \coprod_\NNN \ZZZ \to \coprod_\NNN \RRR$
will be different from the identity on any copy of $\RRR$.
Thus, no such map can be close to $id_{\coprod_\NNN \RRR}$,
since any entourage of the coproduct involves only off-diagonal elements
of finitely many copies of $\RRR$.
\qua

\bigskip

Because of the above example one might look for another more geometric coarse structure on the coproduct.

\begin{defn}\index{disjoint union}\textbf{(disjoint union of coarse spaces)} \\
Let $I$ be a set and assume for any $i\in I$ we are given a coarse space $(X_i,\mathcal{E}_i)$.
Set $X := \bigsqcup_{i\in I} X_i$. We define
$$\mathcal{E}_{\bigsqcup} := \left\{ A\subseteq \bigsqcup_{i\in I} X_i^2 \mid A\cap X_i^2\in\mathcal{E}_i \right\}.$$
\end{defn}
The underlying sets of coproduct and disjoint union are the same.
We just take different coarse structures.
For any finite set $I$, coproduct and disjoint union coincide.

\begin{prop}
If for every $i\in I$ the coarse spaces $(X_i,\mathcal{E}_i)$ and
$(Y_i,\mathcal{F}_i)$ are coarsely equivalent, then
$\bigsqcup_{i\in I} (X_i,\mathcal{E}_i)$ and
$\bigsqcup_{i\in I} (Y_i,\mathcal{F}_i)$ are also coarsely equivalent.
\end{prop}
\textbf{Proof.}
If a coarse equivalence $f_i\colon(X_i,\mathcal{E}_i)\to(Y_i,\mathcal{F}_i)$
is given for all $i\in I$, then $\bigsqcup_{i\in I}f_i\colon
\bigsqcup_{i\in I}(X_i,\mathcal{E}_i) \to \bigsqcup_{i\in I}(Y_i,\mathcal{F}_i)$
is also a coarse equivalence. 
\qua

\subsection*{Products of coarse spaces}

\begin{defn}\label{prodcs}\index{product coarse structure}
For $i\in\{1,\ldots,k\}$ let $(X_i,\mathcal{E}_i)$ be a coarse space.
By $p_i\colon X_1\times\cdots\times X_k \to X_i$ we denote the projection to the $i$-th factor.
The \emph{product coarse structure} is defined as follows.
$$\mathcal{E}_1 * \cdots * \mathcal{E}_k := \left\{E\subseteq (X_1\times\cdots\times X_k)^2 \mid
          p_i\times p_i(E)\in\mathcal{E}_i\text{ for }i\in\{1,\ldots,k\}\right\}$$
If $(X,\mathcal{E})$ is a coarse space, we will sometimes write $\mathcal{E}^{*k}$
for the product coarse structure on $X^k$.
\end{defn}

It is easy to prove that $\mathcal{E}_1 * \cdots * \mathcal{E}_k$ actually is a coarse structure.
The product coarse structure $\mathcal{E}_1 * \cdots * \mathcal{E}_k$ is connected
if and only if the coarse structures $\mathcal{E}_i$ are connected.
Moreover, we have the following formulas:
\begin{eqnarray*}
\Delta_{X_1\times\cdots\times X_k} &=& \Delta_{X_1} \times\cdots\times \Delta_{X_k} \\
(E_1 \times\cdots\times E_k)(E'_1 \times\cdots\times E'_k) &=& E_1E'_1 \times\cdots\times E_kE'_k
\end{eqnarray*}

One remark on our notation:
If $E\subseteq X\times X$, we should not confuse the composition $E^k=E\cdots E$
and the product $E^{\times k} = E\times\cdots\times E$.

\begin{rem}\textbf{(compatibility with other products)}
\begin{itemize}
\item 
Let $(X,d_X)$ and $(Y,d_Y)$ be metric spaces.
Then the formula $$d\left((x_1,y_1),(x_2,y_2)\right):=d_X(x_1,x_2)+d_Y(y_1,y_2)$$ defines a metric on $X\times Y$.
It is easy to see that $\mathcal{E}_d = \mathcal{E}_{d_X} * \mathcal{E}_{d_Y}$.
\item
Let $X$ and $Y$ be Hausdorff spaces with compactifications $\overline{X}$ and $\overline{Y}$ respectively.
Then $\overline{X}\times\overline{Y}$ is a compactification of $X\times Y$.
By $\mathcal{E}_{\overline{X}}$, $\mathcal{E}_{\overline{Y}}$ and $\mathcal{E}_{\overline{X}\times\overline{Y}}$
we denote the coarse structures on $X$, $Y$ and $X\times Y$ induced by $\overline{X}$, $\overline{Y}$ and
$\overline{X}\times\overline{Y}$ respectively.
The definitions easily imply
$\mathcal{E}_{\overline{X}\times\overline{Y}}=\mathcal{E}_{\overline{X}} * \mathcal{E}_{\overline{Y}}$.
\end{itemize}
\end{rem}

\begin{rem}\label{prodce}\textbf{(compatibility with coarse equivalence)}\\ 
Let $(X,\mathcal{E})$, $(X',\mathcal{E}')$ and $(Y,\mathcal{F})$ be coarse spaces
and assume $(X,\mathcal{E})$ and $(X',\mathcal{E}')$ to be coarsely equivalent.
Then the spaces $(X\times Y,\mathcal{E} * \mathcal{F})$ and $(X'\times Y,\mathcal{E}' * \mathcal{F})$
are coarsely equivalent.
\end{rem}

\subsection*{Quotients of coarse spaces}

\begin{defn}\index{quotient coarse structure}
Let $(X,\mathcal{E})$ be a coarse space. Let $\sim$ be an equivalence relation on $X$
and $\ \pi\colon X\to X/\!\!\sim\ $ the corresponding quotient map.
We define a coarse structure on $X/\!\!\sim$ by
$$\mathcal{E}_\sim := \cs\left( \{ \pi\times\pi(E) \mid E\in\mathcal{E} \} \right)\ .$$
If $\mathcal{E}$ is connected, $\mathcal{E}_\sim$ is connected.
The map $\pi$ is coarsely uniform by definition of the coarse structure on the quotient.
If there is an unbounded equivalence class, then $\pi$ is not coarsely proper.
\end{defn}

\subsection*{Attaching coarse spaces}\index{attaching coarse spaces}

Let $(X,\mathcal{E}_X)$, $(Y,\mathcal{E}_Y)$ be coarse spaces,
$A\subseteq X$ and $f\colon A\to Y$ coarsely uniform.
Attaching $X$ to $Y$ using $f$, we get $X\cup_f Y := X\sqcup Y\, /\sim\ $
and a quotient map $\pi\colon X\sqcup Y \to X\cup_f Y$.
Here $\ \sim\ $ is the equivalence relation generated by $a \sim f(a)$ for all $a\in A$.
The map $\pi$ is coarsely uniform, but not necessarily coarsely proper.
We denote the coarse structure on the quotient by $\mathcal{E}_{X\cup_f Y}$.

\begin{lemma} \label{csonY}
$\mathcal{E}_Y = \mathcal{E}_{X\cup_f Y}|_Y$
\end{lemma}
\textbf{Proof.}
We certainly have $\mathcal{E}_Y \subseteq \mathcal{E}_{X\cup_f Y}|_Y$ and need to prove the converse.
Let $L\in\mathcal{E}_{X\cup_f Y}|_Y$.
There are $L'_{X,1},\ldots,L'_{X,k}\in\mathcal{E}_X$ and $L_{Y,1},\ldots,L_{Y,k}\in\mathcal{E}_Y$,
such that $$L\subseteq\pi\times\pi(L'_{X,1}\cup L_{Y,1})\ \cdots\ \pi\times\pi(L'_{X,k}\cup L_{Y,k}).$$

Define $L_{X,i}\cdots L_{X,i+j}:=f\times f(L'_{X,i}\cdots L'_{X,i+j})\in\mathcal{E}_Y$.
This is a partial replacement of the notation of composition which we will apply only in this proof.
We conclude that
$$L\ \subseteq\bigcup_{(\sigma_1,\ldots,\sigma_k)\in\{X,Y\}^n} L_{\sigma_1,1}\cdots L_{\sigma_k,k}\in\mathcal{E}_Y\ .$$
\qua

\bigskip

Since $\pi|_{X\backslash A}$ is coarsely uniform,
we have $\mathcal{E}_X|_{X\backslash A}\subseteq\mathcal{E}_{X\cup_f Y}|_{X\backslash A}$.
But in general these coarse structures do not coincide as the following example shows.
\begin{example}
Let $Y$ be a single point, $X'$ any coarse space such that $X'\times X'$ is not an entourage,
$X=X'\times \RRR_+$ and $A=X'\times\{0\}$. For $r>0$ the set
$\{((x,t),(\widetilde{x},\widetilde{t})) \mid t,\widetilde{t}\in\left(0,r\right], x,\widetilde{x}\in X'\}$
is contained in $\mathcal{E}_{X\cup_f Y}|_{X\backslash A}$ but not in $\mathcal{E}_X|_{X\backslash A}$.
\end{example}

\begin{rem}
If the attaching map $f$ is a coarse embedding, we have
$\mathcal{E}_X|_{X\backslash A} = \mathcal{E}_{X\cup_f Y}|_{X\backslash A}$.
\end{rem}
\textbf{Proof.}
The map $\pi|_X$ is a coarse equivalence and $\pi|_{X\backslash A}$ is a canonical isomorphism.
\qua

\section{Coarse structures induced by metrisable compactifications}\label{secHigson}

We will briefly recall the definitions of the Higson compactification and the Higson corona.
For more detailed explanations see \cite{Roe}, \cite{HR} or \cite{RoeCG}.

Let $X$ be a Hausdorff space with a coarse structure $\mathcal{E}$
such that the closure of every bounded set is compact.
A bounded continuous function $f$ on $X$ is called a \emph{Higson function}\index{Higson function}
with respect to a coarse structure $\mathcal{E}$ if for all entourages $E\in\mathcal{E}$ the function
$$\mathbf{d}f\colon X\times X\to \CCC , \quad (x_1,x_2)\mapsto f(x_1)-f(x_2)$$
restricted to $E$ tends to zero at infinity.
The Higson functions form a unital $C^*$-algebra $C_h(X)$.
It follows from the Gelfand-Naimark theorem
that $C_h(X)$ is the algebra of bounded, continuous functions on a compactification $hX$ of $X$.
We call $hX$ the \emph{Higson compactification} of $(X,\mathcal{E})$.\index{Higson compactification}
The \emph{Higson corona} $\nu X$ is defined to be $hX\backslash X$. \index{Higson corona}

For metrisable compactifications the relation between compactifications and coarse structures
is especially simple:
\begin{prop}
Let $X$ be a Hausdorff space and $K$ a metrisable compactification. The Higson compactification
of the continuously controlled coarse structure on $X$ induced by $K$ is homeomorphic to $K$.
\end{prop}
\textbf{Proof.} This is Proposition 2.48 of \cite{RoeCG}. \qua

\bigskip

Let $X$ be a non-compact Hausdorff space and $hX$ a metrisable compactification of $X$.
By $\mathcal{E}_{hX}$ denote the continuously controlled coarse structure induced by $hX$.

Choose a metric $d$ on $hX$ and define
$$X_i:=\left\{x\in X\mid d(x,\nu X)\geq \frac{1}{i}\right\}$$
for $i\in\NNN$ ($X_0:=\emptyset$).
We have defined a sequence of compact subsets of $X$ such that
$X_0 \subseteq X_1 \subseteq X_2 \subseteq \cdots \subseteq X$
and $\bigcup_{i\in\NNN}X_i=X$.
Moreover $\{hX\backslash X_i\}_{i\in\NNN}$ is a basis of neighborhoods for
$\nu X:=hX\backslash X$.

\begin{lemma}\label{descriptEhX}
For $\rho >0$ define $\Delta_\rho:=\{(x,y)\in X^2\mid d(x,y)<\rho\}$.
We assert that $E\in\mathcal{E}_{hX}$ if and only if there is a sequence $\{\rho_i\}_{i\in\NNN}$
of non-negative real numbers converging to zero such that
$E\backslash X_i^2\subseteq \Delta_{\rho_i}$ for all $i\in\NNN$.
\end{lemma}
The proof of this lemma is straightforward. \qua

\bigskip

Our intention is to construct a coarse structure $\mathcal{E}$ on $\nu X\times \NNN$
such that $(X,\mathcal{E}_{hX})$ is coarsely equivalent to $(\nu X\times \NNN,\mathcal{E})$.

We may assume that $X_1\neq\emptyset$.
Choose a map
\begin{equation}\label{abbf}
f\colon \nu X\times \NNN\to X
\end{equation}
such that $f(\overline{x},n)\in X_n$
has minimal distance to $\overline{x}$ among all points in $X_n$.
Since $X_n$ is compact, this definition makes sense.
Next, we define a coarse structure on $\nu X\times \NNN$ by pulling back the coarse structure $\mathcal{E}_{hX}$. Thus $f$ becomes a coarse embedding.
$$\mathcal{E}:=f^*(\mathcal{E}_{hX})$$
By $\pi_\NNN\colon\nu X\times\NNN\to\NNN$ we denote the projection onto $\NNN$
and by $\mathcal{E}_\NNN$ the continuously controlled coarse structure on $\NNN$
induced by the one-point compactification. Observe that $\mathcal{E}_\NNN$
is the maximal coarse structure which is compatible with the topology.
\begin{lemma}\label{descriptE}
We assert that $E\in\mathcal{E}$ if and only if the following two conditions are satisfied.
\begin{enumerate}
\item[(a)] $\pi_\NNN\times\pi_\NNN(E)\in\mathcal{E}_\NNN$
\item[(b)] There is a sequence $\{\delta_k\}_{k\in\NNN}$ of non-negative real numbers converging to zero such that
$d(\overline{x},\overline{y})<\delta_k$ for all $((\overline{x},i),(\overline{y},j))\in E$ with $\max\{i,j\}>k$.
\end{enumerate}
\end{lemma}
\textbf{Proof.}
Define $$a_i:=\max\{d(X_i,\overline{x})\mid\overline{x}\in\nu X\}$$ for $i\in\NNN$.
The sequence $\{a_i\}_{i\in\NNN}$ is monotone decreasing, since $X_i\subseteq X_{i+1}$ for all $i\in\NNN$.
We claim that this sequence converges to zero.

Assume $a_i>a>0$ for all $i\in\NNN$.
Take the centers of a finite cover of $\nu X$ with balls of radius $\frac{a}{3}$.
For each of these points there is a point in $X$ within distance less than $\frac{a}{3}$.
These finitely many points in $X$ form a compact set and hence are contained in some $X_k$.
We get $a_k\leq\frac{2}{3}a$, which contradicts our assumption.
Hence $\underset{i\to\infty}{\lim} a_i = 0$.

Define $$n_i:=\max\left\{n\in\NNN\mid a_i<\frac{1}{n}\right\}$$ for all $i$ with $a_i<1$.
The equation $\underset{i\to\infty}{\lim}a_i = 0$ implies
$\underset{i\to\infty}{\lim} n_i=\infty$.

For $M\subseteq \NNN\times\NNN$ we define a completion $M^{\compl}$ of $M$ in the following way. Set
$M^{\compl} := \{ (u,v)\in\NNN^2 \mid (x,y)\in M\text{ and } (x\leq u\leq v\leq y \text{ or } y\leq v\leq u\leq x) \}$.
Let $E\subseteq (\nu X\times \NNN)^2$ and define 
$$b_i:=\min\left\{n\in\NNN \mid (n,i)\in\pi_\NNN\!\times\!\pi_\NNN\left(E\cup E^{-1}\right)^{\compl}\right\}.$$

Suppose that $E$ satisfies conditions (a) and (b).
First we will prove that $f\times f(E)\in\mathcal{E}_{hX}$.
Since $\mathcal{E}_\NNN$ is compatible with the topology,
condition (a) implies $\underset{i\to\infty}{\lim}b_i=\infty$.
We define a sequence $\{\rho_i\}_{i\in\NNN}$.
$$ \rho_i:=a_i+a_{b_i}+\delta_i$$
Certainly we have $\underset{i\to\infty}{\lim}\rho_i=0$.

Let $((\overline{x},n),(\overline{y},m))\in E$ and define $x:=f(\overline{x},n)$ and $y:=f(\overline{y},m)$.
If $x\not\in X_i$, then $d(x,\overline{x})\leq a_i$, $d(y,\overline{y})\leq a_{b_i}$
and $d(\overline{x},\overline{y})\leq\delta_i$.
Hence $f\times f(E)\backslash X_i^2\subseteq \Delta_{\rho_i}$.
Now we can apply Lemma~\ref{descriptEhX} to conclude $f\times f(E)\in\mathcal{E}_{hX}$.
It follows $E\subseteq (f\times f)^{-1}(f\times f(E))\in\mathcal{E}$.

Suppose $E\in\mathcal{E}$. There is an entourage $E'\in\mathcal{E}_{hX}$ with $E\subseteq (f\times f)^{-1}(E')$.
For all $i\in\NNN$ there is $m_i\in\NNN$ such that $a_{m_i}<\frac{1}{i}$.
Since $\frac{1}{m_i}\leq a_{m_i}$, we get $m_i>i$ for all $i\in\NNN$.
We may assume that the sequence $\{m_i\}_{i\in\NNN}$ is monotone.

Since $X_i$ is bounded, $E'[X_i]$ is bounded, i.e. $E'[X_i]$ is contained in some $X_j$.
Note that the image of $\nu X\times\{0,\ldots,i\}$ under $f$ is contained in $X_i$.
Observe that $E[\nu X\times\{0,\ldots,i\}] \subseteq f^{-1}\left(E'[X_i]\right) \subseteq f^{-1}\left(X_j\right)$.
It follows that $\pi_\NNN(E[\nu X\times\{0,\ldots,i\}])\subseteq\pi_\NNN(f^{-1}(X_j))\subseteq\{0,\ldots,m_j-1\}$.
Together with an analog fact for $E^{-1}$ this implies that $E$ satisfies condition (a).

Choose a sequence $\{\rho_i\}_{i\in\NNN}$ converging to zero
and such that $E'\backslash X_i^2\subseteq\Delta_{\rho_i}$ for all $i\in\NNN$.
Let $((\overline{x},i),(\overline{y},j))\in E$ with $\max\{i,j\}>k$.
Assume $i\geq k$. It follows that $j\geq b_k$.
We already proved $\pi_\NNN\times\pi_\NNN(E)\in\mathcal{E}_\NNN$. Hence $\underset{i\to\infty}{\lim}b_i=\infty$.
Define $x:=f(\overline{x},i)$ and $y:=f(\overline{y},j)$.
We have $d(x,\overline{x})\leq a_i \leq a_k < \frac{1}{n_k}$, i.e. $x\not\in X_{n_k}$.
Therefore $d(x,y) < \rho_{n_k}$.
Note that $$\delta_k:=a_k+a_{b_k}+\rho_{n_k} \underset{k\to\infty}{ \ -\!-\!\!\!\longrightarrow} 0\ .$$
Obviously $E$ satisfies (b) with $\{\delta_k\}_{k\in\NNN}$ as defined above.
\qua

\bigskip

We take a map
\begin{equation}\label{abbg}
g\colon X\to\nu X\times\NNN
\end{equation}
such that for $x\in X_i\backslash X_{i-1}$
the image  $g(x)=(\overline{x},i)$ is chosen in such a way
that $\overline{x}$ minimizes the distance to $x$ among all points in $\nu X$.

\begin{lemma}
The map $g$ is coarsely uniform.
\end{lemma}
\textbf{Proof.}
Let $E\in\mathcal{E}_{hX}$.
Since $\mathcal{E}_{hX}$ is compatible with the topology,
$g\times g(E)$ satisfies condition (a) of Lemma~\ref{descriptE}.
Let $\{\rho_i\}_{i\in\NNN}$ be a sequence as in Lemma~\ref{descriptEhX} for the entourage $E$
and define $b_i:=\min\{n\in\NNN\mid E[X_i\backslash X_{i-1}]\cap X_n\neq\emptyset\}$.
This implies $\underset{i\to\infty}{\lim}b_i=\infty$.
Defining $\delta_i:=\frac{1}{i-1}+\frac{1}{b_i-1}+\rho_{i-1}$
we see that $g\times g(E)$ satisfies condition (b) of Lemma~\ref{descriptE}.
\qua

\begin{thm}\label{thmXnuXN}
Any map $f$ chosen as in (\ref{abbf}) is a coarse equivalence
and any map $g$ chosen as in (\ref{abbg}) is a coarse inverse of $f$.
\end{thm}
\textbf{Proof.}
We will use some of the notation from the proof of Lemma~\ref{descriptE}.
We first consider the map $g\circ f\colon \nu X\times\NNN\to\nu X\times\NNN$.
Let $(x,k)\in\nu X\times \NNN$.
Observe that $d(f(x,k),x)\leq d(X_k,x)\leq a_k < \frac{1}{n_k}$.
Hence $f(x,k)\in X_k\backslash X_{n_k}$.
Defining $(\widetilde{x},\widetilde{k}):= g\circ f(x,k)$ we get
$n_k+1\leq\widetilde{k}\leq k$. Moreover $d(x,\widetilde{x})\leq 2a_k$.
It follows that conditions (a) and (b) of Lemma~\ref{descriptE} are satisfied for
$\{((x,k),g\circ f(x,k))\mid x\in \nu X\times\NNN\}$.

Now we consider the map $f\circ g\colon X\to X$.
Let $x\in X_i\backslash X_{i-1}$ and $(\overline{x},i):=g(x)$.
Then $d(x,\overline{x})<\frac{1}{i-1}$.
Thus $d(f\circ g(x),x)\leq\frac{2}{i-1}$, i.e. $f\circ g$ is close to $\operatorname{id}_X$.
\qua

\begin{cor}\label{ceiffnuh}
Let $X$ and $Y$ be Hausdorff spaces with metrisable compactifications $hX$ and $hY$ respectively.
The coarse spaces $(X,\mathcal{E}_{hX})$ and $(Y,\mathcal{E}_{hY})$ are coarsely equivalent
if and only if the coronas $\nu X$ and $\nu Y$ are homeomorphic.
\end{cor}
\textbf{Proof.}
To see that the Higson coronas of coarsely equivalent spaces are homeomorphic,
we refer to Corollary 2.42 of \cite{RoeCG}.
\qua

\bigskip

To every compact Hausdorff space $K$ we assign a coarse space $\psi(K)$ as follows.
$$\psi(K):= \left(K\times\left[0,1\right), \mathcal{E}_{K\times[0,1]}\right)$$

\begin{cor}\label{coraboutPsi}
$\psi$ induces a bijection $\Psi$ between homeomorphism classes of metrisable compact spaces\footnote{
Since every metrisable compact set is homeomorphic to a subset of $[0,1]^\NNN$,
the homeomorphism classes of metrisable compact spaces form a set.}
and coarse equivalence classes of spaces whose coarse structure is induced by a metrisable compactification.
The inverse of $\Psi$ is induced by the correspondence $\nu$ assigning to each coarse space its Higson corona.
\end{cor}

We define a "product" of two compact spaces $A$ and $B$ as follows.
\begin{gather*}
A \boxtimes B\, :=\, A\!\times\! B\!\times\! [0,1] \, \cong\,
A\!\times\!\{1\}\!\times\! B\!\times\! [0,1]\, \cup\, A\!\times\! [0,1]\!\times\! B\!\times\!\{1\} \\
= A\!\times\! [0,1]\!\times\! B\!\times\! [0,1] \ \backslash \ 
A\!\times\! \left[0,1\right)\!\times\! B\!\times\! \left[0,1\right)
\end{gather*}
Here $\cong$ stands for homeomorphism.

\begin{prop}\label{prodPsi}
If $A$ and $B$ are metrisable compact spaces, then $$\Psi(A \boxtimes B)=\Psi(A)\times\Psi(B).$$
Here $\times$ denotes the product of coarse spaces as in Definition~\ref{prodcs}.
\end{prop}
\textbf{Proof.}
Because of Remark~\ref{prodce}, the right hand side of the equation makes sense.
The Higson corona of
$\psi(A)\times\psi(B)=A\!\times\! \left[0,1\right)\!\times\! B\!\times\! \left[0,1\right)$
is homeomorphic to $A\boxtimes B$.
Thus, Corollary~\ref{coraboutPsi} implies $\Psi(A \boxtimes B)=\Psi(A)\times\Psi(B)$.
\qua

\begin{example}\label{visualdot}
By $\mathcal{E}_\cdot$ we denote the coarse structure on $\RRR_+$ induced by the one-point compactification and by $\mathcal{E}_{vis}$ the coarse structure induced on $\RRR^n$
by the visual corona $S^{n-1}$.
The usual bounded coarse structure on $\RRR^n$ is contained in $\mathcal{E}_{vis}$.
Proposition~\ref{prodPsi} implies that
$(\RRR_+,\mathcal{E}_\cdot)^n$ and $(\RRR_+^n,\mathcal{E}_{vis}|_{\RRR_+^n})$
are coarsely equivalent.
\end{example}

\clearpage
\thispagestyle{empty}\ \newpage 

\chapter{Translation $C^*$-algebras via projection-valued measures}

To any coarse space $(X,\mathcal{E})$ which admits a countable, uniformly bounded cover,
we will assign translation $C^*$-algebras $E^*(X,\mathcal{E})$ and $C^*(X,\mathcal{E})$.

This may be considered as an analog of assigning to a locally compact space $X$
the algebra $C_0(X)$ of complex valued, continuous function vanishing at infinity.

In \cite{HR}, a representation of the algebra $C_0(X)$ is used to define translation $C^*$-algebras
for coarse spaces which are at the same time nice topological spaces such that topology and coarse structure are compatible.
John Roe proposed in \cite{RoeCG} to use projection-valued measures
in order to define translation $C^*$-algebras for more general coarse spaces.
In this chapter we give a more complete exposition on the subject and eliminate some mistakes.

The idea is to generalize an ad-hoc approach
which has been used for defining translation $C^*$-algebras of discrete coarse spaces.

We prove that the translation $C^*$-algebras are well defined up to isomorphism of $C^*$-algebras,
we compare the different approaches of defining translation $C^*$-algebras and we discuss the matter
whether a coarse map induces maps between the translation $C^*$-algebras of the corresponding coarse spaces.

\section{Geometric Hilbert spaces and the calculus of supports}

Let $(X,\mathcal{E})$ be a coarse space.
\begin{defn}\index{localizing decomposition}
A family $\mathcal{U}$ of subsets of $X$ is called a \emph{localizing decomposition} of $(X,\mathcal{E})$,
if (1) $\mathcal{U}$ is a decomposition, i.e. a cover of $X$ consisting of pairwise disjoint sets,
(2) $\mathcal{U}$ is countable and (3) $\mathcal{U}$ is uniformly bounded,
i.e. $\Delta_\mathcal{U} := \bigcup_{U\in\mathcal{U}} {U\times U} \in\mathcal{E}$.
\end{defn}

\begin{rem}
Let $(X,\mathcal{E})$ be a coarse space and $\mathcal{U}$ a localizing decomposition.
We say that $x,y\in X$ are $\mathcal{U}$-equivalent if they are contained in the same $U\in\mathcal{U}$.
We denote the quotient of $X$ with respect to $\mathcal{U}$-equivalence by $X/_\mathcal{U}$\,.
Using the quotient map $\pi\colon X\to X/_\mathcal{U}$, we define the coarse structure
$\mathcal{E}/_\mathcal{U}:=\{\, \pi\!\times\!\pi(E) \mid E\in\mathcal{E} \}$ on $X/_\mathcal{U}$\,.

In this sence, a localizing decomposition $\mathcal{U}$ induces
a countable ``discretization''\index{discretization} $(X/_\mathcal{U},\mathcal{E}/_\mathcal{U})$
which is coarsely equivalent to $(X,\mathcal{E})$.
\end{rem}

Let $\mathfrak{A}$ be a $\sigma$-algebra over $X$.\footnote{
A $\sigma$-algebra $\mathfrak{A}$ over a set $X$ is a family of subsets of $X$
such that (1) $X\in\mathfrak{A}$,
(2) $X\backslash A\in\mathfrak{A}$ whenever $A\in\mathfrak{A}$
and (3) the union of countably many elements of $\mathfrak{A}$ is in $\mathfrak{A}$.
Compare any book on measure theory.}

\begin{defn}\label{locsigalg}\index{localizing $\sigma$-algebra}
A localizing $\sigma$-algebra over $(X,\mathcal{E})$
is a $\sigma$-algebra $\mathfrak{A}$ together with a topology on $X$ such that
for each entourage $E\in\mathcal{E}$ the closure $\overline{E}$ is still an entourage.
Moreover, there has to exist a localizing decomposition $\mathcal{U}\subseteq\mathfrak{A}$
such that $\bigcup_{U\in\mathcal{U}}\kringel{U}$ is dense in $X$.
We will call $\mathcal{U}$ an $\mathfrak{A}$-localizing decomposition.

For each $x\in X$ we define $\mathfrak{A}(x)\subseteq\mathfrak{A}$ to be all the sets $A\in\mathfrak{A}$
whose interior $\kringel{A}$ contains $x$.
Such an $A\in\mathfrak{A}(x)$ will be called an $\mathfrak{A}$-neighborhood of $x\in X$.
\end{defn}

\begin{example}\label{locsigalgex}\index{localizing $\sigma$-algebra}
Given a coarse space $(X,\mathcal{E})$
which admits a uniformly boun\-ded cover $\mathcal{V}=\{V_n\mid n\in\NNN\}$,
there are several localizing $\sigma$-algebras:
\begin{itemize}
\item[(1)] The $\sigma$-algebra $\sigma(\mathcal{V})$ generated by $\mathcal{V}$
   together with the discrete topology.
\item[(2)] The $\sigma$-algebra $\mathfrak{P}(X)$ consisting of all subsets of $X$
   together with the discrete topology.
\item[(3)] If $X$ is a topological space and the coarse structure $\mathcal{E}$
   is compatible with the topology, then the Borel algebra $\mathfrak{B}(X)$,
   i.e. the smallest $\sigma$-algebra containing all open sets, is localizing.
\end{itemize}
The $\sigma$-algebra $\mathfrak{A}=\{\emptyset,X\}$ is not localizing if $X$ is unbounded.
\end{example}
\textbf{Proof.}
Defining $U_n=V_n\backslash (V_0\cup\cdots\cup V_{n-1})$
we get a localizing decomposition $\mathcal{U}=\{U_n\mid n\in\NNN\}$
which is contained in the $\sigma$-algebras $\sigma(\mathcal{V})$ and $\mathfrak{P}(X)$.

In order to prove (3), note that from $\mathcal{V}$
we can construct a countable, uniformly bounded cover consisting of open sets as follows: 
If $E$ is an open entourage containing $\Delta_X$, take the cover $\{E[V_n]\mid n\in\NNN\}$.
From this cover we get a localizing decomposition as before.
\qua

\begin{defn} \label{defPVM}\index{projection-valued measure}
Let $\mathfrak{A}$ be a $\sigma$-algebra over $X$ and $H$ a Hilbert space.
A \emph{projection-valued measure} is a map
$\lambda\colon\mathfrak{A}\to\mathcal{B}(H)$ with the following properties:
\begin{enumerate}
\item[(1)]For every $A\in\mathfrak{A}$ the operator $\lambda(A)$ is a projection, \\
   i.e. $\lambda(A)^*=\lambda(A)=\lambda(A)^2$.
\item[(2)] $\lambda(\emptyset)=0$ and $\lambda(X)=\id_H$.
\item[(3)] If $A_n\in\mathfrak{A}$ for all $n\in\NNN$ and $A_i\cap A_j=\emptyset$ whenever $i\neq j$, then $$\lambda\left(\bigcup_{n\in\NNN}A_n\right)=\sum_{n\in\NNN}\lambda(A_n).$$
\end{enumerate}
\end{defn}

\begin{rem}
Note that the infinite sum in Definition~\ref{defPVM} has to be considered as a limit
with respect to strong convergence.\footnote{
A sequence $\{T_i\}_{i\in\NNN}$ in $\mathcal{B}(H)$ converges strongly to $T\in\mathcal{B}(H)$,
if for any $v\in H$ the sequence $\{T_i(v)\}_{i\in\NNN}$ converges to $T(v)$ in $H$. }

The sequence $\left\{\sum_{n=0}^k\lambda(A_n)\right\}_{k\in\NNN}$ does in general not converge in the operator norm.
\end{rem}

\begin{lemma} \label{proj_commute}
Let $\lambda\colon\mathfrak{A}\to\mathcal{B}(H)$ be a projection-valued measure over $X$.
For $A,B\in\mathfrak{A}$ we have $\lambda(A)\lambda(B) = \lambda(B)\lambda(A) = \lambda(A\cap B)$.
\end{lemma}
\textbf{Proof.}
First suppose $A\cap B = \emptyset$.
In order to prove $\lambda(A)\lambda(B) = 0$,
observe that the images of $\lambda(A)$ and $\lambda(B)$ are orthogonal,
because $\lambda(A) + \lambda(B)=\lambda(A\cup B)$ is a projection.

The general case now follows easily.
\begin{eqnarray*}
\lambda(A)\lambda(B)
 & = & \big(\lambda(A\backslash B)+\lambda(A\cap B)\big)\big(\lambda(B\backslash A)+\lambda(B\cap A)\big) \\
 & = & \lambda(A\cap B)\lambda(B\cap A) \ = \ \lambda(A\cap B)
\end{eqnarray*}
The first equality uses additivity of projection-valued measures.
The last one requires that a projection $P$ is idempotent, i.e. $P^2=P$.
\qua

\begin{example} \label{measurePVM}
Take the Hilbert space $H=L^2(X,\mu)$
where $\mu$ is a measure on some localizing $\sigma$-algebra $\mathfrak{A}$.
A projection-valued measure $\lambda\colon \mathfrak{A}\to\mathcal{B}(H)$ is given by
$$\lambda(A)f = \chi_A\cdot f \qquad\text{ for } A\in\mathfrak{A} \text{ and } f\in H.$$
Here $\chi_A\colon X\to\CCC$ is the characteristic function of $A$, i.e.
$\chi_A(x)=1$ if $x\in A$ and $\chi(x)=0$ if $x\in X\backslash A$.
\end{example}

\begin{example} \label{decompositionPVM}
Let $(X,\mathcal{E})$ be a coarse space and
$\mathcal{U}=\{U_i\mid i\in\NNN\}$ a localizing decomposition of $X$.
Define $\mathfrak{A}:=\sigma(\mathcal{U})$ to be the $\sigma$-algebra generated by $\mathcal{U}$.
Assume we are given a Hilbert space $H_i$ for each $i\in\NNN$
such that $H_i = \{0\}$ if $U_i=\emptyset$.
In this situation we can define a projection-valued measure
$\lambda\colon\mathfrak{A}\to\mathfrak{B}\left(\bigoplus_{i\in\NNN}H_i\right)$ as follows:
For $A\in\mathfrak{A}$ define
$$\lambda(A) = \text{ orthogonal projection onto } \bigoplus_{\substack{i\in\NNN \\ U_i\subseteq A}}H_i \ .$$
\end{example}

\begin{defn}\index{geometric Hilbert space}
Let $(X,\mathcal{E})$ be a coarse space.
A \emph{geometric Hilbert space} over $(X,\mathcal{E})$ 
is a projection-valued measure $\lambda\colon\mathfrak{A}\to\mathcal{B}(H)$
where $\mathfrak{A}$ is a localizing $\sigma$-algebra over $(X,\mathcal{E})$
and $H$ a separable Hilbert space.
\end{defn}

\begin{rem}
The projection-valued measure in Example~\ref{measurePVM} is a geometric Hilbert space
if and only if $L^2(X,\mu)$ is separable.

The projection-valued measure in Example~\ref{decompositionPVM} is a geometric Hilbert space
if and only if the Hilbert space $H_i$ is separable for all $i\in\NNN$.
\end{rem}

\begin{defn} \label{defSupport}\index{support}
Let $\lambda\colon\mathfrak{A}\to\mathcal{B}(H)$ be a geometric Hilbert space over $(X,\mathcal{E})$.
We define the \emph{support} of $u\in H$ and $T\in\mathcal{B}(H)$ as follows:
\begin{eqnarray*}
\supp(u) & = & \big\{ x\in X \ \mid \ \lambda(A)u\neq 0
               \text{ for all } A\in\mathfrak{A}(x) \big\} \\
\supp(T) & = & \big\{ (x_1,x_2)\in X\!\times\! X \ \mid \ \lambda(A_1)T\lambda(A_2)\neq 0
               \text{ for all } A_i\in\mathfrak{A}(x_i) \big\}
\end{eqnarray*}
\end{defn}

\begin{rem}\label{badsupp}
Assume $\lambda(A)=0$ and $x\in\kringel{A}$.
Then $x\not\in\supp(u)$ for all $u\in H$.
Similarly $(x,y),(y,x)\not\in\supp(T)$ for $T\in\mathcal{B}(H)$ and $y\in X$.

In particular, consider the case where $\{x\}\in\mathfrak{A}$ and $\lambda(\{x\})=0$ for all $x\in X$.
If we consider the discrete topology on $X$ or equivalently
if we set $\mathfrak{A}(x) := \{A\in\mathfrak{A}\mid x\in A\}$,
then $\supp(u) = \emptyset$ for all $u\in H$ and $\supp(T) = \emptyset$ for all $T\in\mathcal{B}(H)$.
This demonstrates why we introduced $\mathfrak{A}$-neighborhoods
and why we do not use the definition of support given in \cite{RoeCG}.
\end{rem}

In order to get a calculus of supports as in Proposition~\ref{calcsupp},
we need $\lambda$ to fulfill an additional condition.

\begin{rem}
Let $\lambda\colon \mathfrak{A}\to\mathcal{B}(H)$ be any geometric Hilbert space.
For $A,A_1,A_2\in\mathfrak{A}$, $u\in H$ and $T\in\mathcal{B}(H)$ the following is always true.
\begin{eqnarray*}
\lambda(A)u = 0 & \Longrightarrow & \kringel{A}\cap\supp(u) = \emptyset \\
\lambda(A_1)\,T\,\lambda(A_2) = 0 & \Longrightarrow & \kringel{A_1}\!\times\! \kringel{A_2} \cap \supp(T) = \emptyset
\end{eqnarray*}
\end{rem}

\begin{defn}
We say that a geometric Hilbert space leads to \emph{computable supports}\index{computable supports}
if for $A,A_1,A_2\in\mathfrak{A}$, $u\in H$ and $T\in\mathcal{B}(H)$ the following is true.
\begin{eqnarray}
A\cap\supp(u) = \emptyset & \Longrightarrow & \lambda(A)u = 0 \label{gl21} \\
A_1\!\times\! A_2 \cap \supp(T) = \emptyset & \Longrightarrow & \lambda(A_1)\,T\,\lambda(A_2) = 0 \label{gl22}
\end{eqnarray}
\end{defn}

\begin{rem}\label{compsuppeasy}
In the situation of Example~\ref{decompositionPVM}, supports are computable
(provided that the Hilbert spaces $H_i$ are separable).
\end{rem}

\begin{lemma} \label{goodSupport}
Let $\lambda\colon\mathfrak{A}\to\mathcal{B}(H)$ be a geometric Hilbert space over $(X,\mathcal{E})$.
Assume that $X$ is $\sigma$-compact.\footnote{
A topological space is $\sigma$-compact if every open cover of $X$ has a countable subcover.}
Then the geometric Hilbert space $\lambda$ leads to computable supports.
\end{lemma}
\textbf{Proof.}
Assume $A\in\mathfrak{A}$ and $A\cap\supp(u) = \emptyset$.
For every $x\in A$ there is $B_x\in\mathfrak{A}(x)$
such that $\lambda(B_x)u = 0$. Set $C_x = B_x\cap A$.
Using Lemma~\ref{proj_commute}, we get
$$ \lambda(C_x)u = \lambda(C_x) \underbrace{\lambda(B_x) u}_{=0} = 0 \ . $$
Take a countable subcover $\mathcal{C}' \subseteq
\mathcal{C}:=\left\{ C_x \mid x\in A \right\}$ of $A$.
Using $\sigma$-additivity of $\lambda$, we get $\lambda(A)u = 0$.

The proof of the corresponding claim about the support of $T$ is similar.
\qua

\bigskip

The assumption of $\sigma$-compactness of $X$ is necessary in Lemma~\ref{goodSupport}.
Compare Remark~\ref{badsupp}.

\begin{prop}\label{calcsupp}
Let $\lambda\colon\mathfrak{A}\to\mathcal{B}(H)$ be a geometric Hilbert space
over $(X,\mathcal{E})$ which leads to computable supports.
Let $u, v\in H$ and $T,S\in\mathcal{B}(H)$.
Take $\Delta_\mathcal{U} = \bigcup_{U\in\mathcal{U}} U\times U\in\mathcal{E}$
for some localizing decomposition $\mathcal{U}$.
We get the following calculus of supports:
\begin{eqnarray}
\supp(u + v) & \subseteq & \supp(u) \cup \supp(v) \label{gl23} \\
\supp(S + T) & \subseteq & \supp(S) \cup \supp(T) \label{gl24} \\
\supp(Tu) & \subseteq & \Delta_\mathcal{U}\supp(T)\Delta_\mathcal{U}[\supp(u)] \label{gl25} \\
\supp(ST) & \subseteq & \Delta_\mathcal{U}\supp(S)\Delta_\mathcal{U}\supp(T)\Delta_\mathcal{U} \label{gl26} \\
\supp(T^*) & = & \supp(T)^{-1} \label{gl27} 
\end{eqnarray}
\end{prop}
\textbf{Proof.}
The first two inclusions (\ref{gl23}) and (\ref{gl24}) follow directly from the definitions.
The same is true for (\ref{gl27}).

In order to prove (\ref{gl25}), suppose $x\not\in \Delta_\mathcal{U}\supp(T)\Delta_\mathcal{U}\supp(u)$
and choose $U_x\in\mathcal{U}$ such that $x\in U_x$.
For every $U\in\mathcal{U}\backslash\{U_x\}$ we have
$$U\cap\supp(u)=\emptyset\qquad\text{or}\qquad U_x\times U\cap \supp(T)=\emptyset.$$
Hence $\lambda(U_x)Tu = \sum_{U\in\mathcal{U}}\lambda(U_x)T\lambda(U)u = 0$,
i.e. $x\not\in\supp(Tu)$.

The prove of (\ref{gl26}) is similar.
Suppose $(x,y)\not\in \Delta_\mathcal{U}\supp(S)\Delta_\mathcal{U}\supp(T)\Delta_\mathcal{U}$
and choose $U_x,U_y\in\mathcal{U}$ containing $x$ and $y$ respectively.
For every $U\in\mathcal{U}$ we have
$$U_x\times U\cap\supp(S) = \emptyset\qquad\text{or}\qquad U\times U_y\cap\supp(T) = \emptyset.$$
Therefore $\lambda(U_x)ST\lambda(U_y) = \sum_{U\in\mathcal{U}}\lambda(U_x)S\lambda(U)T\lambda(U_y)=0$
and it follows that $(x,y)\not\in\supp(ST)$.
\qua

\section{Translation $C^*$-algebras}

In the following let $\lambda\colon\mathfrak{A}\to\mathcal{B}(H)$ be
a geometric Hilbert space over $(X,\mathcal{E})$ with computable supports.
From Proposition~\ref{calcsupp} we get the following corollary.

\begin{cor}
The bounded operators of controlled support on $H$,
i.e. the bounded operators on $H$ whose support is an entourage, 
form a $*$-algebra.
\end{cor}

\begin{defn} \label{defnE}\index{E@$E^*(\ldots)$}\index{translation $C^*$-algebra}
We denote the closure of the $*$-algebra of operators with controlled support in the operator norm
by $E^*_\lambda (X,\mathcal{E})$. We may call it the \emph{big translation $C^*$-algebra}.
\end{defn}

\begin{rem}
The closure of the $*$-algebra of operators with controlled support with respect to strong convergence
is the full algebra of bounded operators $\mathcal{B}(H)$.
\end{rem}
\textbf{Proof.}
Let $\{U_i\}_{i\in\NNN}$ be an $\mathfrak{A}$-localizing decomposition of $(X,\mathcal{E})$
and set $V_i := U_0 \cup\cdots\cup U_i$.
Observe that for $T\in\mathcal{B}(H)$
the sequence $\{\lambda(V_i)T\lambda(V_i)\}_{i\in\NNN}$ converges strongly to $T$.
\qua

\begin{defn}
An operator $T\in\mathcal{B}(H)$ is \emph{pseudolocal}
if $T\lambda(B)-\lambda(B)T$ is a compact operator for every bounded set $B\in\mathfrak{A}$.
\end{defn}

\begin{lemma}
The pseudolocal operators with controlled support form a $*$-algebra.
\end{lemma}
\textbf{Proof.}
The pseudolocal operators are closed under addition, scalar multiplication and adjoints,
i.e. if $S$ and $T$ are pseudolocal operators and $z\in\CCC$,
then $S+T$, $z\cdot T$ and $T^*$ are also pseudolocal.
The composition of pseudolocal operators is pseudolocal, as the following computation shows:
If $S$ and $T$ are pseudolocal operators and $B\in\mathcal{A}$ is a bounded set, then
$$ST\lambda(B)-\lambda(B)ST = \big(S\lambda(B)-\lambda(B)S\big)T + \text{compact operator}$$
is compact.
\qua

\begin{defn} \label{defnD}
By $D^*_\lambda (X,\mathcal{E})$ we denote
the closure of the $*$-algebra of pseudolocal operators with controlled support
in the operator norm.
\end{defn}

\begin{defn}
An operator $T\in\mathcal{B}(H)$ is \emph{locally compact}
if $T\lambda(B)$ and $\lambda(B)T$ are compact operators for every bounded set $B\in\mathfrak{A}$.
\end{defn}

\begin{lemma}
The locally compact operators with controlled support form a $*$-algebra.
\end{lemma}
\textbf{Proof.}
The locally compact operators are closed under addition, scalar multiplication and adjoints.
The composition of locally compact operators is locally compact,
since composing a bounded operator and a compact operator gives a compact operator.
\qua

\begin{defn} \label{defnC}
The \emph{small translation $C^*$-algebra} of $(X,\mathcal{E})$ with respect to $\lambda$
is the closure of the $*$-algebra of locally compact operators with controlled support
in the operator norm.
We denote it by $C^*_\lambda (X,\mathcal{E})$.
This algebra is often called \emph{Roe-$C^*$-algebra}.\index{Roe-$C^*$-algebra}
\end{defn}

\begin{rem}
Observe that every locally compact operator is pseudolocal.
Hence $C^*_\lambda(X,\mathcal{E}) \subseteq D^*_\lambda(X,\mathcal{E}) \subseteq E^*_\lambda(X,\mathcal{E})$.
\end{rem}

\begin{lemma}
The small translation $C^*$-algebra $C^*_\lambda(X,\mathcal{E})$
is an ideal\footnote{In the context of $C^*$-algebras, by an \emph{ideal} we mean
a closed, two-sided ideal which is closed unter taking adjoints.}
in the big translation $C^*$-algebra $E^*_\lambda(X,\mathcal{E})$.
The small translation $C^*$-algebra $C^*_\lambda(X,\mathcal{E})$
is also an ideal in $D^*_\lambda(X,\mathcal{E})$.
\end{lemma}
\textbf{Proof.}
Let $S\in C^*_\lambda(X,\mathcal{E})$ and $T\in E^*_\lambda(X,\mathcal{E})$.
For a bounded set $B\in\mathfrak{A}$ the operator $\lambda(B)S$ is compact.
Therefore $\lambda(B)ST$ is compact.

It remains to prove that $ST\lambda(B_2)$ is compact for any bounded set $B_2\in\mathfrak{A}$.
Let $\mathcal{U}$ be an $\mathfrak{A}$-localizing decomposition of $(X,\mathcal{E})$ and
observe that $B_1:=\Delta_\mathcal{U}\supp(T)[B_2]\in\mathfrak{A}$ is bounded and
$(X\backslash B_1)\times B_2\cap\supp(T) = \emptyset$.
Using (\ref{gl22}), we see that
$$ S\, T\lambda(B_2) 
      = \underbrace{S\lambda(B_1)}_{\text{compact}}T\lambda(B_2) 
      + S\underbrace{\lambda(X\backslash B_1)\, T\lambda(B_2)}_{=0}\ $$
is compact.
\qua

\begin{example}
If $(X,\mathcal{E})$ is a coarse space such that $X$ is bounded, then
$E^*(X,\mathcal{E})$ is the entire algebra of bounded operators on $H$
and $C^*(X,\mathcal{E})$ is the algebra of compact operators.
Moreover, $D^*_\lambda(X,\mathcal{E}) = \{ T\in\mathcal{B}(H) \mid [T,\lambda(A)] \text{ is compact for all } A \in \mathfrak{A}\}$.
\end{example}

\begin{defn} \label{def_ample}\index{ample}
A geometric Hilbert space $\lambda\colon\mathfrak{A}\to\mathcal{B}(H)$ over $(X,\mathcal{E})$
is called \emph{ample} if there is an $\mathfrak{A}$-localizing decomposition $\mathcal{U}\subseteq\mathfrak{A}$
such that $\lambda(U)$ is not compact for $U\in\mathcal{U}$.
Note that for a projection being compact and having finite dimensional image is the same.
\end{defn}

\begin{example} \label{badPVM}
Let $(X,\mathcal{E})$ be a coarse space and $x_0\in X$.
The geometric Hilbert space $\lambda\colon \mathfrak{A}\to\mathcal{B}(H)$ defined by
\begin{displaymath}
\lambda(A) = \left\{
\begin{array}{ll}
0 & \text{if } x_0\not\in A \\
\id & \text{if } x_0\in A
\end{array}
\right.
\end{displaymath}
is not ample if $X$ is unbounded.
\end{example}

\begin{thm} \label{independent_of_lambda}
Let $\lambda\colon\mathfrak{A}\to\mathcal{B}(H_1)$ and $\mu\colon\mathfrak{B}\to\mathcal{B}(H_2)$
be ample, geometric Hilbert spaces over $(X,\mathcal{E})$ with computable supports. Then
$$ C^*_\lambda(X,\mathcal{E}) \cong C^*_\mu(X,\mathcal{E})
\qquad \text{and} \qquad
E^*_\lambda(X,\mathcal{E}) \cong E^*_\mu(X,\mathcal{E}) .$$
\end{thm}

The following two lemmas yield a proof of this theorem.

\begin{lemma}\label{CDgleich}
Let $\lambda\colon\mathfrak{A}\to\mathcal{B}(H)$ be a geometric Hilbert space
over $(X,\mathcal{E})$ with computable supports
and $\mathcal{U}$ an $\mathfrak{A}$-localizing decomposition.
By $\lambda'$ we denote the restriction of $\lambda$
to the $\sigma$-algebra $\sigma(\mathcal{U})$ generated by $\mathcal{U}$.
Then $\lambda'$ together with the discrete topology on $X$ is another geometric Hilbert space
and we have
$C^*_\lambda(X,\mathcal{E})=C^*_{\lambda'}(X,\mathcal{E})$ and
$E^*_\lambda(X,\mathcal{E})=E^*_{\lambda'}(X,\mathcal{E})$.
\end{lemma}
\textbf{Proof.}
Whenever in this proof we make reference to a topology we mean the topology coming with $\mathfrak{A}$.
We first claim $$\supp_\lambda(T)\subseteq\overline{\supp_{\lambda'}(T)}\quad\text{ for } T\in\mathcal{B}(H).$$
This implies $E^*_{\lambda'}(X,\mathcal{E})\subseteq E^*_\lambda(X,\mathcal{E})$.
To prove the claim, set $Y:=\bigcup_{U\in\mathcal{U}}\kringel{U}$.
If $x\in Y$, then $\sigma(\mathcal{U})(x)=\{U\in\mathcal{U}\mid x\in U\}\subseteq\mathfrak{A}(x)$.
Therefore we get $\supp_\lambda(T)\cap {Y\times Y}\subseteq\supp_{\lambda'}(T)\cap{Y\times Y}$.

For $x\in X$ set $\mathcal{U}_x := \{U\in\mathcal{U}\mid x\in\overline{U}\}$ and observe that
$\bigcup_{U\in\mathcal{U}_x}U$ is an $\mathfrak{A}$-neighborhood of $x$.
Thus, for $(x,y)\in\supp_\lambda(T)$ there is $U\in\mathcal{U}_x$ and $V\in\mathcal{U}_y$ such that
$\lambda(U)T\lambda(V)\neq 0$.
Take sequences $\{x_i\}_{i\in\NNN}\subseteq U$ and $\{y_i\}_{i\in\NNN}\subseteq V$
converging to $x$ and $y$ respectively.
Note that $(x_i,y_i)\in\supp_{\lambda'}(T)$ for all $i\in\NNN$ and hence
$(x,y)\in\overline{\supp_{\lambda'}(T)}$.

Set $\Delta_\mathcal{U} = \bigcup_{U\in\mathcal{U}} U\times U$.
We claim
$$\supp_{\lambda'}(T) \subseteq \Delta_\mathcal{U}  \supp_\lambda(T)  \Delta_\mathcal{U}
\quad\text{ for all } T\in\mathcal{B}(H)$$
This claim implies $E^*_\lambda(X,\mathcal{E}) \subseteq E^*_{\lambda'}(X,\mathcal{E})$.
To prove the claim, suppose $U,V\in\mathcal{U}$ and
${U \times V} \cap \supp_\lambda(T) = \emptyset$.
Since supports are computable, we get $\lambda(U) T \lambda(V) = 0$.
Therefore $U \times V \cap \supp_{\lambda'}(T) = \emptyset$.

So far we proved $E^*_\lambda(X,\mathcal{E}) = E^*_{\lambda'}(X,\mathcal{E})$.

In order to prove $C^*_\lambda(X,\mathcal{E}) \subseteq C^*_{\lambda'}(X,\mathcal{E})$,
just note that each bounded set in $\sigma(\mathcal{U})$ is also a bounded set in $\mathfrak{A}$.

Each bounded set $B\in\mathfrak{A}$ is contained in $\widetilde{B}:=\Delta_\mathcal{U}[B]$
which is a bounded set in $\sigma(\mathcal{U})$.
Let $T$ be locally compact with respect to $\sigma(\mathcal{U})$.
Using Lemma~\ref{proj_commute}, we see that $\lambda(B) T = \lambda(B)\lambda(\widetilde{B}) T$ is compact,
since $\lambda(\widetilde{B}) T$ is compact.
\qua

\begin{lemma}\label{seclemthm}
Given localizing decompositions $\mathcal{U}$ and $\mathcal{V}$ of the coarse space $(X,\mathcal{E})$
let $\lambda\colon\sigma(\mathcal{U})\to\mathcal{B}(H_1)$ and $\mu\colon\sigma(\mathcal{V})\to\mathcal{B}(H_2)$
be ample, geometric Hilbert spaces over $(X,\mathcal{E})$. Then there is an isometric isomorphism
$\rho\colon H_1\to H_2$ such that $ad_\rho\colon\mathcal{B}(H_1)\to\mathcal{B}(H_2), \ T\mapsto \rho T \rho^{-1}$
restricts to isomorphisms $E^*_\lambda(X,\mathcal{E}) \overset{\cong}{\longrightarrow} E^*_\mu(X,\mathcal{E})$ and
$C^*_\lambda(X,\mathcal{E}) \overset{\cong}{\longrightarrow} C^*_\mu(X,\mathcal{E})$.
\end{lemma}
\textbf{Proof.}
Defining $\mathcal{W} = \{ U\cap V \mid U\in\mathcal{U}, V\in\mathcal{V} \} \ \backslash \ \{ \emptyset \}$
we get a localizing decomposition refining $\mathcal{U}$ and $\mathcal{V}$.
Set $\mathcal{W}_U = \{ W\in\mathcal{W} \mid W\subseteq U \}$ for $U\in\mathcal{U}$ and
$\mathcal{W}_V = \{ W\in\mathcal{W} \mid W\subseteq V \}$ for $V\in\mathcal{V}$.
We will write $H_U$ for the image of $\lambda(U)$ and $H_V$ for the image of $\mu(V)$.
For every $U\in\mathcal{U}$ we choose a decomposition $H_U = \bigoplus_{W\in\mathcal{W}_U}H_{1,W}$
such that $H_{1,W}$ is infinite dimensional for all $W\in\mathcal{W}$.
Similarly, we choose a decomposition $H_V = \bigoplus_{W\in\mathcal{W}_V}H_{2,W}$ for each $V\in\mathcal{V}$.
Now we extend the geometric Hilbert spaces $\lambda$ and $\mu$
to $\widetilde{\lambda}\colon\sigma(\mathcal{W}) \to \mathcal{B}(H_1)$
and $\widetilde{\mu}\colon\sigma(\mathcal{W}) \to \mathcal{B}(H_2)$ respectively.
We define $\widetilde{\lambda}(W)$ to be the projection onto $H_{1,W}$.
Similarly $\widetilde{\mu}(W)$ will be the projection onto $H_{2,W}$.
Remember that $\lambda$, $\widetilde{\lambda}$, $\mu$ and $\widetilde{\mu}$ have computable supports.
We choose an isometric isomorphism $\rho\colon H_1\to H_2$ such that
the image of $\rho|_{H_{1,W}}$ is $H_{2,W}$ for all $W\in\mathcal{W}$.
Then $ad_\rho$ is an isomorphism of $C^*$-algebras.

Using $\rho\widetilde{\lambda} = \widetilde{\mu}\rho$, we see
that $\supp_{\widetilde{\lambda}}(T) = \supp_{\widetilde{\mu}}(\rho T\rho^{-1})$
for all $T\in\mathcal{B}(H_1)$.
Hence $ad_\rho$ restricts to an isomorphism
$E^*_{\widetilde{\lambda}}(X,\mathcal{E}) \overset{\cong}{\longrightarrow} E^*_{\widetilde{\mu}}(X,\mathcal{E})$.

Moreover, we note that for $B\in\sigma(\mathcal{W})$ and $T\in\mathcal{B}(H_1)$ 
the operator $\widetilde{\lambda}(B) T$ is compact if and only if
$\widetilde{\mu}(B) ad_\rho(T)$ is compact.
Therefore $ad_\rho$ also restricts to an isomorphism
$C^*_{\widetilde{\lambda}}(X,\mathcal{E}) \overset{\cong}{\longrightarrow} C^*_{\widetilde{\mu}}(X,\mathcal{E})$.
Applying Lemma~\ref{CDgleich} now completes the proof.
\qua

\bigskip

There are some choices involved in order to get the isometric isomorphism $\rho\colon H_1\to H_2$
in the proof of Lemma~\ref{seclemthm}. Therefore, the isomorphisms
$E^*_\lambda(X,\mathcal{E}) \overset{\cong}{\longrightarrow} E^*_\mu(X,\mathcal{E})$ and 
$C^*_\lambda(X,\mathcal{E}) \overset{\cong}{\longrightarrow} C^*_\mu(X,\mathcal{E})$
are not uniquely determined and without further specifications there is no natural choice.

\bigskip

We compare our definition of small translation $C^*$-algebra
with the corresponding definition in Section 6.3 of \cite{HR}.

\begin{prop}
Let $X$ be a locally compact, separable and metrizable space with a coarse structure $\mathcal{E}$
which is compatible with the topology. Suppose that 
$$\rho\colon C_0(X)\to \mathcal{B}(H)$$
is a non-degenerate\footnote{
  A $C^*$-representation $\rho$ of a $C^*$-algebra $A$ on a Hilbert space $H$
  is non-degenerate if $\{ \rho(a)(v) \mid a\in A, v\in H \}$ is dense in $H$.}
and ample\footnote{
  A $C^*$-representation $\rho$ is called ample,
  if the operator $\rho(f)$ is not compact if $f\neq 0$.} 
representation of the $C^*$-algebra $C_0(X)$ on the separable Hilbert space $H$.
There is an ample geometric Hilbert space $\lambda$ over $(X,\mathcal{E})$
with computable supports and for any such $\lambda$
$$C^*_\rho(X,\mathcal{E}) \cong C^*_\lambda(X,\mathcal{E})\ .$$
\end{prop}
\textbf{Proof.}
According to Theorem~\ref{independent_of_lambda} we are free to replace $\lambda$
by any ample geometric Hilbert space with computable supports.
We will construct $\lambda$ from $\rho$.

Using Borel functional calculus (compare Remark 1.5.7 of \cite{HR}),
we see that $\rho$ extends to a representation $\widetilde{\rho}$
of the $C^*$-algebra $B(X)$ of bounded Borel functions on $X$.
It is easy to check that
$$\lambda \colon \mathfrak{B}(X) \to \mathcal{B}(H), \quad A \mapsto \widetilde{\rho}(\chi_A)$$
is a projection-valued measure.

Since $X$ is separable and its coarse structure is compatible with the topology,
there exists a countable, uniformly bounded, open cover of $X$.
We may even assume that each $U\in\mathcal{U}$ has non-empty interior.
(Compare Claim 6.3.14 of \cite{HR}.)
Hence $\lambda$ is a geometric Hilbert space.

Let $U\in\mathcal{U}$ and $f\in C_0(U)$, $f\neq 0$.
Since $\rho$ is ample, $\lambda(U)\rho(f)=\rho(f)$ is not compact.
Hence $\lambda(U)$ is not compact and $\lambda$ is ample.
Applying Lemma~\ref{goodSupport} we see that $\lambda$ leads to computable supports.

Now that we have a convenient geometric Hilbert space,
we compare the support of a bounded operator $T$ with respect to $\lambda$
(as defined in Definition~\ref{defSupport})
and the support with respect to $\rho$ (as defined in \cite{HR}).
Observe that
$$\supp_\rho(T) = \left\{ (x_1,x_2)\in X^2 \mid
                     \forall_{\substack{U_1\in\mathfrak{B}(X)(x_1) \\ U_2\in\mathfrak{B}(X)(x_2)}}
                     \exists_{\substack{f_1\in C_0(\kringel{U}_1) \\ f_2\in C_0(\kringel{U}_2)}}
                     \rho(f_1) T \rho(f_2) \neq 0
                  \right\}\ .$$

Let $x_1,x_2\in X$ and $U_i\in\mathfrak{B}(X)(x_i)$.
Suppose there are $f_i\in C_0(\kringel{U_i})$ such that $\rho(f_1) T\rho(f_2)\neq 0$.
Observe the following relations between $\rho$ and $\lambda$.
\begin{eqnarray*}
\rho(f_1)\lambda(U_1) = \rho(f_1)\widetilde{\rho}(\chi_{U_1}) = \widetilde{\rho}(f_1\cdot\chi_{U_1}) = \rho(f_1) \\
\lambda(U_2)\rho(f_2) = \widetilde{\rho}(\chi_{U_2})\rho(f_2) = \widetilde{\rho}(\chi_{U_2}\cdot f_2) = \rho(f_2)
\end{eqnarray*}
It follows
$$\rho(f_1)\lambda(U_1) T\lambda(U_2)\rho(f_2) = \rho(f_1) T\rho(f_2) \neq 0$$
and  hence $\lambda(U_1) T \lambda(U_2) \neq 0$.
This proves $\supp_\rho(T)\subseteq\supp_\lambda(T)$.

Choose a metric $d$ which induces the given topology on $X$.
For $U\subseteq X$ and $k\in\NNN$ define
$$f_{U,k}\colon X\to \RRR, \quad x\mapsto \max\{\, 0,\, 1-k\cdot d(x,U)\, \}\ .$$

Let $(x_1,x_2)\in\supp_\lambda(T)$ and $U_i\in\mathfrak{B}(X)(x_i)$.
Since $(X,d)$ is a metric space, there is $k_i\in\NNN$ and $V_i\in\mathfrak{B}(X)(x_i)$
such that $f_i := f_{V_i,k_i}\in C_0(\kringel{U_i})$.
We have the following relations between $\rho$ and $\lambda$.
\begin{eqnarray*}
\lambda(V_1)\rho(f_1) = \widetilde{\rho}(\chi_{V_1})\rho(f_1) = \widetilde{\rho}(\chi_{V_1}\cdot f_1) = \lambda(V_1) \\
\rho(f_2)\lambda(V_2) = \rho(f_2)\widetilde{\rho}(\chi_{V_2}) = \widetilde{\rho}(f_2\cdot\chi_{V_2}) = \lambda(V_2)
\end{eqnarray*}
These relations imply $$\lambda(V_1)\rho(f_1) T \rho(f_2)\lambda(V_2) = \lambda(V_1) T \lambda(V_2) \neq 0$$
and hence $\rho(f_1) T \rho(f_2) \neq 0$.
It follows $\supp_\lambda(T)\subseteq\supp_\rho(T)$.

We assume that $T\in C^*_\rho(X,\mathcal{E})$.
Note that $f_A:=f_{A,1}\in C_0(X)$ for a bounded set $A\subseteq X$.
Since $\rho(f_A) T$ is compact, we conclude that $\lambda(\overline{A}) T = \lambda(\overline{A})\rho(f_A) T$ is compact.
In the same way we get compactness for $T \lambda(\overline{A})$.
This proves $T\in C^*_\lambda(X,\mathcal{E})$.

Now let $T\in\mathcal{B}(H)$ be locally compact with respect to $\lambda$.
If $f\in C_0(X)$, then $M_{f,n}:=\left\{x\in X\mid \abs{f(x)}\geq \frac{1}{n}\right\}$
is compact for all $n\in\NNN$ and thus bounded.
Continuity of the representation $\widetilde{\rho}$ yields
$$ \xymatrix@1{\rho(f) \lambda(M_{f,n}) T \ar[rrr]_(.6){n\to\infty}^(.6){\norm{\cdot}} & & & \rho(f) T}\ .$$
Since all operators on the left hand side are compact, $\rho(f) T$ is compact.
In the same way we get compactness for $T \rho(f)$.
\qua

\begin{cor}
Up to isomorphism, the $C^*$-algebra $C_\rho(X)$
does not depend on the choice of the representation $\rho$
(as long as $\rho$ is an ample and non-degenerate representation on a separable Hilbert space).
\end{cor}

\section{Induced maps}

Let $(X,\mathcal{E}_X)$ and $(Y,\mathcal{E}_Y)$ be coarse spaces which admit localizing decompositions
and let $f \colon (X,\mathcal{E}_X) \to (Y,\mathcal{E}_Y)$ be a coarse map.
We would like to obtain induced maps $E^*(f)$ and $C^*(f)$ on the translation $C^*$-algebras
such that the following diagramm commutes.
$$ \xymatrix{
C^*(X,\mathcal{E}_X) \ar[rr]^{C^*(f)}\ar@{^(->}[d] & & C^*(Y,\mathcal{E}_Y) \ar@{^(->}[d] \\
E^*(X,\mathcal{E}_X) \ar[rr]^{E^*(f)} & & E^*(Y,\mathcal{E}_Y)
} $$


\begin{question}
Let $f \colon (X,\mathcal{E}_X) \to (Y,\mathcal{E}_Y)$ be a coarse map
and $\mathfrak{A}$ a localizing $\sigma$-algebra with respect to $\mathcal{E}_X$.
Is there a way of defining an induced localizing $\sigma$-algebra $f_*(\mathfrak{A})$?
\end{question}

An attempt to define an induced localizing $\sigma$-ring\footnote{
   A $\sigma$-ring $\mathfrak{R}$ over a set $X$ is a family of subsets of $X$
   such that (1) $\emptyset\in\mathfrak{R}$,
   (2) $A\backslash B\in\mathfrak{R}$ whenever $A, B\in\mathfrak{R}$
   and (3) the union of countably many elements of $\mathfrak{R}$ is in $\mathfrak{R}$.
}
was made in Remark 4.26 of \cite{RoeCG}.
Roe defines $$f_*(\mathfrak{A}):=\{ S\subseteq f(X)\subseteq Y \mid f^{-1}(S)\in\mathfrak{A}\}\ .$$
Now $f_*(\mathfrak{A})$ is indeed a $\sigma$-ring,
but it does not have to be localizing as the following example tells us.

\begin{example}\label{badindsigalg}
Take $(X,\mathcal{E}_X)=(Y,\mathcal{E}_Y)=(\RRR,\mathcal{E}_{eucl})$ and
the $\sigma$-algebra $\mathfrak{A}$ generated by the set of half open intervals
$\{\left[n,n\!+\!1\right[\mid n\in\ZZZ\})$.
For the coarse map $f\colon X \to Y$
$$ f(x) = \left\{ \begin{array}{ll}
   n & \text{if } x\in\left[n-\frac{1}{4},n+\frac{1}{4}\right] \text { for some } n\in\ZZZ \\
   2\cdot x -n -\frac{1}{2} & \text{if } x\in\left[n+\frac{1}{4},n+\frac{3}{4}\right] \text { for some } n\in\ZZZ
\end{array} \right. $$
we get $f_*(\mathfrak{A})=\{\emptyset, \RRR\}$ which is not localizing
(in the sense of Definition 4.20 in \cite{RoeCG}), since $\RRR$ is unbounded.

Roe's definition of $f_*(\mathfrak{A})$ does give a ``localizing $\sigma$-ring''
if the coarse map $f$ ``respects an $\mathfrak{A}$-localizing decomposition'',
i.e. if $f$ has the following additional property:
There is an $\mathfrak{A}$-localizing decomposition $\mathcal{U}$ of $X$ such that
$U_1,U_2\in\mathcal{U}$ and $U_1 \neq U_2$ implies $f(U_1)\cap f(U_2) = \emptyset$.
This is the case for example if $f$ is injective.
\end{example}

We start constructing geometric Hilbert spaces $\lambda_X$ and $\lambda_Y$
such that $f$ does induce maps between the corresponding translation $C^*$-algebras.

\begin{lemma}\label{lodata}
There are geometric Hilbert spaces $\lambda_X\colon \sigma(\mathcal{U}_X) \to \mathcal{B}(H_X)$
and  $\lambda_Y\colon \sigma(\mathcal{U}_Y) \to \mathcal{B}(H_Y)$,
a map $F\colon\mathcal{U}_X\to\mathcal{U}_Y$
and a partial isometry $\varphi\colon H_X \to H_Y$ such that
\begin{itemize}
\item $\mathcal{U}_X$ and $\mathcal{U}_Y$ are localizing decompositions of $X$ and $Y$ respectively,
\item $\lambda_X$ and $\lambda_Y$ are ample and lead to computable supports,
\item the map $F\colon\mathcal{U}_X\to\mathcal{U}_Y$ satisfies $F(U)\subseteq f(U)$ for all $U\in\mathcal{U}_X$,
\item the map $F\colon\mathcal{U}_X\to\mathcal{U}_Y$ is almost injective,\\
i.e. $F(U) = F(U')$ implies $U=U'$ or $F(U) = \emptyset$,
\item the partial isometry satisfies $\varphi(\lambda_X(U)(H_X)) \subseteq \lambda_Y(F(U))(H_Y)$\\
for all $U\in\mathcal{U}_X$.
\end{itemize}
\end{lemma}
\textbf{Proof.}
Take any localizing decomposition $\mathcal{U}_X$ of $(X,\mathcal{E}_X)$.
We enumerate the elements of $\mathcal{U}_X$, denote them by $U_0,U_1,U_2,\ldots$
and set $W_i := f(U_i)$.
For $i\in\NNN$ define $V_i = W_i \backslash \bigcup_{j<i}W_j$ and note that
$\mathcal{V} := \{V_0,V_1,V_2,\ldots\}$ is a localizing decomposition of $\im(f)$.
The empty set might be an element of $\mathcal{V}$.
For $i\in\NNN$ define $F(U_i) := V_i$.
Let  $\widetilde{\mathcal{U}}_Y$ be a localizing decomposition of $Y$ and define
$$ \mathcal{U}_Y := \mathcal{V}\ \cup\ \left\{ U\backslash\im(f) \mid U\in\widetilde{\mathcal{U}}_Y\right\}.$$
Note that $\mathcal{U}_Y$ is itself a localizing decomposition of $Y$.

Choose a separable Hilbert space $H_U$ for each $U\in\mathcal{U}_X$ with $U\neq\emptyset$
and similarly choose $H_V$ for each $V\in\mathcal{U}_Y$ with $V\neq\emptyset$.
Define $H_{\emptyset} = \{0\}$ if $\emptyset\in\mathcal{U}_X$ or $\emptyset\in\mathcal{U}_Y$.
Define $H_X = \bigoplus_{U\in\mathcal{U}_X} H_U$ and $H_Y = \bigoplus_{V\in\mathcal{U}_Y} H_V$.
Choose a partial isometry $\varphi\colon H_X \to H_Y$ such that
$\varphi\left(H_U\right)\subseteq H_{F(U)}$ for all $U\in\mathcal{U}_X$.
We consider the geometric Hilbert space
$\lambda_X\colon\sigma(\mathcal{U}_X)\to\mathcal{B}\left( H_X \right)$
which is determined by defining $\lambda_X(U)$ to be the projection onto $H_U$ for any $U\in\mathcal{U}_X$.
Likewise we define $\lambda_Y\colon\sigma(\mathcal{U}_Y)\to\mathcal{B}\left( H_Y \right)$.
Note that $\lambda_X$ and $\lambda_Y$ are ample if $H_U$ and $H_V$
are infinite dimensional for all $U\in\mathcal{U}_X$ and $V\in\mathcal{U}_Y$.
Moreover $\lambda_X$ and $\lambda_Y$ lead to computable supports.
(Compare Remark~\ref{compsuppeasy}.)
\qua

\begin{lemma}\label{DCmaps}
In the situation of Lemma~\ref{lodata}, the partial isometry $\varphi$ induces the map
$$ ad_\varphi\colon \mathcal{B}(H_X) \rightarrow \mathcal{B}(H_Y) , \quad T\mapsto\varphi T \varphi^*$$
and $ad_\varphi$ restricts to maps between the translation $C^*$-algebras.
\begin{eqnarray}
\xymatrix@1{
E^*_{\lambda_X}(X,\mathcal{E}_X) \ \ar[rr]^(.5){ad_\varphi} & & \ E^*_{\lambda_Y}(Y,\mathcal{E}_Y)
} \label{glD} \\
\xymatrix@1{
C^*_{\lambda_X}(X,\mathcal{E}_X) \ \ar[rr]^(.5){ad_\varphi} & & \ C^*_{\lambda_Y}(Y,\mathcal{E}_Y)
} \label{glC}
\end{eqnarray}
\end{lemma}
\textbf{Proof.}
Let $T\in\mathcal{B}(H_X)$. We will prove
$$ \supp_{\lambda_Y}\left(ad_\varphi(T)\right)\subseteq f\times f\left( \supp_{\lambda_X}(T) \right)\ .$$
This implies $ad_\varphi\left(E^*_{\lambda_X}(X,\mathcal{E}_X)\right) \subseteq E^*_{\lambda_Y}(Y,\mathcal{E}_Y)$.

Let $U_1,U_2\in\mathcal{U}_X$ such that $U_1\times U_2\not\subseteq\supp_{\lambda_X}(T)$.
The definition of $\lambda_X$ implies $U_1\times U_2\cap \supp_{\lambda_X}(T) = \emptyset$
and hence $\lambda_X(U_1)T\lambda_X(U_2)=0$, since supports are computable.
Observe that
\begin{eqnarray*}
&   & \lambda_Y(F(U_1)) \ ad_\varphi(T) \ \lambda_Y(F(U_2))
\ = \ \lambda_Y(F(U_1)) \ \varphi \ T \ \varphi^* \ \lambda_Y(F(U_2)) \\
& = & \lambda_Y(F(U_1)) \ \varphi \left( \sum_{U\in\mathcal{U}_X}\lambda_X(U) \right) 
       T \left( \sum_{U\in\mathcal{U}_X}\lambda_X(U) \right) \varphi^* \ \lambda_Y(F(U_2)) \\
& = & \lambda_Y(F(U_1)) \ \varphi \ \underbrace{\lambda_X(U_1) \ 
       T \ \lambda_X(U_2)}_{=0} \ \varphi^* \ \lambda_Y(F(U_2))
\ = \ 0 \ .
\end{eqnarray*}

If $(y_1,y_2)\in\supp_{\lambda_Y}\left( ad_\varphi(T) \right)$,
there are $U_1,U_2\in\mathcal{U}_X$ such that $y_i\in F(U_i)$.
Hence $\lambda_Y(F(U_1))ad_\varphi(T)\lambda_Y(F(U_2))\neq 0$.
As we have just seen, this implies $U_1\times U_2 \subseteq \supp_{\lambda_X}(T)$ and therefore
$$(y_1,y_2)\ \in\ F(U_1)\times F(U_2)\ \subseteq\ f(U_1)\times f(U_2)\ \subseteq\ f\times f\left(\supp_{\lambda_X}(T)\right) .$$

It remains to prove
$ad_\varphi\left(C^*_{\lambda_X}(X,\mathcal{E}_X)\right) \subseteq C^*_{\lambda_Y}(Y,\mathcal{E}_Y)$.
For this purpose suppose $T\in C^*_{\lambda_X}(X,\mathcal{E}_X)$. Let $B\in\sigma(\mathcal{U}_Y)$ be a bounded set.
Then $\lambda_X\left(\Delta_{\mathcal{U}_X}[f^{-1}(B)]\right) T$ is compact, since
$\Delta_{\mathcal{U}_X}[f^{-1}(B)]\in\sigma(\mathcal{U}_X)$ is bounded
and $T$ is locally compact. Therefore
$$ \lambda_Y(B) ad_\varphi(T) = \lambda_Y(B) \varphi T \varphi^*
= \lambda_Y(B)\varphi\underbrace{\lambda_X\left(\Delta_{\mathcal{U}_X}[f^{-1}(B)]\right) T}_{\text{compact}} \varphi^* $$
and $ad_\varphi(T)\lambda_Y(B)$ are compact operators.
This proves that $ad_\varphi(T)$ is locally compact.
\qua

\bigskip

We denote the maps (\ref{glD}) and (\ref{glC}) by
$E^*_{\lambda_X,\lambda_Y,\varphi}(f)$ and $C^*_{\lambda_X,\lambda_Y,\varphi}(f)$ respectively.
Note that these maps do not only depend on $\lambda_X$ and $\lambda_Y$,
but also on the partial isometry $\varphi$.

\begin{prop}
Let $(X,\mathcal{E}_X)$ and $(Y,\mathcal{E}_Y)$ be coarse spaces which
admit countable, uniformly bounded covers.
If $(X,\mathcal{E}_X)$ and $(Y,\mathcal{E}_Y)$ are coarsely equivalent, then
$C^*(X,\mathcal{E}_X) \cong C^*(Y,\mathcal{E}_Y)$
and $E^*(X,\mathcal{E}_X) \cong E^*(Y,\mathcal{E}_Y)$.
\end{prop}
\textbf{Proof.}
Let $f\colon (X,\mathcal{E}_X)\to(Y,\mathcal{E}_Y)$ be a coarse equivalence
with coarse inverse $g\colon (Y,\mathcal{E}_Y)\to(X,\mathcal{E}_X)$.
Take a localizing decomposition $\mathcal{W}$ of $\im(f\circ g)\subseteq Y$
which does not contain the empty set.
Define
$$ E := \Delta_\mathcal{W} \cup \{ (y,f\circ g(y)) \mid y\in Y\backslash \im(f\circ g) \} $$
and observe that $\mathcal{U}_Y := \{ E[W] \mid W\in\mathcal{W} \}$ is a localizing decomposition
of the coarse space $(Y,\mathcal{E}_Y)$.

Let $U\in\mathcal{U}_Y$. Since there is $W\in\mathcal{W}$ with $W\subseteq U$,
the inverse image of $U$ under $f$ is not empty.
Define $\mathcal{U}_X := \{ f^{-1}(U) \mid U\in\mathcal{U}_Y \}$.
Choose an infinite dimensional, separable Hilbert space $H_U$ for each $U\in\mathcal{U}_Y$
and set $H:=\bigoplus_{U\in\mathcal{U}_Y}H_U$.
Consider the geometric Hilbert space $\lambda_Y\colon\sigma(\mathcal{U}_Y)\to\mathcal{B}(H)$
which is determined by setting $\lambda(U)$ to be the projection onto $H_U$.
Furthermore, consider the geometric Hilbert space $\lambda_X\colon\sigma(\mathcal{U}_X)\to\mathcal{B}(H)$
which is given by defining $\lambda(f^{-1}(U))$ to be the projection onto $H_U$.

If we define $F\colon\mathcal{U}_X\to\mathcal{U}_Y$ by $F(f^{-1}(U)) = U$ for all $U\in\mathcal{U}_Y$
and $\varphi = \id_H$, the conclusions of Lemma~\ref{lodata} are satisfied.
Moreover, the map $ad_\varphi$ is just the identity on $\mathcal{B}(H)$.
Hence, Lemma~\ref{DCmaps} yields
$E^*_{\lambda_X}(X,\mathcal{E}_X) \subseteq E^*_{\lambda_Y}(Y,\mathcal{E}_Y)$ and
$C^*_{\lambda_X}(X,\mathcal{E}_X) \subseteq C^*_{\lambda_Y}(Y,\mathcal{E}_Y)$.

The reverse inclusion for the big translation $C^*$-algebras follows from
$$\supp_{\lambda_X}(T) = (f\times f)^{-1} \left( \supp_{\lambda_Y}(T) \right) \quad\text{ for all } T\in\mathcal{B}(H)$$
and the fact that $f$ is a coarse embedding.

Observe that $\lambda_Y(B)T = \lambda_X(f^{-1}(B))T$ for each bounded sets $B\in\sigma(\mathcal{U}_Y)$.
We conclude that an operator $T\in\mathcal{B}(H)$ is locally compact with respect to $Y$
if and only if $T$ is locally compact with respect to $X$.
\qua

\chapter{General theory of asymptotic dimension}

The following definitions are an asymptotic analog of the covering dimension of topological
spaces. Asymptotic dimension was first introduced in \cite{Gromov}.
For dimension theory of topological spaces see \cite{HuWa}, \cite{Engelking} and \cite{Fedorchuk}.

\bigskip

Finiteness of asymptotic dimension seems to play an important role
for some isomorphism conjectures in K-theory.
Guoliang Yu proved in \cite{Yufin} that the Baum-Connes assembly map is injective
for groups of finite asymptotic dimension admitting a finite classifying space.
This implies the Novikov conjecture on the homotopy invariance of higher signatures.
Arthur Bartels proved injectivity of the assembly map in algebraic K-theory
for the same class of groups. Compare \cite{bartels}.

\section{Asymptotic dimension of pseudometric spaces}

\begin{defn} \label{metric_asdim}
\textbf{(asymptotic dimension of pseudometric spaces)} \\
Let $(X,d)$ be a pseudometric space, $\mathcal{M}$ a collection of subsets of $X$ and $n \in \NNN$.
\begin{itemize}
\item
The \emph{multiplicity}\index{multiplicity} of $\mathcal{M}$
is defined to be the maximal number of sets with non-empty
intersection and will be denoted by $\mu(\mathcal{M})$.
\item
The \emph{mesh}\index{mesh} of $\mathcal{M}$ is defined to be
the supremum of the diameter of sets from $\mathcal{M}$.
\item
The \emph{Lebesgue number}\index{Lebesgue number} of a cover $\mathcal{U}$ of $X$
is the largest positive number $L$ such that for every $x \in X$
the open ball of radius $L$ with center $x$ is contained in some $U \in \mathcal{U}$.
We will write $L(\mathcal{U})$ for the Lebesgue number of $\mathcal{U}$. 
\item\index{asdim}
$\asdim(X,d) \leq n$ if for all $L>0$ there exist $D>0$ and a cover $\mathcal{U}$ of $X$ such that
   \begin{enumerate}
   \item[(1)] \ $\mu(\mathcal{U}) \leq n+1$,
   \item[(2)] \ $L(\mathcal{U}) \geq L$ and
   \item[(3)] \ $\mesh(\mathcal{U}) \leq D$,
                i.e. $\mathcal{U}$ is uniformly bounded.\index{uniformly bounded}
   \end{enumerate}
\item
$\asdim(X,d)=n$ if $\asdim(X,d) \leq n$ and $\asdim(X,d) \not\leq n-1$.
\item
$\asdim(X,d)=\infty$ if there is no $n \in \NNN$ with $\asdim(X,d)\leq n$.
\end{itemize}
\end{defn}
In Section~\ref{generalizeasdim} we will generalize this definition to coarse spaces.

As in dimension theory of topological spaces, there are some alternative definitions.
For details see \cite{Gromov} and \cite{Dran}.
We will give a second definition of asymptotic dimension,
which has been generalized to coarse spaces in \cite{RoeCG}.
We will prove that both definitions coincide, in the more general setting of coarse spaces.
\begin{defn} \label{metric_asdim_2}
Let $(X,d)$ be a pseudometric space.
\begin{itemize}
\item
Let $L>0$.
A family $\mathcal{V}$ of subsets of $X$ is called $L$-disjoint,\index{disjoint@$L$-disjoint}
if the distance of two sets from $\mathcal{V}$ is always bigger than $L$.
\item\index{asdim}
$\asdim(X,d) \leq n$ if for all $L>0$ there is a cover $\mathcal{U}$ of $X$ such that
   \begin{enumerate}
   \item[(1)] \ the cover $\mathcal{U}$ consists of $n+1$ families $\mathcal{U}_1,\ldots,\mathcal{U}_{n+1}$,
   \item[(2)] \ each family $\mathcal{U}_i$ is $L$-disjoint and
   \item[(3)] \ the cover $\mathcal{U}$ is uniformly bounded.
   \end{enumerate}
\end{itemize}
\end{defn}

\begin{rem}[monotony] \label{asdim_subset}\index{monotony}
Let $(X,d)$ be a pseudometric space.

If $(A,d|_A)$ is a subspace of $(X,d)$, then $\asdim(A,d|_A)\leq\asdim(X,d)$.
\end{rem}

\begin{example}\index{tree}
Let $T$ be a tree and $d$ the natural metric on $T$.
(All edges are supposed to be of length $1$.)
Then $\asdim(T,d) \leq 1$ and $\asdim(T,d) = 1$ if and only if $T$ is not bounded.
\end{example}
\textbf{Proof.}
Choose $x_0\in vert(T)$. Let $L>0$ and define $L'$ to be
the smallest natural number bigger than $2\cdot L$.
Consider the map
$$f \colon T \to \NNN, \ \ 
x \mapsto \left[\frac{d(x_0,x)}{L'}\right].$$
We define $A:=\overline{f^{-1}(2\cdot\NNN)}$ and
$B:=\overline{f^{-1}(2\cdot\NNN+1)}$, thus $X = A \cup B$.
There is an equivalence relation on $A$:
Let $a_1, a_2 \in A$. With $\gamma_i$ we denote the geodesic path
from $x_0$ to $a_i$.
$$ a_1 \sim a_2 \iff f(a_1)=f(a_2) \ \text{ and } \
\gamma_1(\tau)=\gamma_2(\tau)$$
where $\tau := L'\cdot\left(f(a_1)-\frac{1}{2}\right)$.
If $a_1\not\sim a_2$, then $d(a_1,a_2)\geq L'$.
Furthermore, we see $\diam([a]) \leq 3\cdot L'$.
Here $[a]$ denotes the equivalence class of $a\in A$.
Of course we have the same equivalence relation on $B$.
The open $L$-neighborhoods of the equivalence clases give an open
cover of $T$ with Lebesgue number $L$ and multiplicity two.
\qua

\begin{prop} \label{asdim_ofR}
Assume $X\subseteq\RRR^n$ contains arbitrarily big balls, i.e.
for each $n\in\NNN$ there is a point $x_n\in X$
such that the ball $B_{\RRR^n}(x_n,n)$ with center $x_n$ and radius $n$ is contained in $X$.
In this case $\asdim(X,d_{eucl.})=n$.
\end{prop}
\textbf{Proof.}
First we prove $\asdim(\RRR^n,d_{eucl.})\leq n$.
Then an application of Remark~\ref{asdim_subset} leads to $\asdim(X,d_{eucl.})\leq n$.

By $Q_a(x) := \left\{y\in\RRR^n \mid \left|y_i-x_i\right|<\frac{a}{2}
\text{ for all } i\in\{1,\ldots,n\} \right\}$
we denote the specified $n$-cube around $x\in\RRR^n$ with edges of length $a>0$.
We define $v := (1,\ldots,1)\in\ZZZ^n$ and consider
the following $n+1$ families of disjoint, open cubes:
$$\mathcal{Q}_i := \left\{ Q_a\left(a\cdot\left(z+\frac{i}{n+1}\cdot v \right)\right) \mid 
z\in\ZZZ^n\right\} \quad \text{ for } i\in\{0,\ldots,n\}$$
The cover $\mathcal{U}_a:=\mathcal{Q}_0 \cup\cdots\cup \mathcal{Q}_n$
is uniformly bounded and open and has multiplicity $\mu(\mathcal{U}_a)=n+1$.
We calculate the Lebesgue number of this cover.
The boundaries of the covering sets decompose $\RRR^n$ into cubes with
edges of length $\frac{a}{n+1}$. Any $n$-cube is limited by $2n$ faces.
Now take $x\in\RRR^n$. Choose the small closed cube $x$ is belonging
to. At most $n$ of its limiting faces are at a distance not exceeding
$\frac{a}{2\cdot(n+1)}$. There is $i\in\{0,\ldots,n\}$ such that a cube of
$\mathcal{Q}_i$ contains all these faces.
Therefore $B_{\RRR^n}\left(x,\frac{a}{2\cdot(n+1)}\right)$
is contained in this cube.
If we want the cover $\mathcal{U}_a$ to have Lebesgue number $L(\mathcal{U}_a)\geq L$,
we just choose $a\geq 2\cdot(n+1)\cdot L$.

Now suppose $\asdim(X,d_{eucl.})=:k<n$.
There is a uniformly bounded, open cover $\mathcal{U}$ of $X$ with multiplicity
$\mu(X)\leq k+1<n+1$. Let $d$ be the mesh of the cover.
Let $\varepsilon>0$. We get an open cover of $B_{\RRR^n}(0,1)$ with
sets of diameter at most $\varepsilon$ and of multiplicity not exceeding $k+1$
by translating a ball of diameter $\frac{d}{\varepsilon}$ contained in $X$
to the origin and multiplying with $\frac{\varepsilon}{d}$.
Applying Theorem~1.6.12 of \cite{Engelking} yields
$\dim\left(B_{\RRR^n}(0,1)\right)\leq k<n$, but this is wrong.
Thus $\asdim(X,d_{eucl.})\geq n$.
\qua

\begin{cor} \ \
$\asdim(\RRR^n,d_{eucl.}) = n$
\end{cor}

\section{Asymptotic dimension of coarse spaces}\label{generalizeasdim}

In \cite{RoeCG}, John Roe generalized Definition~\ref{metric_asdim_2}.\
We are now going to generalize Definition~\ref{metric_asdim}.
\begin{defn}\textbf{(asymptotic dimension of coarse spaces)} \label{asdimGr} \\
Let $(X,\mathcal{E})$ be a coarse space.
\begin{itemize}
\item Let $L\in\mathcal{E}$ be an entourage and $\mathcal{U}$ a cover of $X$.
We say that $\mathcal{U}$ has \emph{appetite}\index{appetite} $L$ if
$\underset{x\in X}{\forall} \ \underset{U\in\mathcal{U}}{\exists} \ L(x)\subseteq U$.
\item We call a cover $\mathcal{U}$ \emph{uniformly bounded}\index{uniformly bounded}
if $\Delta_\mathcal{U} := \bigcup_{U\in\mathcal{U}}U\times U$ is an entourage.
\item Let $n\in\NNN$. We say
$\asdim(X,\mathcal{E})\leq n$\index{asdim}\ if for every\footnote{
    An entourage $L\in\mathcal{E}$ is called symmetric\index{symmetric entourage} if $L=L^{-1}$.
    We need to consider only symmetric entourages which contain the diagonal,
    because for any entourage $L\in\mathcal{E}$
    we have $L\subseteq L\cup L^{-1}\cup \Delta_X\in\mathcal{E}$.}
entourage $L\in\mathcal{E}$ there exists a cover $\mathcal{U}$ of $X$ such that
\begin{enumerate}
\item[(1)] the multiplicity $\mu(\mathcal{U})$ is at most $n+1$,
\item[(2)] $\mathcal{U}$ has appetite $L$ and
\item[(3)] $\mathcal{U}$ is uniformly bounded.
\end{enumerate}
\end{itemize}
\end{defn}

\begin{rem}
If $(X,d)$ is a pseudometric space and $\mathcal{E}_d$ the corresponding bounded coarse structure,
then $\asdim(X,d)$ as defined in Definition~\ref{metric_asdim} and $\asdim(X,\mathcal{E}_d)$ coincide.
\end{rem}

\begin{defn} \label{asdimgroup}
Let $G$ be a finitely generated group.
Lemma~\ref{wordlemma} yields that different word metrics on $G$
all induce the same coarse structure $\mathcal{E}_G$.
We define the asymptotic dimension of the finitely generated group $G$.
$$\asdim(G):=\asdim(G,\mathcal{E}_G)$$
\end{defn}

\begin{prop}
A cover $\mathcal{U}$ is uniformly bounded if and only if
\begin{equation} \label{glUB}
\underset{D\in\mathcal{E}}{\exists} \ \underset{U\in\mathcal{U}}{\forall} \
                       \underset{x\in X}{\exists} \ U\subseteq D(x) \ .
\end{equation}
\end{prop}
\textbf{Proof.}
Let $(X,\mathcal{E})$ be a coarse space and let $\mathcal{U}$ be a cover of $X$.
Defining $D:=\Delta_\mathcal{U}$ we see that (\ref{glUB}) follows from being uniformly bounded.

Suppose conversely that we have an entourage $D$ as in (\ref{glUB}).
Now we conclude that $\bigcup_{U\in\mathcal{U}} U\times U \subseteq DD^{-1}\in\mathcal{E}$,
i.e. $\mathcal{U}$ is uniformly bounded.
\qua

\begin{example}\label{tcs}
Let $X$ be a Hausdorff space. We call the collection
$$\mathcal{T} := \cs\big(\{K\subseteq X\times X\mid K \text{ compact }\}\big)$$
the \emph{trivial coarse structure} on $X$. \index{trivial coarse structure}
Observe that this coarse structure is not compatible with the topology
if $X$ is not compact and points in $X$ are not open.
The asymptotic dimension of $(X,\mathcal{T})$ is zero.
\end{example}
\textbf{Proof.}
The projections $\pi_i\colon X\times X \to X, (x_1,x_2)\mapsto x_i$ are continuous.
Let $L\in\mathcal{E}$ be a symmetric entourage and define $K:=\pi_1(\overline{L\backslash\Delta_X})$.
Then $\mathcal{U}:=\{K\}\cup\{\{x\}\mid x\not\in K\}$ is a uniformly bounded cover with appetite $L$ and multiplicity zero.
\qua

\begin{rem}
Note that the covers in Definition~\ref{asdimGr} do not have to be open.
If we demand the covers to be open, we obtain $\asdim_{open}(X,\mathcal{T}) = \infty$
for every connected, non-compact Hausdorff space $X$, because a cover of $X$ cannot be open and uniformly bounded with respect to $\mathcal{T}$ at the same time.
\end{rem}

\subsection*{Comparison with Roe's definition}

The following is a redraft of Roe's definition of asymptotic dimension.
\begin{defn} \label{asdimRoe}\index{asdim}
Let $(X,\mathcal{E})$ be a coarse space.
We say $\asdim_{Roe}(X,\mathcal{E})\leq n$ if for every entourage $L\in\mathcal{E}$
there is a cover $\mathcal{U}$ of $X$ such that
\begin{enumerate}
\item[(1)] \ $\mathcal{U} = \mathcal{U}_1\ \cup\ \cdots\ \cup\ \mathcal{U}_{n+1}$,
\item[(2)] \ each of the families $\mathcal{U}_1,\ldots,\mathcal{U}_{n+1}$ is $L$-disjoint\index{disjoint@$L$-disjoint}\\
(i.e. whenever $A,B\in\mathcal{U}_i$ and $A\neq B$, then ${A\times B}\ \cap\ L=\emptyset$) and
\item[(3)] \ $\mathcal{U}$ is uniformly bounded.
\end{enumerate}
\end{defn}

There is a small difference between Definition~\ref{asdimRoe} and the definition given in \cite{RoeCG}.
In Roe's original definition the cover $\mathcal{U}$ is supposed to be countable.
We will not make any assumptions on the cardinality of $\mathcal{U}$.

\begin{rem}
Let $(X,d)$ be a pseudometric space and $\mathcal{E}$ the corresponding bounded coarse structure.
It is easy to see that $\asdim_{Roe}(X,\mathcal{E})$ and the asymptotic dimension of $(X,d)$ as defined in Definition~\ref{metric_asdim_2} coincide.
\end{rem}

A third version of asymptotic dimension will appear in Theorem~\ref{asdimgleich}.
\begin{defn} \label{asdimfam}\index{asdim}
$\asdim_{fam}(X,\mathcal{E})\leq n$ if for every entourage $L\in\mathcal{E}$
there is a cover $\mathcal{U}$ of $X$ such that
(1) $\mathcal{U} = \mathcal{U}_1\ \cup\ \cdots\ \cup\ \mathcal{U}_{n+1}$
where each of the families $\mathcal{U}_i$ consists of disjoint sets,
(2) $\mathcal{U}$ has appetite $L$ and
(3) $\mathcal{U}$ is uniformly bounded.
\end{defn}

\begin{thm} \label{asdimgleich}
Let $(X,\mathcal{E})$ be a coarse space. Then
$$\asdim(X,\mathcal{E}) = \asdim_{Roe}(X,\mathcal{E})
    = \asdim_{fam}(X,\mathcal{E})\ .$$ 
\end{thm}
\textbf{Proof.}
We first prove $\asdim_{Roe}\geq\asdim_{fam}$.
Assume $\asdim_{Roe}(X,\mathcal{E})=n\in\NNN$.
Let $L$ be a symmetric entourage which contains the diagonal.
For $L^2:=LL\in\mathcal{E}$, there exists a cover $\mathcal{U}$ as in Definition~\ref{asdimRoe}.
Since $A\times B\ \cap\ L^2=\emptyset$ is equivalent to $L[A]\cap L[B]=\emptyset$, the cover
$\mathcal{U}_L:=\{L[U]\mid U\in\mathcal{U}\}$ meets all conditions required in Definition~\ref{asdimfam}.
Note that
$\bigcup_{U\in\mathcal{U}}L[U]\times L[U] \subseteq L\left(\bigcup_{U\in\mathcal{U}}U\times U\right)L^{-1}\in\mathcal{E}$.

In a second step we have to prove $\asdim_{fam}\geq\asdim$,
but this is obvious, since condition (1) of Definition~\ref{asdimfam}
implies condition (1) of Definition~\ref{asdimGr}.

It remains to prove $\asdim\geq\asdim_{Roe}$.
For this purpose we need to construct a uniformly bounded cover $\mathcal{V}$
consisting of $L$-disjoint families from a uniformly bounded cover $\mathcal{U}$
with appetite $L^{n+1}$.
The idea is to take all intersections of $n+1$ sets from $\mathcal{U}$ as one family,
the intersections of exactly $n$ sets from $\mathcal{U}$ as a second family, etc.
However, we have to ensure these families to be $L$-disjoint.

Assume that $\asdim(X,\mathcal{E})=n\in\NNN$.
Let $L\in\mathcal{E}$ be a symmetric entourage that contains the diagonal.
Let $\mathcal{U}$ be a uniformly bounded cover of $X$ with appetite $L^{n+1}$ and multiplicity at most $n+1$.
For an entourage $E$ and $U\subseteq X$ we define $Int_E(U):=\{x\in X\mid E(x)\subseteq U\}$.
Observe that $E_1\subseteq E_2$ implies $Int_{E_2}(U)\subseteq Int_{E_1}(U)$.
Some more definitions are needed to get $\mathcal{V}$.
\begin{eqnarray*}
\mathcal{U}_i &:=& \{U_1\cap\cdots\cap U_i\mid U_1,\ldots,U_i\in\mathcal{U}\text{ pairwise distinct }\} \\
S_i &:=& \bigcup_{U\in\mathcal{U}_i}Int_{L^{n+2-i}}(U) \\
\mathcal{V}_i &:=& \left\{Int_{L^{n+2-i}}(U)\backslash S_{i+1}\mid U\in\mathcal{U}_i\right\} \\
\mathcal{V} &:=& \mathcal{V}_1 \cup\cdots\cup\mathcal{V}_{n+1}
\end{eqnarray*}
Now $\mathcal{V}$ is a cover of $X$. Actually, $\mathcal{V}$ is a refinement of the cover $\mathcal{U}$.
Therefore $\mathcal{V}$ is uniformly bounded.

It remains to prove that each of the families $\mathcal{V}_1,\ldots,\mathcal{V}_{n+1}$ is $L$-disjoint.
Let $A,B\in\mathcal{V}_i$ such that $A\neq B$.
There are $A_1,\ldots,A_i,B_1,\ldots,B_i\in\mathcal{U}$ such that
$A=Int_{L^{n+2-i}}(A_1\cap\cdots\cap A_i)\backslash S_{i+1}$ and
$B=Int_{L^{n+2-i}}(B_1\cap\cdots\cap B_i)\backslash S_{i+1}$.
The sets $A_1,\ldots,A_i$ are supposed to be pairwise distinct as are the sets $B_1,\ldots,B_i$.

Let $(a,b)\in A\times B\ \cap\ L$ and observe the following facts.
\begin{eqnarray}
a,b & \not\in & S_{i+1} \label{glSi1} \\
a\ \in\ A &\subseteq & Int_{L^{n+1-i}}(A_1\cap\cdots\cap A_i) \label{glA} \\
b\ \in\ B &\subseteq & Int_{L^{n+1-i}}(B_1\cap\cdots\cap B_i) \label{glB} \\
a\ \in\ L[B] &\subseteq & L[Int_{L^{n+2-i}}(B_1\cap\cdots\cap B_i)] \label{glLB} \\
b\ \in\ L[A] &\subseteq & L[Int_{L^{n+2-i}}(A_1\cap\cdots\cap A_i)] \label{glLA}
\end{eqnarray}
Since
$L[Int_{L^j}(U)]
 = \left\{x\mid \exists_{y\in X}L^j(y)\subseteq U,x\in L(y)\right\}
\subseteq \{ x\mid L^{j-1}(x)\subseteq U\} = Int_{L^{j-1}}(U)$,
we get the following conclusions from (\ref{glLB}) and (\ref{glLA}):
\begin{eqnarray*}
a &\in & Int_{L^{n+1-i}}(B_1\cap\cdots\cap B_i) \\
b &\in & Int_{L^{n+1-i}}(A_1\cap\cdots\cap A_i)
\end{eqnarray*}
Finally $a,b\in Int_{L^{n+2-(i+1)}}(A_1\cap\cdots\cap A_i\cap B_1\cap\cdots\cap B_i)$. Since $A\neq B$,
we know that the set $\{A_1,\ldots,A_i,B_1,\ldots,B_i\}$ contains at least $i+1$ elements.
Thus $a,b\in S_{i+1}$, but this is a contradiction to (\ref{glSi1}).
\qua

\begin{rem} \label{LebesguezahlBuchfuehrung}
For further reference we record the results of our constructions in the previous proof.
Let $(X,\mathcal{E})$ be a coarse space, $L\in\mathcal{E}$ and $n\in\NNN$.
\begin{itemize}
\item If we are given a uniformly bounded cover $\mathcal{U}$ of $X$ such that
   $\mathcal{U}=\mathcal{U}_1\cup\cdots\cup\mathcal{U}_{n+1}$ and each of the families
   $\mathcal{U}_1,\ldots,\mathcal{U}_{n+1}$ is $L^2$-disjoint,
   we get the uniformly bounded cover $\mathcal{U}_L:=\{L[U]\mid U\in\mathcal{U}\}$
   consisting of $n+1$ families of disjoint subsets of $X$
   and with appetite $L$.
\item Given a uniformly bounded cover $\mathcal{U}$ of $X$ with multiplicity $n+1$
   and appetite $L^{n+1}$, we constructed a uniformly bounded cover
   $\mathcal{V}$ such that $\mathcal{V}=\mathcal{V}_1\cup\cdots\cup\mathcal{V}_{n+1}$
   and each of the families $\mathcal{V}_1,\ldots,\mathcal{V}_{n+1}$ is $L$-disjoint.
\end{itemize}
Suppose $X$ is a metric space and $\mathcal{E}$ its bounded coarse structure.
Given a uniformly bounded cover with multiplicity $m$ and Lebesgue number $L$, the previous construction yields
a uniformly bounded cover which consists of $m$ families of disjoint subsets of $X$
and with Lebesgue number $\frac{L}{2n+2}$.
\end{rem}

\subsection*{Asymptotic dimension via anti-\v Cech systems}

The asymptotic dimension of proper metric spaces can be characterized
using anti-\v Cech systems. Compare Theoreom 9.9 of \cite{RoeCG}.
We will see that this characterisation does not work for all coarse spaces.
However, it does work for coarse structures induced from metrizable compactifications.

We recall some definitions from John Roe's book.
\begin{defn} Let $(X,\mathcal{E})$ be a coarse space.
\begin{itemize}
\item
A uniformly bounded cover $\mathcal{U}$ of $X$ is called \emph{uniform}\index{uniform cover}
if each bounded subset of $X$ meets only finitely many elements of $\mathcal{U}$.
\item
A collection of uniform covers is called an anti-\v Cech system\index{anti-\v Cech system}
if for every entourage $L\in\mathcal{E}$ it contains a cover having appetite $L$.
\item
We define a partial order on an anti-\v Cech system.
$\mathcal{U}\leq\mathcal{V}$ if $\mathcal{V}$ has appetite $\Delta_\mathcal{U}$.
This way, any anti-\v Cech system becomes a directed set.
\end{itemize}
\end{defn}

\begin{example}
Consider the coproduct $\coprod_{n\in\NNN}(\RRR,\mathcal{E}_{eucl.})$.
According to Pro\-position~\ref{asdiminfcoprod}, this coarse space has asymptotic dimension one.
On the other hand, for every uniformly bounded cover there exists a copy of $\RRR$
where the cover consists of the sets $\{x\}$ for $x\in\RRR$.
Hence, there don't exist uniform covers for this coarse space
and the following proposition can not be true in general.\footnote{
Even if this coarse space does not admit an anti-\v Cech system,
it does have a coarsening sequence in the sense of \cite{MitchAdd},
i.e. this space belongs to the category where coarse homology theory can be defined.}
\end{example}

\begin{prop} \label{asdimanticech}
Let $X$ be a Hausdorff space and let $K$ be a metrisable compactification of $X$.\footnote{
This implies that $X$ is locally compact and paracompact.
Note that a compact set $K$ is metrisable if and only if its topology has a countable base.}
Let $\mathcal{E}$ be the coarse structure induced by $K$.
Then $\asdim(X,\mathcal{E})\leq n$ if and only if there is an anti-\v Cech system
consisting of covers with multiplicity at most $n+1$.
\end{prop}
\begin{rem} \label{asdimanticechrem}
If we take just uniformly bounded covers in the definition of anti-\v Cech systems,
Proposition~\ref{asdimanticech} is true for all coarse spaces - just by definitions.
\end{rem}

\textbf{Proof of Proposition~\ref{asdimanticech}.}
Because of Remark~\ref{asdimanticechrem}, we only need to prove that we can choose all covers to be uniform.

Let us call a cover $\mathcal{V}$ a nice refinement of a cover $\mathcal{U}$ if we have an injective map
$i\colon\mathcal{V}\to\mathcal{U}$ such that $V\subseteq i(V)$ for all $V\in\mathcal{V}$.

Let $L$ be a symmetric entourage that contains the diagonal and 
let $\mathcal{U}$ be a uniformly bounded cover of $X$ consisting of $n+1$ families each of them being $L^2$-disjoint.
Note that every nice refinement of $\mathcal{U}$ has the same properties.

Since $K$ is metrisable, we can write $X=\bigcup_{i\in\NNN}K_i$
with $K_i$ compact and $K_1\subseteq K_2\subseteq \cdots\subseteq X$.
We proceed inductively in order to define a nice refinement of $\mathcal{U}$
which finally will lead to a uniform cover of appetite $L$.
Set $\mathcal{U}_0 := \mathcal{U}$. 
Observe that $L[K_i]$ can be covered by a finite subcover $\mathcal{U}_{fin,i}$ of $\mathcal{U}_{i-1}$.
Define $$\mathcal{U}_i:= \mathcal{U}_{fin,i}\ \cup\
\{\ U\backslash L[K_i]\ \mid\ U\in\mathcal{U}_{i-1}\backslash\mathcal{U}_{fin,i}\ \}\ .$$
The cover $\mathcal{U}_i$ is a nice refinement of $\mathcal{U}_{i-1}$.
Finally, we get a nice refinement
$$\mathcal{U}_{red}:= \left\{ U \mid \exists_{k\in\NNN} \forall_{i\geq k} U\in\mathcal{U}_{fin,i} \right\}$$
of the original cover $\mathcal{U}$.
In order to see that $\mathcal{U}_{red}$ is indeed a cover of $X$, let $x\in K_i\subseteq X$.
The union of all $U\in\mathcal{U}_{fin,i}$ containing $x$ is bounded
and therefore contained in $K_j$ for some $j\in\NNN$.
It follows that any $U\in\mathcal{U}_{fin,j}$ which contains $x$ is an element of $\mathcal{U}_{red}$.

Observe that $\mathcal{V}_L := \{L[U]\mid U\in\mathcal{U}_{red}\}$ is
a uniformly bounded cover of $X$ with appetite $L$ and multiplicity at most $n+1$.
Moreover $\mathcal{V}_L$ is uniform.
In order to see this, let $B\subseteq X$ be bounded. This implies $B\subseteq K_i$ for some $i\in\NNN$.
Hence $B \cap U \neq \emptyset$ for most $N_i$ different $U\in\mathcal{U}_{red}$
where $N_i$ is the cardinality of $\mathcal{U}_{fin,i}$. 

Thus, the covers $\mathcal{V}_L$ for symmetric entourages $L$ which contain the diagonal
form the desired anti-\v Cech system.
\qua

\section[Basic properties]{Basic properties of asymptotic dimension}

\begin{thm} \label{thmasdimemb}\index{coarse embedding}
If $f\colon(X,\mathcal{E})\to(Y,\mathcal{F})$ is a coarse embedding, then
$$\asdim(X,\mathcal{E}) \leq \asdim(Y,\mathcal{F}).$$
\end{thm}
\textbf{Proof.}
Suppose that $n:=\asdim(Y,\mathcal{F})<\infty$.
Let $E\in\mathcal{E}$ be an entourage and set $F:=f\times f(E)$. Note that $E\subseteq (f\times f)^{-1}(F)$.
There is a uniformly bounded cover $\mathcal{U}$ of $Y$ with appetite $F$ and multiplicity at most $n+1$.
The inverse image of $\mathcal{U}$ is a uniformly bounded cover of $X$ with appetite $E$
and the same multiplicity as $\mathcal{U}$. 
\qua

\begin{cor}[monotony of asymptotic dimension] \label{monotony}\index{monotony}
Let $(X,\mathcal{E})$ be \linebreak a coarse space and $A\subseteq X$.
Remember that the inclusion map is a coarse embedding.
Hence $\asdim(A,\mathcal{E}|_A)\leq\asdim(X,\mathcal{E})$.
\end{cor}

\begin{cor}\label{asdimci}\index{coarse invariance}
\textbf{(coarse invariance of asymptotic dimension)} \\
Asymptotic dimension is a coarse invariant, i.e. given a coarse equivalence
$f\colon(X,\mathcal{E}_X)\to(Y,\mathcal{E}_Y)$,
we have $\asdim(X,\mathcal{E}_X)=\asdim(Y,\mathcal{E}_Y)$.
\end{cor}

\begin{defn}
Let $(X,\mathcal{E})$ be a coarse space and $A\subseteq X$.
We call $A$ a \emph{substantial part}\index{substantial part}
of $X$ if $\asdim(A,\mathcal{E}|_A)=\asdim(X,\mathcal{E})$.
\end{defn}

For a coarsely uniform map $(X,\mathcal{E})\to (Y,\mathcal{F})$ which is also injective,
there is no relation between the asymptotic dimensions of $(X,\mathcal{E})$
and $(Y,\mathcal{F})$.

\begin{example}
Let $(X,\mathcal{E})$ be a coarse space.
Observe that the power set $\mathcal{P}(X\!\times\! X)$ is a coarse structure on $X$.
The map $id\colon (X,\mathcal{E})\to (X,\mathcal{P}(X\!\times\! X))$ is coarsely uniform,
but $\asdim(X,\mathcal{E}) \geq 0 = \asdim(X,\mathcal{P}(X\!\times\! X))$.
\end{example}

\begin{example} \label{examplern}\index{one-point compactification}
Let $n\in\{2,3,\ldots\}$.
By $\mathcal{E}_\cdot$ we denote the coarse structure coming from the one-point compactification of $\RRR^n$
and by $\mathcal{E}_{eucl.}$ the bounded coarse structure corresponding to the euclidean metric of $\RRR^n$.
The map $id \colon (\RRR^n,\mathcal{E}_{eucl.}) \to (\RRR^n,\mathcal{E}_\cdot)$ is coarse, but
$\asdim(\RRR^n,\mathcal{E}_{eucl.}) = n > 1 = \asdim(\RRR^n,\mathcal{E}_\cdot)$.
\end{example}
\textbf{Proof.}
Since $\mathcal{E}_{eucl.}\subseteq\mathcal{E}_\cdot$, the map $id$ is coarsely
uniform. A set $B$ is bounded with respect to $\mathcal{E}_{eucl.}$ if and
only if $B$ is precompact. The same is true for $\mathcal{E}_\cdot$.
Thus $id$ is coarsely proper.

It remains to prove $\asdim(\RRR^n,\mathcal{E}_\cdot)=1$.
Note that
\begin{eqnarray*}
\abs{\ \cdot\ }\ \colon(\RRR^n,\mathcal{E}_\cdot) & \to & (\RRR_+,\mathcal{E}_\cdot) \\
x & \mapsto & \abs{x}
\end{eqnarray*}
is a coarse equivalence with inverse
$ \imath\colon (\RRR_+,\mathcal{E}_\cdot) \to (\RRR^n,\mathcal{E}_\cdot)\ ,\ r \mapsto (r,0,\ldots,0)$.
In fact $\abs{\ \cdot \ }\circ \imath = id_{\RRR^n}$
and $\imath\circ\abs{\ \cdot\ }$ is close to $id_{\RRR^n}$ with respect to $\mathcal{E}_\cdot$.
This yields $\asdim(\RRR^n,\mathcal{E}_\cdot)=\asdim(\RRR_+,\mathcal{E}_\cdot)$.
It remains to calculate the asymptotic dimension of $(\RRR_+,\mathcal{E}_\cdot)$.
Since $(\RRR_+,\mathcal{E}_\cdot)$is a ray, we can apply Corollary~\ref{asdimcoarsecells}.

As an alternative, we refer to Example 9.7 of \cite{RoeCG}.
\qua

\begin{prop}\label{cnegal}
Let $(X,\mathcal{E})$ be a coarse space and $\mathcal{E}_{cn} = \cncs(\mathcal{E})$
the connected coarse structure generated by $\mathcal{E}$.
Then $\asdim(X,\mathcal{E}) = \asdim(X,\mathcal{E}_{cn})$.
\end{prop}
\textbf{Proof.}
Suppose $n:=\asdim(X,\mathcal{E}_{cn})<\infty$.
Let $E\in\mathcal{E}\subseteq\mathcal{E}_{cn}$.
There is a cover $\mathcal{U}$ of $X$  with appetite $E$
and multiplicity at most $n+1$ which is uniformly bounded with respect to $\mathcal{E}_{cn}$.
Each $U\in\mathcal{U}$ can be written as the disjoint union of finitely many sets $U_1,\ldots,U_k$
which are bounded with respect to $\mathcal{E}$ and such that the union of any two of the sets $U_1,\ldots,U_k$
is not bounded with respect to $\mathcal{E}$.
Define $\operatorname{comp}(U) := \{U_1,\ldots,U_k\}$
and observe that $\mathcal{U}' := \bigcup_{U\in\mathcal{U}} \operatorname{comp}(U)$
is a cover of $X$ with multiplicity at most $n+1$.
Furthermore, $\mathcal{U}'$ has appetite $E$ and is uniformly bounded with respect to $\mathcal{E}$.
Hence $\asdim(X,\mathcal{E})\leq n$.

Set $n:=\asdim(X,\mathcal{E})$.
Let $E\in\mathcal{E}_{cn}$ be a symmetric entourage. This implies that
$E = E'\cup\left(A_1\times A_{\sigma(1)}\right)\cup\cdots\cup\left(A_k\times A_{\sigma(k)}\right)$
with $E'\in\mathcal{E}$, $k\in\NNN$, $A_1,\ldots,A_k$ bounded subsets of $X$ (not necessarily pairwise distinct)
and $\sigma$ a permutation of the set $\{1,\ldots,k\}$ with the additional property $\sigma\circ\sigma = \id$.
Set $M:=\Delta_X\cup A_1^2\cup\cdots\cup A_k^2$
and observe that $E'' := (E'\cup\Delta_X)M\in\mathcal{E}$.
Let $\mathcal{U}''$ be a cover of $X$ which is uniformly bounded with respect to $\mathcal{E}$
and which has appetite $E''$ and multiplicity at most $n+1$.
Note that there are sets $U_1,\ldots,U_k\in\mathcal{U}''$ such that $E'[A_i]\cup A_i\subseteq U_i$.
We define the cover $\mathcal{U} := \mathcal{U}'' \cup \{U_1\cup\cdots\cup U_k\} \backslash \{U_1,\ldots,U_k\}$
of $X$. Observe that $\mathcal{U}$ is uniformly bounded with respect to $\mathcal{E}_{cn}$
and has multiplicity at most $n+1$.

Moreover, $\mathcal{U}$ has appetite $E$.
To see this, let $x\in X$.
If $x\in A_i$, then $E(x) = E'(x) \cup A_{\sigma(i)} \subseteq U_i\cup U_{\sigma(i)}
\subseteq U_1\cup\cdots\cup U_k$.
If $x\not\in\{A_1\cup\cdots\cup A_k\}$, then $E(x) = E'(x) \subseteq E''(x)$.
\qua

\begin{prop}\textbf{(asymptotic dimension of finite unions)}\index{union of coarse spaces}\label{asdimfinunions}\\
Let $(X,\mathcal{E})$ be a coarse space and $A,B\subseteq X$ with $A\cup B=X$.
$$\asdim(X,\mathcal{E}) = \max \{\asdim(A,\mathcal{E}|_A),\asdim(B,\mathcal{E}|_B)\}$$
\end{prop}
\textbf{Proof.}
The proof of $\ \geq\ $ follows from monotony.
To see $\ \leq\ $, we generalize an argument of Bell and Dranishnikov (see \cite{BD1}).

Let $n$ be the maximum of $\asdim(A,\mathcal{E}|_A)$ and $\asdim(B,\mathcal{E}|_B)$
and take a symmetric entourage $L\in\mathcal{E}$ which contains $\Delta_X$.
For $\mathcal{U}\subseteq\mathcal{P}(X)$ and $V\subseteq X$ we define
$$N_L(V,\mathcal{U})\ :=\ V\ \cup \bigcup_{\substack{U\in\mathcal{U}\\ L\cap U\times V\neq\emptyset}} U\ .$$
There is a uniformly bounded cover $\mathcal{U}=\mathcal{U}_1\cup\cdots\cup\mathcal{U}_{n+1}$ of $A$
consisting of $L$-disjoint families $\mathcal{U}_i$.
Moreover, there is a uniformly bounded cover $\mathcal{V}=\mathcal{V}_1\cup\cdots\cup\mathcal{V}_{n+1}$ of $B$
consisting of $(L\Delta_\mathcal{U}L\Delta_\mathcal{U}L)$-disjoint families $\mathcal{V}_i$. For $i\in\{1,\ldots,n+1\}$ set
$$\mathcal{W}_i := \big\{ N_L(V,\mathcal{U}_i) \mid V\in\mathcal{V}_i \big\} \cup 
                   \big\{ U\in\mathcal{U}_i \mid L\cap U\times V = \emptyset \text{ for all } V\in\mathcal{V}_i \big\}.$$
Observe that $N_L(V,\mathcal{U}_i) \subseteq \Delta_\mathcal{U}L[V]$.
Hence, we get a uniformly bounded cover $\mathcal{W}=\mathcal{W}_1\cup\cdots\cup\mathcal{W}_{n+1}$ of $X$
where $\mathcal{W}_i$ is $L$-disjoint for $1\leq i\leq n+1$. This proves $\asdim(X,\mathcal{E})\leq n$.
\qua

\begin{prop}\label{asdiminfcoprod}\index{coproduct}\textbf{(asymptotic dimension of coproducts)}\\
Let $\Lambda$ be any set and $(X_\lambda,\mathcal{E}_\lambda)$ a coarse space for every $\lambda\in\Lambda$.
Define $X:=\underset{\lambda\in\Lambda}{\coprod} X_\lambda$.
If $(X,\mathcal{E})$ and $(X,\mathcal{E}_{cn})$ are the coproducts
in the categories $\mathcal{D}$ and $\mathcal{D}_{cn}$ respectively, then
$$\asdim\left(X,\mathcal{E}_{cn}\right) = \asdim\left(X,\mathcal{E}\right)
= \underset{\lambda\in\Lambda}{\sup} \asdim(X_\lambda,\mathcal{E}_\lambda) \ .$$
\end{prop}
\textbf{Proof.}
Set $n:=\underset{\lambda\in\Lambda}{\sup} \asdim(X_\lambda,\mathcal{E}_\lambda)$.
Monotony of asymptotic dimension implies $\asdim(X,\mathcal{E}_{cn}) \geq n$
and $\asdim(X,\mathcal{E}) \geq n$.

We will now prove $\asdim(X,\mathcal{E}) \leq n$.
Take an entourage $L\in\mathcal{E}$ which contains $\Delta_X$.
Then there are $\lambda_1,\ldots,\lambda_k\in\Lambda$ and $L_{\lambda_i}\in\mathcal{E}_{\lambda_i}$
such that $L=L_{\lambda_1}\cup\cdots\cup L_{\lambda_k}\cup\Delta_X$.
For $i\in\{1,\ldots,k\}$ choose a uniformly bounded cover $\mathcal{U}_{\lambda_i}$ of $X_{\lambda_i}$
with appetite $L_{\lambda_i}$ and multiplicity at most $\asdim(X_{\lambda_i},\mathcal{E}_{\lambda_i})+1$.
For $\lambda\in\Lambda\backslash\{\lambda_1,\ldots,\lambda_k\}$
set $\mathcal{U}_\lambda:=\{\{x\} \mid x\in X_\lambda\}$.
The union $\mathcal{U}:=\bigcup_{\lambda\in\Lambda}\mathcal{U}_\lambda$ is a uniformly bounded cover of $X$
with appetite $L$ whose multiplicity does not exceed $n+1$.

For the equality $\asdim(X,\mathcal{E}_{cn})=\asdim(X,\mathcal{E})$ compare Proposition~\ref{cnegal}.
\qua

\begin{prop}\textbf{(asymptotic dimension of disjoint unions)}\index{disjoint union}\\
Let $\Lambda$ be any set and $(X_\lambda,\mathcal{E}_\lambda)$ a coarse space for every $\lambda\in\Lambda$.
$$\asdim\left(\underset{\lambda\in\Lambda}{\bigsqcup} X_\lambda,\mathcal{E}_{\bigsqcup}\right)
= \underset{\lambda\in\Lambda}{\sup} \asdim(X_\lambda,\mathcal{E}_\lambda)$$
\end{prop}
\textbf{Proof.}
The proof of $\ \geq\ $ follows from monotony of asymptotic dimension.

Set $X := \underset{\lambda\in\Lambda}{\bigsqcup} X_\lambda$
and $n := \underset{\lambda\in\Lambda}{\sup} \asdim(X_\lambda,\mathcal{E}_\lambda)$.
In order to prove the inequality $\asdim(X,\mathcal{E}_{\bigsqcup})\leq n$, let $L\in\mathcal{E}_{\bigsqcup}$.
Set $L_\lambda := L\cap X_\lambda^2 \in \mathcal{E}_\lambda$ and
choose uniformly bounded covers $\mathcal{U}_\lambda$ of $X_\lambda$ with appetite $L_\lambda$
and multiplicity at most $n+1$.
Now $\mathcal{U} := \bigcup_{\lambda\in\Lambda}\mathcal{U}_\lambda$ is a uniformly bounded cover of $X$
with appetite $L$ and multiplicity at most $n+1$.
\qua

\begin{prop}[asdim of certain direct limits]\label{asdimcolim}\index{direct limit}
Let $j\in I$. If $I$ is a directed set and
if for all $k\geq j$ the map $f_k\colon X_k\to\underrightarrow{\lim}X_i$ is injective, then
$$\asdim\left(\underrightarrow{\lim}(X_i,\mathcal{E}_i)\right) = \sup_{k\geq j}\asdim(X_k,\mathcal{E}_k)\ .$$
\end{prop}
\textbf{Proof.}
With Proposition~\ref{limpor} in mind, the proof is not difficult.
Observe that injectivity of $f_k$ for all $k\geq j$ implies that $f_i$ is a coarse embedding for $i\geq j$.
The proof of $\ \geq\ $ now is given by applying Theorem~\ref{thmasdimemb}.

We write $(X,\mathcal{E})$ for the direct limit.
Let $L\in\mathcal{E}$ be an entourage which contains $\Delta_X$.
Proposition~\ref{limpor} implies that $L$ is the union of $\Delta_X$
and an entourage $L_i\in\mathcal{E}_i$ for some $i\in I$.
We may assume $i\geq j$.
If $\mathcal{U}_i$ is a uniformly bounded cover of $X_i$ with appetite $L_i$
and multiplicity at most $\asdim(X_i,\mathcal{E}_i)+1$,
then $\mathcal{U} := \mathcal{U}_i \cup \{ \{x\} \mid x\in X\backslash X_i \}$
is a uniformly bounded cover of $X$ with appetite $L$.
Moreover $\mathcal{U}$ has the same multiplicity as $\mathcal{U}_i$.
\qua

\begin{prop}\textbf{(asymptotic dimension of products)}\\
Let $(X,\mathcal{E}_X)$ and $(Y,\mathcal{E}_Y)$ be coarse spaces.
\begin{eqnarray*} \index{product coarse structure}
\asdim(X\times Y,\mathcal{E}_X * \mathcal{E}_Y) & \leq & \asdim(X,\mathcal{E}_X) + \asdim (Y,\mathcal{E}_Y) \\
\asdim(X,\mathcal{E}_X) & \leq & \asdim(X\times Y,\mathcal{E}_X * \mathcal{E}_Y) \quad \text{ if } Y\neq\emptyset
\end{eqnarray*}
\end{prop}
\textbf{Proof.}
Compare \cite{RoeCG} for the special case of bounded coarse structures and Remark~\ref{asdimprodthmsc}
for the case of continuously controlled coarse structures induced by metrisable compactifications.

Set $n:=\asdim(X)$ and $m:=\asdim(Y)$. Let $E\in\mathcal{E}_X * \mathcal{E}_Y$.
There are symmetric entourages $E_X\in\mathcal{E}_X$ and $E_Y\in\mathcal{E}_Y$ containing the diagonals
$\Delta_X$ and $\Delta_Y$ respectively such that $E\subseteq E_X\times E_Y$.

There is a uniformly bounded cover $\mathcal{U}$ of $X$  with appetite $E_X^{n+m+1}$
and multiplicity $\mu(\mathcal{U})\leq n+1$.
There is also a uniformly bounded cover $\mathcal{V}$ of $Y$ with appetite $E_Y^{n+m+1}$
and multiplicity $\mu(\mathcal{V})\leq m+1$.
We get a uniformly bounded cover $\mathcal{U}\times\mathcal{V}:=\{U\times V\mid U\in\mathcal{U},V\in\mathcal{V}\}$
of $X\times Y$ with appetite $E^{n+m+1}$ and multiplicity $\leq (n+1)\cdot (m+1)=n\cdot m+n+m+1$.
Thus, we need to improve the multiplicity.

We proceed similar as in the proof of Theorem~\ref{asdimgleich}
and begin with some definitions. Let $k\in\{2,\ldots,n+m+2\}$.
\begin{eqnarray*}
\mathcal{A}_k & := & \{U_1\cap\cdots\cap U_p\times V_1\cap\cdots\cap V_q \mid \\
&& \qquad p+q=k, U_i\in\mathcal{U}, V_i\in\mathcal{V} \text{ pairwise distinct }\} \\
B_k & := & \bigcup_{A\in\mathcal{A}_k}Int_{E^{n+m+3-k}}(A) \quad \text{ and } \quad B_{n+m+3}:=\emptyset \\
\mathcal{W}_k & := & \{Int_{E^{n+m+3-k}}(U)\, \backslash\, B_{k+1} \mid U\in\mathcal{A}_k \} \\
\mathcal{W} & := & \mathcal{W}_2\cup\cdots\cup\mathcal{W}_{n+m+2}
\end{eqnarray*}
Notice that $\mathcal{W}$ is a uniformly bounded cover of $X\times Y$
consisting of the $n+m+1$ disjoint families $\mathcal{W}_2,\ldots,\mathcal{W}_{n+m+2}$.
It remains to prove that $\mathcal{W}_k$ is $E$-disjoint for $k\in\{2,\ldots,n+m+2\}$.

For this purpose let $M, N\in\mathcal{W}_k$ with $M\neq N$ and suppose
$M\!\times\! N\, \cap\, E\neq\emptyset$.
Choose $((x_M,y_M),(x_N,y_N))\in M\!\times\! N\, \cap\, E$.
There are $p_M,q_M\in\NNN$ with $p_M+q_M=k$ and
$M_1,\ldots,M_{p_M}\in\mathcal{U}$, $M'_1,\ldots,M'_{q_M}\in\mathcal{V}$ such that
$$M=Int_{E^{n+m+3-k}}(\underbrace{M_1\cap\cdots\cap M_{p_M}\times M'_1\cap\cdots\cap M'_{q_M}}_{=:M\! M'})
\, \backslash\, B_{k+1}.$$
Similarly, there are $p_N,q_N\in\NNN$ with $p_N+q_N=k$ and
$N_1,\ldots,N_{p_N}\in\mathcal{U}$, $N'_1,\ldots,N'_{q_N}\in\mathcal{V}$ such that
$$N=Int_{E^{n+m+3-k}}(\underbrace{N_1\cap\cdots\cap N_{p_N}\times N'_1\cap\cdots\cap N'_{q_N}}_{=:N\! N'})
\, \backslash\, B_{k+1}.$$
Observe that the following relations hold.
\begin{eqnarray*}
(x_M,y_M) & \!\!\not\in\!\! & B_{k+1} \\
(x_M,y_M) & \!\! \in\!\! & M \ \subseteq \ Int_{E^{n+m+2-k}}(M\! M') \\
(x_M,y_M) & \!\! \in\!\! & E[N] \ \subseteq \ E[Int_{E^{n+m+3-k}}(N\! N')] \ \subseteq \ Int_{E^{n+m+2-k}}(N\! N') 
\end{eqnarray*}
It follows that $(x_M,y_M)\in Int_{E^{n+m+3-(k+1)}}(M\! M' \cap N\! N')$.
Since $M\neq N$, the set
$\{M_1,\ldots,M_{p_M},M'_1,\ldots,M'_{q_M},N_1,\ldots,N_{p_N},N'_1,\ldots,N'_{q_N}\}$
contains at least $k+1$ different elements.
Hence $(x_M,y_M)\in B_{k+1}$. But this is a contradiction to what we found before.
\qua

\bigskip

The equality $\asdim(X\times Y) = \asdim(X) + \asdim(Y)$ is not true in general.
Compare Corollary~\ref{asdimprodnoteq}.

\bigskip

There are no general restrictions on the asymptotic dimensions of quotients
\index{quotione coarse structure} as the following three examples show.

\begin{example}
Let $X=\RRR^2$ be the real plane and $\mathcal{E}$ the bounded coarse structure
induced by the euclidean metric on $\RRR^2$.
Define an equivalence relation by $(x,y)\sim (x',y')$ if and only if $x=x'$.
The quotient $X/\!\!\sim$ is clearly $\RRR$ and $\mathcal{E}_\sim$
is the bounded coarse strucure on $\RRR$ again induced by the euclidean metric on $\RRR$.
In this case $\asdim(X,\mathcal{E}) > \asdim(X/\!\!\sim,\mathcal{E}_\sim)$.
\end{example}

\begin{example}
Let $\mathbb{F}_2=\langle a,b\mid\rangle$ be the free group with two generators.
We say that two words are equivalent if the number of letters $a$
counted with exponents as signs is equal in both words and the same for the letter $b$.
The quotient is the group $\ZZZ^2$.
In this case $\asdim(\mathbb{F}_2)=1 < 2=\asdim(\ZZZ^2)$.
\end{example}

\begin{rem}
Let $(X,\mathcal{E})$ be a coarse space. Assume there is an equivalence relation $\sim$ on $X$.
If the equivalence relation happens to be an entourage,
i.e. $\{ (x,y) \mid x\sim y \}\in\mathcal{E}$,
then the projection onto the quotient is a coarse equivalence.
\end{rem}

Let $(X,\mathcal{E}_X)$, $(Y,\mathcal{E}_Y)$ be coarse spaces,
$A\subseteq X$ and $f\colon A\to Y$ coarsely uniform.
Attaching $X$ to $Y$ using $f$ we get $X\cup_f Y := X\sqcup Y\, /\ a\!\sim\! f(a)$
and a quotient map $\pi\colon X\sqcup Y \to X\cup_f Y$.
We denote the coarse structure on the quotient by $\mathcal{E}_{X\cup_f Y}$.

\begin{lemma} \label{asdimattach}
Whenever $\asdim(X,\mathcal{E}_X) = \asdim(X\backslash A,\mathcal{E}_{X\cup_f Y}|_{X\backslash A})$,
we can conclude that $\asdim(X\cup_f Y) = \max\{\asdim(X),\asdim(Y)\}$.
\end{lemma}
\textbf{Proof.}
Note that $X\cup_f Y = X\backslash A \cup Y$.
Thus 
$\asdim(X\cup_f Y,\mathcal{E}_{X\cup_f Y})
   = \max\left\{\asdim(X\backslash A,\mathcal{E}_{X\cup_f Y}|_{X\backslash A}),
   \asdim(Y,\mathcal{E}_{X\cup_f Y}|_Y)\right\}$
by Proposition~\ref{asdimfinunions}.

In Lemma~\ref{csonY} we proved $\mathcal{E}_{X\cup_f Y}|_Y=\mathcal{E}_Y$.
Together with the assumption this implies the claim.
\qua

\clearpage
\thispagestyle{empty}\ \newpage 

\chapter{Asymptotic dimension of coarse cell complexes}

\subsection*{Coarse cell complexes}

In order to define coarse cell complexes, we need to describe their building blocks.
Compare the discussion in \cite{MitchAdd} for a better understanding of the following definition of a ray.

\begin{defn} \label{defray}\index{ray}
Let $\mathcal{E}$ be a connected coarse structure on $\RRR_+$
which is compatible with the standard topology on $\RRR_+$.
We call $R=(\RRR_+,\mathcal{E})$ a \emph{ray}
if the coarse structure $\mathcal{E}$ satisfies the following two conditions.
\begin{itemize}
\item
If $M$ and $N$ are entourages, the same is true for
$$M+N := \{(u+x,v+y) \mid (u,v)\in M,(x,y)\in N\}\ .$$
\item
If $M$ is an entourage, so is
$$M^{\compl} := \{(u,v)\mid (x,y)\in M\text{ and } (x\leq u\leq v\leq y \text{ or } y\leq v\leq u\leq x)\}\ .$$
\end{itemize}
\end{defn}

\begin{prop}\label{rayeucl}
By $\mathcal{E}_{eucl.}$ we denote the usual bounded coarse structure on $\RRR_+$.
The coarse space $(\RRR_+,\mathcal{E}_{eucl.})$ is a ray
and $\mathcal{E}_{eucl.}\subseteq\mathcal{E}$ whenever $(\RRR_+,\mathcal{E})$ is a ray.
\end{prop}
For a proof consult Proposition 2.5 of \cite{MitchAdd}.

\begin{prop}\label{maxray}
By $\mathcal{E}_\cdot$ we denote the continuously controlled coarse structure on $\RRR_+$
induced by the one-point compactification.
The coarse space $(\RRR_+,\mathcal{E}_\cdot)$ is a ray
and $\mathcal{E}\subseteq\mathcal{E}_\cdot$ whenever $(\RRR_+,\mathcal{E})$ is a ray.
\end{prop}
\textbf{Proof.}
The inclusion $\mathcal{E}\subseteq\mathcal{E}_\cdot$ is true
for every coarse structure which is compatible with the topology of $\RRR_+$.
It remains to prove that $(\RRR_+,\mathcal{E}_\cdot)$ is a ray.
This is straightforward.
\qua

\begin{defn} \label{defcell}
Let $R$ be a ray and $n\in\NNN$.
\begin{eqnarray*}
S^{n-1}R &:=& (R\; {\scriptscriptstyle\coprod}\, R)^n\times\{0\} \\
D^nR &:=& (R\; {\scriptscriptstyle\coprod}\, R)^n\times R
\end{eqnarray*}
We call $S^nR$ a \emph{coarse $n$-sphere} and $D^nR$ a \emph{coarse $n$-cell}.
\index{coarse n-cell@coarse $n$-cell}\index{coarse n-sphere@coarse $n$-sphere}
\end{defn}

There are two different coarse structures on $R\; {\scriptscriptstyle\coprod}\, R$.
There is the coproduct coarse structure $\mathcal{E}_{\scriptscriptstyle\coprod}$ which is not connected
and there is the connected coarse structure $\cncs(\mathcal{E}_{\scriptscriptstyle\coprod})$
generated by $\mathcal{E}_{\scriptscriptstyle\coprod}$.
The arguments in this chapter work for both coarse structures,
but only equipped with the connected coarse structure
the spaces $S^nR$ and $D^nR$ deserve the names sphere and cell respectively.

\begin{defn}\label{defccc}\index{coarse cell complex}\textbf{(coarse cell complex)} \\
We call $(Y,\mathcal{E})$ a \emph{coarse cell complex}
if it is obtained by ``inductively'' attaching coarse cells to a disjoint union of coarse cells.

More precisely, suppose that for all $k\in\NNN$
we are given a coarse space $(Y_k,\mathcal{E}_k)$ and a set $I_k$.
Furthermore, for all $k\in\NNN$ and $i\in I_k$ there is
\begin{itemize}
\item a number $n_{k,i}\in\NNN$,
\item a ray $R_{k,i}$ and
\item a coarsely uniform map $f_{k,i}\colon S^{n_{k,i}-1}R_{k,i} \to Y_{k-1}$.
\end{itemize}
Define $f_k := \coprod_{i\in I_k} f_{k,i}\ \colon \coprod_{i\in I_k} S^{n_{k,i}-1}R_{k,i} \to Y_{k-1}$
and observe that $f_k$ is coarsely uniform. If
$$Y_0 = \coprod_{i\in I_0} D^{n_{0,i}}R_{0,i}\qquad\text{ and }\qquad
Y_k = \left(\coprod_{i\in I_k} D^{n_{k,i}}R_{k,i}\right) \cup_{f_k} Y_{k-1}\ ,$$
then $\underrightarrow{\lim} (Y_k,\mathcal{E}_k)$ is a coarse cell complex.

Define $\delta(Y) := \sup \{ n_{k,i} \mid k\in\NNN, i\in I_k \}$
and call $\delta(Y)$ the \emph{cell dimension}\index{cell dimension} of the coarse cell complex $Y$.
\end{defn}

\subsection*{Asymptotic dimension of coarse cell complexes}

Our goal is to prove the following theorem.
\begin{thm}\label{asdimcoarsecomplex}\index{coarse cell complex}
Let $Y$ be a coarse cell complex. 
Then $\asdim Y = \delta(Y) + 1$ where $\delta(Y)$ is the cell dimension\index{cell dimension} of $Y$.
\end{thm}
In order to get Theorem~\ref{asdimcoarsecomplex}, we start with some lemmas.

\bigskip

Let $R$ be a ray, $X'$ any coarse space and $A'\subseteq X'$.
Define $X:=D^{n-1}R\sqcup X'$ and $A:=S^{n-2}R\sqcup A'$ and denote the coarse structure of $X$ by $\mathcal{E}_X$.
Notice that $R^n\subseteq D^{n-1}R$. We will consider the set
$$P_n:=\left\{ (x_1,\ldots,x_n)\in R^n \mid x_n>0 \text{ and } x_i\leq x_n\text{ for } i<n \right\}
   \subseteq X\backslash A.$$

\begin{lemma}\label{csonPn}
For any coarse space $(Y,\mathcal{E}_Y)$ and any coarsely uniform map
$f\colon (A,\mathcal{E}_X|_A)\to (Y,\mathcal{E}_Y)$, we have
$\mathcal{E}_X|_{P_n} = \mathcal{E}_{X\cup_f Y}|_{P_n}$.
\end{lemma}
\textbf{Proof.}
It follows directly from the definition of the quotient coarse structure
that $\mathcal{E}_X|_{P_n} \subseteq \mathcal{E}_{X\cup_f Y}|_{P_n}$.
We need to prove the converse inclusion.
Let $E\in\mathcal{E}_{X\cup_f Y}|_{P_n}$.
There are symmetric entourages $L_i\in\mathcal{E}_{D^{n-1}R}$ containing the diagonal $\Delta_{D^{n-1}R}$
such that $E \subseteq L_1^A\cdots L_k^A$ where $L_i^A:=L_i\cup{A\!\times\! A}$.
If we set $L:=L_1\cdots L_k$, we get $E \subseteq L_1^A\cdots L_k^A \subseteq L\cup L(A\!\times\! A)L$.
Thus, it remains to prove $L(A\!\times\! A)L\, \cap\, P^2_n\in\mathcal{E}_X |_{P_n}$.

Since the projection $p_n\colon P_n\to R, (x_1,\ldots,x_n)\mapsto x_n$ is coarsely uniform,
the image $\widetilde{L}:=p_n\times p_n(L)$ of $L$ is in $\mathcal{E}_R$
and compatibility with the topology implies $\widetilde{L}[0]\subseteq [0,r_E]$.
We conclude  $L(A\!\times\! A)L\subseteq\{ (x_1,\ldots,x_n) \mid x_n\leq r_E \}^2$
and hence $L(A\!\times\! A)L \cap P^2_n \subseteq \left([0,r_E]^n\right)^2 \in \mathcal{E}_X|_{P_n}$.
\qua

\begin{lemma}\label{asdimPngeqn}
$\asdim(P_n,\mathcal{E}_X|_{P_n}) \geq n$
\end{lemma}
\textbf{Proof.}
Let $\mathcal{U}$ be a uniformly bounded cover of $P_n$ with appetite
$\Delta_{P_n,1} := \{(x,y)\in P_n\times P_n \mid d_{eucl.}(x,y)\leq 1\}$.
Note that by Proposition~\ref{rayeucl} the set $\Delta_{P_n,1}$ is an entourage.
We will prove that the multiplicity of $\mathcal{U}$ is at least $n\! +\! 1$.

By $p_i\colon P_n\to R, (x_1,\ldots,x_n)\mapsto x_i$,
we denote the projection of $P_n$ onto the $i$-th factor.
Set $E_i:=p_i\times p_i(\Delta_\mathcal{U})\in\mathcal{E}_R$ and $\Delta_1 := \{ (x,y)\in R^2 \mid \abs{x-y} \leq 1 \}$.
Choose $r>1$ such that $E_n\Delta_1 E_{n-1}E_{n-2}\cdots E_2E_1[0]\subseteq\left[0,r\right)$.
Define $$a_n:=(0,\ldots,0,1)\in P_n\qquad \text{ and }\qquad
a_j:=(\underbrace{0,\ldots,0}_{j},r,\ldots,r)\in P_n.$$
The points $a_0,\ldots,a_n$ are the vertices of the $n$-simplex
$$S:=\{\lambda_0 a_0 + \cdots + \lambda_n a_n
   \mid \lambda_0,\ldots,\lambda_n\geq 0 \text{ and } \lambda_0 +\cdots + \lambda_n = 1\}.$$
By $S_j$ we denote the $(n\! -\! 1)$-face of $S$ opposite to $a_j$.
For $x=(x_1,\ldots,x_n)$ the following conclusions hold.
\begin{eqnarray*}
x\in S_0 & \Rightarrow & x_1=0 \\
x\in S_j & \Rightarrow & x_j=x_{j+1} \qquad\qquad\qquad (j = 1,\ldots,n\! -\! 2) \\
x\in S_{n-1} & \Rightarrow & 0\leq x_n-x_{n-1}\leq 1 \\
x\in S_n & \Rightarrow & x_n=r
\end{eqnarray*}
We claim that for each $U\in\mathcal{U}$ there is a face $S_i$
such that $U\cap S_i=\emptyset$.
Suppose there is $U\in\mathcal{U}$ with $U\cap S_i\neq\emptyset$
for all $i\in\{0,\ldots,n\}$.
Take an element $x^{(0)}\in U\cap S_0$.
Then $x^{(0)}_1=0$ and hence $p_1(U)\subseteq E_1[0]$.
Now take $x^{(1)}\in U\cap S_1$.
We have $x^{(1)}_2=x^{(1)}_1$ and therefore $p_2(U)\subseteq E_2E_1[0]$.
Inductively, we get $p_j(U)\subseteq E_j\cdots E_1[0]$ for $j\in\{1,\ldots,n-1\}$.
Finally, we take $x^{(n-1)}\in U\cap S_{n-1}$.
This time we have $x^{(n-1)}_n\in \Delta_1(x^{(n-1)}_{n-1})$ and hence
$p_n(U)\subseteq E_n\Delta_1E_{n-1}\cdots E_1[0]\subseteq\left[0,r\right)$.
On the other hand, since $U\cap S_n\neq \emptyset$, we get $r\in p_n(U)$.
This is a contradiction.

Of course $\mathcal{U}_S:=\{ U\in\mathcal{U}\mid U\cap S\neq\emptyset \}$ is a cover of $S$.
Take any map $g\colon\mathcal{U}_S\to\{a_0,\ldots a_n\}$,
such that $U\in\mathcal{U}_S$ does not intersect the $(n\! -\! 1)$-face opposite to $g(U)$.
Defining $U_i := \bigcup_{U\in g^{-1}(a_i)}U$ we obtain the cover $\{U_0,\ldots U_n\}$ of $S$
with Lebesgue number at least $1$ and such that $U_i$ contains $a_i$ and $U_i\cap S_i=\emptyset$.

We intend to apply Theorem~\ref{SpernersLemma}. Therefore,
take a triangulation $K$ of $S$ such that the diameter of any simplex of $K$ is at most $1$.
Choose a map $h\colon \operatorname{vert}(K) \to \{a_0\ldots, a_n\}$
such that for each $v\in h^{-1}(a_i)$ the star of $v$,
i.e. the union of all simplices containing $v$, is contained in $U_i$.
Applying Theorem~\ref{SpernersLemma} we get an $n$-simplex of $K$ which is contained in $U_0\cap\cdots\cap U_n$.
Hence $U_0\cap\cdots\cap U_n \neq \emptyset$.
This yields $\mu(\mathcal{U})\geq n+1$.
\qua

\begin{thm}[Sperner's Lemma]\label{SpernersLemma}\index{Sperner's Lemma}
Let $K$ be a triangulation of a closed $n$-simplex $S$.
Suppose we are given a map $h\colon vert(K) \to vert(S)$ with the following property:
If $v\in vert(K)$ lies in any face of $S$, then $h(v)$ lies in the same face.
Under these assumptions, there is at least one $n$-simplex $T$ of the triangulation $K$
whose $n+1$ vertices correspond under $h$ to the $n+1$ different vertices of $S$.
\end{thm}
A proof of Sperner's Lemma can be found in \cite{Fedorchuk}.
Sperner's Lemma is a combinatorial version of the fact
that the sphere $S^{n-1}$ is not a retract of the closed ball of radius $1$ around the origin in $\RRR^n$.

\begin{lemma}\label{asdimXleqn}
We can estimate the asymptotic dimension of coarse spheres and coarse cells as follows.
Let $R$ be a ray and $n\in\NNN$. Then
$$ \asdim D^{n-1}R\leq n \quad\text{ and }\quad \asdim S^{n-1}R\leq n. $$
\end{lemma}
\textbf{Proof.}
Notice that Proposition~\ref{asdiminfcoprod} yields $$\asdim D^{n-1}R = \asdim R^n = \asdim S^{n-1}R .$$
For $i\in\{1,\ldots,n\}$ consider the projection $p_i\colon R^n\to R, (x_1,\ldots,x_n) \mapsto x_i$ and 
remember that $\Delta_1:=\{(x,y)\in R\times R\mid \abs{x-y}\leq 1\}\in\mathcal{E}_R$.
Let $E\in\mathcal{E}^{*n}_R$ be a symmetric entourage.
Set $E_i:=p_i\times p_i(E)\cup \Delta_1\in\mathcal{E}_R$ and
$\widetilde{E}:=(E_1\cup\cdots\cup E_n)^{\compl}\in\mathcal{E}_R$.
This implies $E\subseteq \widetilde{E}^{\times n}\in\mathcal{E}_R^{*n}$.
We define some bounded sets out of which we will build a cover of $R^n$.
\begin{eqnarray*}
K_0 &:=& \{0\} \\ 
K_{i+1} &:=& \widetilde{E}[K_i]\qquad\ \text{ for }i\in\NNN \\
K_i &:=& \emptyset\qquad\qquad\ \text{ for }i<0 \\
\kappa_i &:=& \sup K_i\qquad\text{ for } i\in\NNN
\end{eqnarray*}
Observe that $\left[0,\kappa_i\right)\subseteq K_i$. 
Now we define a uniformly bounded cover $\mathcal{U}$ of $R^n$
consisting of $n+1$ families each of them being $E$-disjoint.
\begin{eqnarray*}
U_i &:=& K_i\ \backslash\ K_{i-n} \\
\mathcal{U}_i &:=& \{U_{i_1}\times\cdots\times U_{i_n}\mid i_1,\ldots,i_n\equiv i\text{ mod }n+1\} \\
\mathcal{U} &:=& \mathcal{U}_0\ \cup\ \cdots\ \cup\ \mathcal{U}_n 
\end{eqnarray*}
We first prove that $\mathcal{U}$ is indeed a cover of $R^n$.
Let $(x_1,\ldots,x_n)\in R^n$. Choose $i_1,\ldots,i_n\in\NNN$ such that
$(x_1,\ldots,x_n)\in K_{i_1}\backslash K_{i_1-1} \times \cdots \times K_{i_n}\backslash K_{i_n-1}$.
There is $i\in\{1,\ldots,n+1\}$ such that for all $k\in\{1,\ldots,n\}$
we have $i_k \not\equiv i\text{ mod }n+1$.
Define $j_k := \min\{j\in\NNN\mid j\geq i_k\text{ and }j-n\equiv i\text{ mod }n+1\}$.
It follows that $j_k-n < i_k \leq j_k$ and hence $K_{i_k}\backslash K_{i_k-1}\subseteq\mathcal{U}_{j_k}$.
This implies $(x_1,\ldots,x_n)\in\mathcal{U}_{j_1}\times\cdots\times\mathcal{U}_{j_n}\in\mathcal{U}_{i-1}$.

Since $\Delta_1\subseteq\widetilde{E}$, the closure of $\widetilde{E}$
(with respect to the usual topology on $\RRR^2$) is contained in $\widetilde{E}^3$.
This implies $(\kappa_i,\kappa_{i+1})\in\widetilde{E}^3$ for $i\in\NNN$
and $(x,\kappa_i)\in\widetilde{E}^3$ for $i\in\NNN$ and $x\in K_i\backslash K_{i-1}$.
Hence $\Delta_{\mathcal{U}}\subseteq \left( \widetilde{E}^{\times n} \right)^{3n+6}$,
i.e. $\Delta_{\mathcal{U}}$ is uniformly bounded.

Observe that for $i\in\{0,\ldots,n\}$ the family $\mathcal{U}_i$ is $E$-disjoint.
This is a direct consequence of the definitions.
\qua

\begin{cor} \label{asdimcoarsecells}
If $R$ is a ray, then $\asdim D^nR = \asdim S^nR = n+1$.
\end{cor}
\textbf{Proof.}
Apply Lemma~\ref{asdimPngeqn} and Lemma~\ref{asdimXleqn}.
\qua

\begin{lemma}\label{laststep}
Let $I$ be a set. For $i\in I$ let $R_i$ be a ray and $n_i\in\NNN$.
Define $X := \coprod_{i\in I}D^{n_i}R_i$ and $A := \coprod_{i\in I}S^{n_i-1}R_i$.
Let $(Y,\mathcal{E}_Y)$ be any coarse space and $\mathcal{E}_X$ the coarse structure of $X$.
Set $n := \sup \{\ n_i \mid i\in I \}$.
If $f\colon (A,\mathcal{E}_X|_A) \to (Y,\mathcal{E}_Y)$ is coarsely uniform, then
$$ \asdim\left(X\backslash A,\mathcal{E}_{X\cup_f Y}|_{X\backslash A}\right)
\leq \max \left\{ n+1, \asdim\left(\im f, \mathcal{E}_Y|_{\im f}\right) \right\} . $$
\end{lemma}
\textbf{Proof.}
Suppose $n<\infty$ and set
$m := \max \left\{ n+1, \asdim\left(\im f, \mathcal{E}_Y|_{\im f}\right) \right\}$.
For every entourage $\breve{E}\in\mathcal{E}_{X\cup_f Y}|_X$ there is
$E\in\mathcal{E}_X$ and $L\in\mathcal{E}_Y|_{\im f}$ such that
$\breve{E} \subseteq E \cup E\breve{L}E$ where $\breve{L} := (f\times f)^{-1}(L)$.
Note that far away from $A$ the entourages $E$ and $\breve{E}$ coincide,
i.e. $E\backslash E[A]^2 =\breve{E} \backslash E[A]^2$.
We may assume $\Delta_A\subseteq \breve{L}$ and $\Delta_X\subseteq E$.
Note that $D^{n_i}R_i$ consists of at most $2^n$ copies of $R_i^{n_i+1}$.
Choose one of these copies for each $i\in I$ and denote their union by $X'$.
Furthermore, set $A' := A\cap X'$.
With Proposition~\ref{asdimfinunions} in mind it is enough to prove that
$ \asdim\left(X'\backslash A', \mathcal{E}_{X\cup_f Y}|_{X'\backslash A'}\right) \leq m $.

Similarly as in Lemma~\ref{asdimXleqn} we can choose $E_i\in\mathcal{E}_{R_i}$ for $i\in I$
such that $E|_{R_i^{n_i+1}}\subseteq E_i^{\times(n_i+1)}$ and $E_i = E_i^{\compl}$.
Since $X$ is a coproduct, Proposition~\ref{propcdcoprod} implies
that we can choose $E_i$ in such a way that there is a finite set $I_0\subseteq I$ such that
$E_i\subseteq\Delta_{R_i}$ for all $i\in I\backslash I_0$.
Moreover, we may assume $\Delta_1 := \{(x,y)\in\RRR^2_+ \mid \abs{x-y} \leq 1\}\subseteq E_i$ for $i\in I_0$.

Define
\begin{eqnarray*}
K_{i,0} &:=& \{0\} \!\;\qquad\qquad\text{ for } i\in I_0, \\
K_{i,k} &:=& \emptyset \!\;\quad\qquad\qquad\text{ for } k<0, i\in I_0 \text{ and } \\
K_{i,k+1} &:=& E_i^5[K_{i,k}] \qquad\text{ for } k\in\NNN,i\in I_0.
\end{eqnarray*}
Furthermore, we define a uniformly bounded cover
$\mathcal{U} = \mathcal{U}_0\cup\cdots\cup\mathcal{U}_m$ of $X'$ as follows:
Let $k\in\NNN$, $l\in\{1,\ldots,m+1\}$ and $i\in I_0$.
\begin{eqnarray*}
U_{i,k} &:=& E_i^2\left[ K_{i,k}\backslash K_{i,k-m} \right] \\
\mathcal{U}_l &:=& \{U_{i,k_1}\times\cdots\times U_{i,k_{n_i+1}} \mid i\in I_0, k_1,\ldots,k_{n_i+1}\equiv l \text{ mod } m+1 \}
\end{eqnarray*}
So far we did not cover the parts of $X'$ corresponding to $I\backslash I_0$.
Define $$\mathcal{U}_0 \ := \ \mathcal{U}_{m+1} \cup \{ \{ x \} \mid x\in R_i^{n_i+1}, i\in I\backslash I_0 \} . $$
The proof that $\mathcal{U}$ is a uniformly bounded cover of $X'$
is similar to an argument in the proof of Lemma~\ref{asdimXleqn}.
Observe that $\mathcal{U}$ has appetite $(E|_{X'})^2$ and that each of the families
$\mathcal{U}_0,\ldots,\mathcal{U}_m$ is $E$-disjoint.

We want to construct a nice cover $\mathcal{V}$ of $A'$.
For $M\subseteq X\times X$ define $f_*M := f\times f(M\cap A\!\times\! A)$.
Notice that $f_*E, f_*\Delta_\mathcal{U}\in\mathcal{E}_Y$.
Take a uniformly bounded cover
$\widetilde{\mathcal{V}}:=\widetilde{\mathcal{V}}_0\cup\cdots\cup\widetilde{\mathcal{V}}_m$
of $\im f$ with appetite $f_*E\cup (f_*E)L(f_*E)$
and such that each family $\widetilde{\mathcal{V}}_k$ is $f_*\Delta_\mathcal{U}$-disjoint.
By taking the inverse image under $f|_{A'}$
we get a cover $\mathcal{V} = \mathcal{V}_0\cup\cdots\cup\mathcal{V}_m$ of $A'$ with appetite
$(E\breve{L}E)|_{A'}=(E\cup E\breve{L}E)|_{A'}$ and such that each of the families
$\mathcal{V}_0,\ldots,\mathcal{V}_m$ is $\Delta_\mathcal{U}$-disjoint.
In general $\mathcal{V}$ will not be uniformly bounded with respect to $\mathcal{E}_X|_{A'}$.

Now we construct a cover of $E[A]\cap X'\backslash A'$ starting from the cover $\mathcal{V}$ of $A'$.
Let $k\in\{0,\ldots,m\}$. For $V\in\mathcal{V}_k$ define
$$ V^\flat := \bigcup_{i\in I} (V\cap R_i^{n_i}) \times E_i[0] \qquad\text{and}\qquad
V^\sharp := \left(V^\flat \cup \bigcup_{{U\in\mathcal{U}_k} \atop {U\cap V^\flat\neq\emptyset}}U\right) \backslash A \ . $$
Observe that $\mathcal{V}_k^\sharp := \{ V^\sharp \mid V\in\mathcal{V}_k \}$
is uniformly bounded with respect to $\mathcal{E}_{X\cup_f Y}|_{X'\backslash A'}$.
From $\mathcal{U}_k$ and $\mathcal{V}_k^\sharp$ we build the family
$$ \mathcal{W}_k := \{ U\backslash A \mid U\in\mathcal{U}_k, U\cap V^\flat = \emptyset \text{ for all } V\in\mathcal{V}_k \}\
                    \cup \mathcal{V}_k^\sharp . $$
Note that any $U\in\mathcal{U}_k$ does intersect $V^\flat$ for at most one $V\in\mathcal{V}_k$.
We obtain a cover $\mathcal{W} := \mathcal{W}_0 \cup \cdots \cup \mathcal{W}_m$ of $X'\backslash A'$
which is uniformly bounded with respect to $\mathcal{E}_{X\cup_f Y}|_{X'\backslash A'}$
and has appetite $\breve{E}|_{X'\backslash A'}$.
\qua

\bigskip

\textbf{Proof of Theorem~\ref{asdimcoarsecomplex}.}
Let $Y_k$, $I_k$, $R_{k,i}$, $n_{k,i}$ and $f_{k,i}$ be as in Definition~\ref{defccc}.
By induction we will prove
\begin{equation}\label{finccc}
\asdim Y_k = \sup\{n_{j,i} \mid j\leq k, i\in I_j\} + 1.
\end{equation}

The case $k=0$ follows from Proposition~\ref{asdiminfcoprod} and Corollary~\ref{asdimcoarsecells}.

Set $X:=\coprod_{i\in I_k}D^{n_{k,i}}R_{k,i}$ and $A:=\coprod_{i\in I_k}S^{n_{k,i}-1}R_{k,i}$.
By Corollary~\ref{asdimcoarsecells} the asymptotic dimension of $X$ is $s_k+1$
where $s_k:=\sup\{n_{k,i}\mid i\in I_k\}$.
We need information on the asymptotic dimension of
$\left(X\backslash A,\ \mathcal{E}_{X\cup_{f_k} Y_{k-1}}|_{X\backslash A}\right)$.
Observe that $$X\backslash A=\coprod_{i\in I_k}D^{n_{k,i}}R_{k,i}\backslash S^{n_{k,i}-1}R_{k,i}\ .$$

Using Lemma~\ref{csonPn}, Lemma~\ref{asdimPngeqn} and monotony of asymptotic dimension,
we get $\asdim\left(X\backslash A,\ \mathcal{E}_{X\cup_{f_k} Y_{k-1}}|_{X\backslash A}\right) \geq s_k+1$.
Applying Lemma~\ref{laststep} we get
$\asdim\left(X\backslash A,\ \mathcal{E}_{X\cup_{f_k} Y_{k-1}}|_{X\backslash A}\right) \leq 
   \max\{s_k+1, \asdim Y_{k-1}\}$.
Lemma~\ref{csonY} assures us that there is no ambiguity about the coarse structure of $Y_{k-1}$.
Notice that $Y_k = X\cup_{f_k} Y_{k-1} = X\backslash A \cup Y_{k-1}$.
Hence, Proposition~\ref{asdimfinunions} completes the proof of (\ref{finccc}).

Since all the maps $Y_k \to Y$ are injective,
the claim of Theorem~\ref{asdimcoarsecomplex} now follows from Proposition~\ref{asdimcolim}.
\qua

\bigskip

We conclude this chapter by discussing some possible modifications
in the definition of a coarse cell complex and the effect on the asymptotic dimension.

\begin{rem}
We might consider building a coarse cell out of different rays, e.g.
$R_1 {\scriptscriptstyle\coprod}\, R_1 \times \cdots \times R_n {\scriptscriptstyle\coprod}\, R_n \times R_{n+1}$.

Using Proposition~\ref{maxray}, the arguments given in this chapter can be easily adjusted
to prove $\asdim Y = \delta(Y) + 1$
for a coarse cell complex $Y$ built with modified cells.
\end{rem}

\begin{rem}
If we want to use the disjoint union $\bigsqcup$ instead of the coproduct $\coprod$ in Definition~\ref{defccc},
we have to demand the maps $f_k$ to be coarsely uniform, since this is no longer automatic.

If for all $k\in\NNN$ the map $f_k$ is coarsely uniform,
we can prove the formula of Theorem~\ref{asdimcoarsecomplex} in the same way as before.
In the proof of Lemma~\ref{laststep} we do not have to distinguish
two subsets of $I$. We may define $I_0 := I$ in this case.
We do not need the fact that $I_0$ is a finite set if we work with disjoint unions instead of coproducts.
\end{rem}

\subsection*{Metric coarse cell complexes}

If $(X,d_X)$ and $(Y,d_Y)$ are pseudometric spaces, $A\subseteq X$ 
and $f\colon A\to Y$ is a (not necessarily continuous) map,
then $X\cup_f Y$ is defined as a pseudometric space.
We will use $\mathbb{S}^n := \RRR^{n+1}\times \{0\}$ and $\mathbb{D}^n := \RRR^n\times \RRR_+$
with the restrictions of the euclidean metric as coarse $n$-sphere and coarse $n$-cell respectively.

\begin{defn}
Let $I$ be any set and $X_i$ and $Y_i$ pseudometric spaces.
A family of maps $\{f_i\colon X_i\to Y_i\}_{i\in I}$ is said to be
\emph{coarsely uniform in a uniform way}
if there is a (not necessarily continuous) map $s\colon\RRR_+\to\RRR_+$
such that
$$ \underset{r>0}{\forall}\quad \underset{i\in I}{\forall}\quad \underset{x,\widetilde{x}\in X_i}{\forall}\quad
d(x,\widetilde{x})\leq r\ \Rightarrow\ d(f_i(x),f_i(\widetilde{x}))\leq s(r) . $$
\end{defn}

\begin{defn}\label{defmccc}\index{metric coarse cell complex}
Suppose that for all $k\in\NNN$ there is a pseudometric space $(Y_k,d_k)$ and a set $I_k$.
Moreover, suppose that for all $k\in\NNN$ and $i\in I_k$ there is a number $n_{k,i}\in\NNN$
and a map $f_{k,i}\colon \mathbb{S}^{n_{k,i}-1} \to Y_{k-1}$
such that for each $k\in\NNN$ the family $\{f_{k,i}\}_{i\in I_k}$ is coarsely uniform in a uniform way.
Define $f_k := \bigsqcup_{i\in I_k} f_{k,i}\ \colon \bigsqcup_{i\in I_k} \mathbb{S}^{n_{k,i}-1} \to Y_{k-1}$
and observe that $f_k$ is coarsely uniform. If
$$Y_0 = \bigsqcup_{i\in I_0} \mathbb{D}^{n_{0,i}}\qquad\text{ and }\qquad
Y_k = \left(\bigsqcup_{i\in I_k} \mathbb{D}^{n_{k,i}}\right) \cup_{f_k} Y_{k-1}\ ,$$
then $(Y,d) := \underrightarrow{\lim} (Y_k,d_k)$ is a \emph{metric coarse cell complex}.
\end{defn}

There are two coarse structures on $Y$.
We write $\mathcal{E}_d$ for the bounded coarse structure of $Y$ induced by the pseudometric $d$.
Similarly, we write $\mathcal{E}_{d_k}$ for the bounded coarse structure of $Y_k$ induced by the pseudometric $d_k$.
Denote the coarse structure of $\underrightarrow{\lim} (Y_k,\mathcal{E}_{d_k})$ by $\mathcal{E}$.

Using Proposition~\ref{limpor}, we see that $\mathcal{E}\subseteq\mathcal{E}_d$.
The following example yields that in general $\mathcal{E}$ and $\mathcal{E}_d$ do not coincide.
\begin{example}
Consider any metric coarse cell complex $(Y,d)$ with $I_k\neq\emptyset$ for all $k\in\NNN$.
Define $\Delta_1 := \left\{ (y_1,y_2)\in Y\times Y \mid \abs{y_1-y_2}\leq 1 \right\}$
and observe that $\Delta_1\in\mathcal{E}_d$.
Proposition~\ref{limpor} implies $\Delta_1\not\in\mathcal{E}$.
\end{example}

\begin{rem}\label{remmccc}
The arguments given in this chapter can be easily adjusted to prove
$$ \asdim Y_k = \sup\{n_{j,i} \mid j\leq k, i\in I_j\} + 1 $$
for metric coarse cell complexes.
Here we use the notations of Definition~\ref{defmccc}.
\end{rem}

\begin{prop}
$\asdim(Y,\mathcal{E}) = \delta(Y) + 1$.
\end{prop}
\textbf{Proof.}
The claim follows from Remark~\ref{remmccc} and Proposition~\ref{asdimcolim}.
\qua

\begin{prop}\label{asdimfmccc}
If the metric coarse cell complex $Y$ can be built in finitely many steps,
i.e. if there is $k_0\in\NNN$ such that $I_k=\emptyset$ for all $k\geq k_0$,
then $\asdim(Y,\mathcal{E}_d) = \delta(Y) + 1$.
\end{prop}
Proposition~\ref{asdimfmccc} follows directly from Remark~\ref{remmccc}.

\bigskip

The following example of a metric coarse cell complex yields
that the formula of Theorem~\ref{asdimcoarsecomplex} is not true in general.
\begin{example}
We inductively define a metric coarse cell complex $(Y,d)$
and start with $Y_0 := \RRR^2\times\RRR_+ = \mathbb{D}^2\subseteq\RRR^3$.
For each $k\in\NNN\backslash\{0\}$ and $i\in I := \{0\}\times\ZZZ\times\NNN\subseteq\RRR^3$
we attach a coarse $1$-cell $\mathbb{D}^1$ to $Y_0$ using the map
$$ \widetilde{f}_{k,i} \colon\ \mathbb{S}^0 \to Y_0, \ (x,0) \mapsto (2^kx,0,0) + i\ . $$
Denote the metric space obtained after attaching a cell for all $k\in\{1,\ldots,n\}$ and $i\in I$
as just explained by $Y_n$.
Observe that the inclusion $i_k\colon Y_0\hookrightarrow Y_k$ is not isometric.
On the other hand, all the attaching maps $f_{k,i} := i_{k-1}\circ \widetilde{f}_{k,i}$
are $2$-Lipschitz.\footnote{
   Let $X$ and $Y$ be pseudometric spaces.
   A map $f\colon X\to Y$ is called \emph{$\lambda$-Lipschitz}\index{Lipschitz}
   if $d(f(x),f(x')) \leq \lambda\cdot d(x,x')$ for all $x,x'\in X$.
   The map $f$ is called \emph{contracting}\index{contracting} if it is $1$-Lipschitz. }
Hence $Y:=\underrightarrow{\lim} Y_k$ is a metric coarse cell complex.

The metric coarse cell complex $Y$ is coarsely equivalent to the consecutively defined metric coarse cell complex $Y'$:
We start with $Y'_0 := \RRR\times\RRR_+ = \mathbb{D}^1$
and for each $k\in\NNN\backslash\{0\}$ and $y\in\ZZZ\times\NNN\subset Y'$ we attach a coarse cell $\mathbb{D}^1$
using the constant attaching map $f'_{k,y}\colon\mathbb{S}^0\to Y'_0$
whose image is just $\{y\}$.

Finally, we obtain $\asdim(Y,\mathcal{E}_d) = \asdim(Y') = 2 < 2+1 = \delta(Y) +1$.
\end{example}

Note that in the previous example the inclusion $Y_k \hookrightarrow Y$ is not a coarse embedding.

\begin{rem}
Let $Y$ be a metric coarse cell complex such that all attaching maps $f_{k,i}$ are contracting.
In this case $(Y_k,d_k)$ is a subspace of $(Y,d)$, i.e. the inclusion map $Y_k \hookrightarrow Y$ is isometric.
Hence $\asdim(Y,\mathcal{E}_d) \geq \delta(Y)+1$.
\end{rem}
\textbf{Proof.}
Compare Proposition~\ref{asdim_ofR}.
\qua


\begin{question}
Does $\asdim Y \leq \delta(Y)+1$ hold for any metric coarse cell complex $Y$?
Does the inequality hold if all attaching maps are contracting?
\end{question}

\clearpage
\thispagestyle{empty}\ \newpage 

\chapter{Asymptotic dimension and the covering dimension of the Higson corona}

We will need dimension theory of topological spaces in this section,
in particular some knowledge about the covering dimension $\dim$.
Compare \cite{HuWa}, \cite{Engelking} and \cite{Fedorchuk}.

We refer to the beginning of Section~\ref{secHigson}
for the definitions of the Higson compactification and the Higson corona.

\bigskip

For a proper metric space $X$, we can recover the bounded coarse structure from its Higson compactification $hX$,
since the continuously controlled coarse structure with respect to $hX$ and the bounded coarse structure coincide.
(This is Proposition 2.47 of \cite{RoeCG}.)

Dranishnikov, Keesling and Uspenskij proved in \cite{DKU} and \cite{Dran} the following theorem.

\begin{thm}\label{thmDKU}
Let $(X,d)$ be a metric space of finite asymptotic dimension and $\nu X$ its Higson corona.
Then $\asdim(X,d) = \dim(\nu X)$.
\end{thm}

The following examples show that Theorem~\ref{thmDKU} cannot be generalized to arbitrary coarse spaces.

\begin{example}\index{Stone-\v Cech compactification}\index{trivial coarse structure}
Let $X$ be a non-compact, normal Hausdorff space which can be written as a union of countably many compact sets.
The Stone-\v Cech compactification $\beta X$ induces on $X$ the trivial coarse structure $\mathcal{T}$,
and the Higson compactification of $(X,\mathcal{T})$ is homeomorphic to $\beta X$.

We get the following results about dimensions:
On the one hand we know that $\asdim(X,\mathcal{T})=0$.
On the other hand $\dim(\beta X\backslash X)=0$
if and only if $\beta X\backslash X$ is totally disconnected.
If $X$ has the property that for every compact subset $K\subseteq X$ the space $X\backslash K$ is connected,
then $\beta X\backslash X$ is also connected.
Hence, for such $X$ we have $\dim(\beta X\backslash X)\neq 0$.
\end{example}
\textbf{Proof.}
Obviously $\mathcal{T}$ is contained in the coarse structure induced by any compactification.
We write $\mathcal{E}_{\beta X}$ for the coarse structure on $X$ induced by $\beta X$.
It remains to prove that $\mathcal{E}_{\beta X} \subseteq \mathcal{T}$.
There are compact sets $K_0\subset K_1\subset K_2\subset\cdots\subset X$
such that $X=\bigcup_{i\in\NNN}K_i$. Let $E\in\mathcal{E}_{\beta X}$.
Assume that for all $i\in\NNN$ there is $(x_i,y_i)\in E \backslash (\Delta_X\cup K_i\!\times\! K_i)$.
This implies $x_i\neq y_i$.
We may assume that $\{x_i\mid i\in\NNN\}\cap\{y_i\mid i\in\NNN\}=\emptyset$,
since we can choose subsequences with this property.
Observe that $\{x_i\mid i\in\NNN\}\subseteq X$ and $\{y_i\mid i\in\NNN\}\subseteq X$ are closed sets.
Let $(x,y)\in(\beta X\backslash X)^2$ be an element in the closure of the sequence $\{(x_i,y_i)\}_{i\in\NNN}$.
Since the closure of $E$ in $\beta X\times \beta X$ is contained in $X\!\times\! X\cup \Delta_{\beta X}$,
we conclude $x=y$.
Using Urysohn's Lemma, we get a continuous function $f\colon X\to [0,1]$ with
$f(x_i)=1$ and $f(y_i)=0$ for all $i\in\NNN$.
We may as well think of $f$ as a continuous function with domain $\beta X$.
Since inverse images of closed sets are closed, we get $f(x)=1$ and $f(y)=0$.
This is a contradiction to $x=y$.
Hence $E\subseteq \Delta_X\cup K\!\times\! K$ for a compact set $K\subseteq X$,
i.e. $E\in\mathcal{T}$.

Every bounded, continuous function on $X$ is a Higson function
with respect to $\mathcal{T}$. Hence, the Higson compactification of $(X,\mathcal{T})$ is $\beta X$.

Compare Example~\ref{tcs}
and \cite{Fedorchuk}, Chapter 1, Section 3.2, Theorem 5 for the statements on dimensions.

Assume that $\beta X\backslash X$ is the disjoint union of two open and non-empty sets $A$ and $B$.
Then the characteristic function $\chi_A$ is a continuous projection
in $C(\beta X\backslash X) \cong C_b(X) / C_0(X)$.\footnote{
   Let $M$ be  a locally compact Hausdorff space.
   We denote the space of continuous functions on $M$ with values in $\CCC$ by $C(M)$.
   The space of continuous, bounded functions on $M$ with values in $\CCC$ is called $C_b(M)$
   and by $C_0(M)\subseteq C_b(M)$ we denote the subspace of those functions vanishing at infinity.}
Choose a function $f\in C_b(X)$ whose image under $\pi\colon C_b(X)\to C_b(X) / C_0(X)$ is $\chi_A$.
Observe that $g:=\frac{1}{2}(f+f^*)$ is a self-adjoint preimage of $\chi_A$ with respect to $\pi$.
In particular $g$ is real-valued.
Consider $g$ as an element of $C(\beta X)$.
There is a compact set $K\subseteq X$ such that $g(x)\neq\frac{1}{2}$ for all $x\in X\backslash K$.
Hence $$X\backslash K
= \left\{x\in X\backslash K\mid g(x)<\frac{1}{2}\right\}\cup \left\{x\in X\backslash K\mid g(x)>\frac{1}{2}\right\}$$
is not connected.
\qua

\begin{example}
Let $X$ be a proper metric space and $\mathcal{E}_\cdot$ the continuously controlled coarse structure
induced by the one-point compactification of $X$. In this case $\nu X$ is just a single point and
$\asdim(X,\mathcal{E}_\cdot) = 1 = \dim(\nu X) + 1$.
\end{example}
\textbf{Proof.}
Compare Example 9.7 of \cite{RoeCG} for the calculation of asymptotic dimension.
If $X=\RRR^n$, we can apply Example~\ref{examplern}.
\qua

\begin{example}
Consider $\RRR^n$ with the visual corona $S^{n-1}$ and denote by
$\mathcal{E}_{vis}$ the coarse structure coming from the corresponding compactification.
Then $\asdim(\RRR^n,\mathcal{E}_{vis}) = n = \dim(S^{n-1}) + 1$.
\end{example}
\textbf{Proof.}
Example~\ref{visualdot}, the calculation of the asymptotic dimension of coarse cells
and Proposition~\ref{asdimfinunions} yield the desired formula for the asymptotic dimension.
\qua

\bigskip

Let $hX$ be a metrisable compactification of $X$
and denote the corresponding continuously controlled coarse structure by $\mathcal{E}_{hX}$.
The last two examples suggest a relation between the asymptotic dimension of $(X,\mathcal{E}_{hX})$
and the covering dimension of the corona $\nu X := hX\backslash X$.
The goal of this chapter is to prove the following theorem.

\begin{thm}\label{asdimdim}
$$\asdim(X,\mathcal{E}_{hX}) = \dim(\nu X) + 1$$
\end{thm}

The proof of Theorem~\ref{asdimdim} is given by the next two lemmas.

\begin{lemma}\label{lemmausingdimth}
$$\dim(\nu X) + 1 \leq \asdim(X,\mathcal{E}_{hX})$$
\end{lemma}
\textbf{Proof.}
Consider the map $f\colon \nu X\times\RRR_+ \to \nu X\times\NNN,\ (x,t) \mapsto (x,[t])$ and
denote by $\mathcal{E}'$ the pull-back of the coarse structure $\mathcal{E}$ described in Lemma~\ref{descriptE}.
Observe that $f$ is a coarse equivalence.
The coarse structure $\mathcal{E}'$ has a description
which is very similar to the description of $\mathcal{E}$ in Lemma~\ref{descriptE}.
Theorem~\ref{thmXnuXN} now implies
that $(X,\mathcal{E}_{hX})$ and $(\nu X\times \RRR_+,\mathcal{E}')$ are coarsely equivalent.

Let $\asdim(X,\mathcal{E}_{hX})=n$.
We have a uniformly bounded cover of $\nu X\times\RRR_+$ with appetite
$$ E := \left\{ \big((x,t),(y,s)\big)\in(\nu X\times\RRR_+)^2 \mid d(x,y) < \frac{1}{\max\{t,s\}},\ \abs{t-s} < 1 \right\} $$ 
and multiplicity $\leq n+1$.
Since $E((x,t))$ is open for all $(x,t)\in \nu X\times\RRR_+$, this cover has an open refinement
which also has appetite $E$ and multiplicity at most $n+1$.
Denote this cover by $\mathcal{U}$.

Let $\varepsilon > 0$. From $\mathcal{U}$ we will construct an open cover
$\mathcal{U}_\varepsilon$ of $\nu X\times [0,1]$
with multiplicity $\leq n+1$ and $\mesh(\mathcal{U}_\varepsilon)\leq\varepsilon$.
This implies $\dim(\nu X\times [0,1])\leq n$.

We start constructing $\mathcal{U}_\varepsilon$.
Observe that $L := \pi_{\RRR_+}\times\pi_{\RRR_+} (\Delta_\mathcal{U})\subseteq \RRR_+^2$
is symmetric and contains $\Delta_1 := \{(t,s)\in\RRR_+^2 \mid \abs{t-s}<1\}$.
Choose $t_\varepsilon >0$ such that
$d(x,y)\leq\frac{\varepsilon}{2}$ for all $((x,t),(y,s))\in\Delta_\mathcal{U}$
with $\max\{s,t\}\geq t_\varepsilon$.
Consider the continuous map $f_\varepsilon\colon\RRR_+\to\RRR_+$ with
$f_\varepsilon(0) = t_\varepsilon$ and
$$ f_\varepsilon\left(\frac{n\cdot\varepsilon}{4}\right)
 = \sup L^n(t_\varepsilon) \quad\text{ for }\quad n=\{1,2,\ldots\} $$
which is affine on the interval $\left[\frac{n\cdot\varepsilon}{4},\frac{(n+1)\cdot\varepsilon}{4}\right]$
for each $n\in\NNN$. Observe that $f_\varepsilon$ is strictly monotone increasing.
Define $\mathcal{U}_\varepsilon$ to be the pull-back of $\mathcal{U}$
under the continuous map $\id_{\nu X}\times f_\varepsilon$.
It is clear that $\mathcal{U}_\varepsilon$ has the same multiplicity as the cover $\mathcal{U}$.

Take $V\in\mathcal{U}$ and set $U:=(\id_{\nu X}\times f_\varepsilon)^{-1}(V)$.
It remains to check that the diameter of $U$ is at most $\varepsilon$.
Let $(x,t),(y,s)\in U$.
The definition of $t_\varepsilon$ implies $d(x,y)\leq\frac{\varepsilon}{2}$.
Suppose $s\leq t$.
Notice that $f_\varepsilon(t)\in L(f_\varepsilon(s))$.
This implies $\abs{t-s}\leq\frac{\varepsilon}{2}$.
Together, we get $d((x,t),(y,s))\leq\varepsilon$.
This completes the construction of the cover $\mathcal{U}_\varepsilon$.

Using some results about the covering dimension of products
(compare \cite{Fedorchuk}, Chapter 2, Section 6.3, Theorem 18 and Chapter 5, Section 5.3, Corollary 5),
we get $\dim(\nu X) + 1 = \dim(\nu X\times [0,1])\leq n$.
\qua

\begin{lemma}
$$ \asdim(X,\mathcal{E}_{hX}) \leq \dim(\nu X) + 1 $$
\end{lemma}
\textbf{Proof.}
Assume $\dim(\nu X) + 1 = n$.
We will prove $\asdim(\nu X\times \NNN, \mathcal{E})\leq n$
where $\mathcal{E}$ is the coarse structure described in Lemma~\ref{descriptE}.

Let $E\in\mathcal{E}$ be a symmetric entourage containing the diagonal $\Delta_{\nu X\times\NNN}$ and
$\{\delta_k\}_{k\in\NNN}$ a sequence of non-negative real numbers converging to zero
such that $d(x,x')<\delta_k$ for $((x,m),(x',m'))\in E$ with $\max\{m,m'\} \geq k$.
We may assume this sequence to be monotone.

For each $k\in\NNN_+$ there is a finite cover 
of $\nu X$
with open balls of radius $\frac{1}{2k}$ and multiplicity at most $n$.
Ostrand's theorem\footnote{
   We use the following version of Ostrand's theorem:
   A normal space $X$ satisfies $\dim X\leq n$
   if and only if for every locally finite open cover $\mathcal{U}$
   there exists an open cover $\mathcal{V}=\mathcal{V}_1\cup\cdots\cup\mathcal{V}_{n+1}$
   where each of the families $\mathcal{V}_i$ consists of disjoint sets
   and for each $V\in\mathcal{V}$ there exists $U\in\mathcal{U}$ such that $V\subseteq U$.
   Compare \cite{Engelking}.
}
provides a finite open cover $\mathcal{V}_k$ of $\nu X$
with $\mesh(\mathcal{V}_k)\leq\frac{1}{k}$ consisting of $n$ disjoint families
$\mathcal{V}_{k,1},\ldots,\mathcal{V}_{k,n}$.
For $V\in\mathcal{V}_{k,j}$ define $j_V := j$.

Since $\nu X$ is compact, $\mathcal{V}_k$ has positive Lebesgue number $L(\mathcal{V}_k)>0$.
We define inductively a sequence $\{l_i\}_{i\in\NNN}$. Set $l_0:=1$ and for $i>0$ take $l_i\in\NNN$
such that $\frac{1}{l_i}\leq L(\mathcal{V}_{l_{i-1}})$ and $l_i > l_{i-1}$.
We choose maps $\phi_i\colon\mathcal{V}_{l_i}\to\mathcal{V}_{l_{i-1}}$ such that
$W\subseteq\phi_i(W)$ for all $W\in\mathcal{V}_{l_i}$.
For $V\in\mathcal{V}_{l_{i-1}}$ we define $\breve V :=\bigcup_{W\in\phi_i^{-1}(V)} W$.
Observe that $\breve V\subseteq V$.

Consider the projection $\pi_\NNN\colon \nu X\times \NNN\to \NNN$
and define $E_\NNN:=\pi_\NNN\times\pi_\NNN(E)\cup\Delta_1$
where $\Delta_1:=\left\{(i,j)\in\NNN^2\mid \abs{i-j}\leq 1 \right\}$.
Observe that $E_\NNN$ is an entourage of the continuously controlled coarse structure
induced by the one-point compactification of $\NNN$.
Set $K_0:=\{0\}$ and $K_k:=E_\NNN[K_{k-1}]$
and notice that $K_0\subset K_1\subset\cdots\subset\NNN$ and $\bigcup_{k\in\NNN}K_k = \NNN$.

We define a sequence $\{k_i\}_{i\in\NNN}$. Set $k_{-2}:=0$.
For all $i\in\ZZZ$ with $i\geq -1$ choose $k_i\in\NNN$ such that $k_i>k_{i-1}+2n$ and
$\delta_m <\frac{1}{l_{i+1}}$ for all $m\not\in K_{k_i-2}$.

We need some last preparations. For $i\in\NNN\backslash\{0\}$ and $j\in\{1,\ldots,n\}$ we set
$A_{i,j} := K_{k_i}\backslash K_{k_{i-1}+2j-2}$ and $B_{i,j} := K_{k_i+2j}\backslash K_{k_i}$.

Finally, we define a cover $\mathcal{U}$ of $\nu X\times \NNN$.
\begin{eqnarray*}
U_0 &:=& \bigcup_{V\in\mathcal{V}_{l_0}} V\times K_{k_0+2j_V} \\
\mathcal{U}_i &:=& \left\{V\times A_{i,j_{\phi_i(V)}}\ \cup\ \breve V\times B_{i,j_V} \mid V\in\mathcal{V}_{l_i}\right\} \\
\mathcal{U} &:=& \{U_0\}\ \cup\bigcup_{i\in\NNN\backslash\{0\}}\mathcal{U}_i
\end{eqnarray*}
To check that $\mathcal{U}$ is uniformly bounded, let $((x,m_x),(y,m_y))\in\Delta_\mathcal{U}$.
Suppose that $m_x,m_y\not\in K_{k_{-1}}$.
There are natural numbers $k_x$ ,$k_y$, $i_x$ and $i_y$ such that the following is true.
\begin{eqnarray*}
m_x \ \in \ K_{k_x}\backslash K_{k_x-1} & \subseteq & K_{k_{i_x}}\backslash K_{k_{i_x-1}} \\
m_y \ \in \ K_{k_y}\backslash K_{k_y-1} & \subseteq & K_{k_{i_y}}\backslash K_{k_{i_y-1}}
\end{eqnarray*}
We have $\abs{i_x-i_y}\leq 1$. Hence $\pi_\NNN\times\pi_\NNN(\Delta_\mathcal{U})\in\mathcal{E}_\NNN$
where $\mathcal{E}_\NNN$ is the continuously controlled coarse structure on $\NNN$
induced by the one-point compactification.
Without loss of generality, assume that $m_x\leq m_y$.
Set $d_m := 1$ if $m\in K_{k_1}$.
Set $d_m := \frac{1}{l_{i-2}}$ if $m\in K_{k_i}\backslash K_{k_{i-1}}$ and $i>1$.
Observe that $\lim_{m\to\infty}d_m=0$
and $d(x,y)\leq \frac{1}{l_{i_y-2}} = d_{m_y}$.

Now we prove that $\mathcal{U}$ has multiplicity at most $n+1$.
Let $(x,m)\in\nu X\times\NNN$.
If $m\in K_{k_0}$, then $(x,m)\in U_0$
and $(x,m)$ is not contained in any other covering set.
Suppose that $m\in K_{k_i}\backslash K_{k_{i-1}}$ with $i\in\NNN\backslash\{0\}$.
If $(x,m)\in U$ for some $U\in\mathcal{U}$, then $U\in\mathcal{U}_i\cup\mathcal{U}_{i-1}$.
Observe that $x$ is contained in at most one set of each of the collections
$\mathcal{V}_{l_i,1},\ldots,\mathcal{V}_{l_i,n},\mathcal{V}_{l_{i-1},1},\ldots,\mathcal{V}_{l_{i-1},n}$.
Moreover, $m\not\in A_{l,j}$ if $l\neq i$, and $m\not\in B_{l,j}$ if $l\neq i-1$.
Hence $(x,m)$ is contained in at most $2n$ sets of the cover $\mathcal{U}$.
There is $l\in\NNN\backslash\{0\}$ such that $m\in K_{k_{i-1}+2l}\backslash K_{k_{i-1}+2l-2}$.
If $l\geq n$, then $m\not\in B_{i-1,j}$ for $j\in\{1,\ldots,n-1\}$.
If $1<l<n$, then $m\not\in B_{i-1,j}$ for $j\in\{1,\ldots,l-1\}$
and $m\not\in A_{i,j}$ for $j\in\{l+1,\ldots,n\}$.
If $l=1$, then $m\not\in A_{i,j}$ for $j\in\{2,\ldots,n\}$.
Hence, there are at most $n+1$ sets of the cover $\mathcal{U}$ containing $(x,m)$.

Finally we prove that $\mathcal{U}$ has appetite $E$.
Let $(x,m)\in\nu X\times\NNN$.
Set $K_k := \emptyset$ for $k<0$.
There is a unique $k\in\NNN$ such that $m\in K_k\backslash K_{k-1}$.
It follows that $\pi_\NNN(E((x,m)))\subseteq K_{k+1}\backslash K_{k-2}$.
There is a unique $i\in\{-1\}\cup\NNN$ with $k_{i-1}\leq k < k_i$.
We conclude that $\pi_\NNN(E((x,m)))\subseteq K_{k_i}\backslash K_{k_{i-1}-2}$.
If $i\leq 0$, then $E((x,m))\subseteq \nu X\times K_{k_0}\subseteq U_0$.

Suppose $i>0$.
Observe that $\pi_{\nu X}(E((x,m)))$ is contained in the ball $B_{\delta_m}(x)$ of radius $\delta_m$ around $x$.
There is $W\in\mathcal{V}_{l_i}$ such that $B_{\delta_m}(x)\subseteq W$.
Set $V:=\phi_i(W)$ and observe that $W\subseteq V\in\mathcal{V}_{l_{i-1}}$.

It remains to prove that $W\times K_{k+1}\backslash K_{k-2}$ is contained in some $U\in\mathcal{U}$.

If $k\geq k_{i-1}+2j_V$, then
\begin{eqnarray*}
W\times K_{k+1}\backslash K_{k-2}
& \subseteq & W\times K_{k_i}\backslash K_{k_{i-1}+2j_V-2} \\
& = & W\times A_{i,j_V}
\ \subseteq\ W\times A_{i,j_V} \cup \breve W\times B_{i,j_W}
\ \in\ \mathcal{U}_i\ .
\end{eqnarray*}

If $k<k_{i-1}+2j_V$ and $i>1$, then
\begin{eqnarray*}
W\times K_{k+1}\backslash K_{k-2}
& \subseteq & W\times K_{k_{i-1}+2j_V}\backslash K_{k_{i-1}-2} \\
& \subseteq & V\times A_{i-1,j_{\phi_{i-1}(V)}} \cup \breve V \times B_{i-1,j_V}
\ \in\ \mathcal{U}_{i-1}\ .
\end{eqnarray*}

If $k<k_{i-1}+2j_V$ and $i=1$, then
$$ W\times K_{k+1}\backslash K_{k-2}
\ \subseteq\ V\times K_{k_0+2j_V}
\ \subseteq\ U_0 \ \in\ \mathcal{U}\ . $$

This completes the proof that the cover $\mathcal{U}$ has appetite $E$.
\qua

\begin{rem}\label{asdimprodthmsc}
For any coarse spaces $X$ and $Y$ whose coarse structures are induced by metrisable compactifications,
we get another proof of $$\asdim(X\times Y) \leq \asdim(X) + \asdim(Y).$$
\end{rem}
\textbf{Proof.}
By Corollary~\ref{coraboutPsi} the coarse space $X$ is coarsely equivalent to $\psi(\nu X)$ and the same for $Y$.
This and Proposition~\ref{prodPsi} imply that
$X\times Y$ is coarsely equivalent to $\psi(\nu X\boxtimes\nu Y)$.
Hence, the Higson corona of $X\times Y$ is homeomorphic to $\nu X\boxtimes\nu Y$.
Using Theorem~\ref{asdimdim} and the product theorem of dimension theory
(i.e. \cite{Fedorchuk}, Chapter 2, Section 6.3, Theorem 18), we get
$\asdim(X\!\times\! Y) = \dim(\nu X\boxtimes\nu Y)+1
= \dim(\nu X\!\times\!\nu Y\!\times\! [0,1])+1 \leq \dim(\nu X)+\dim(\nu Y)+2
= \asdim(X)+\asdim(Y)$.
\qua

\begin{cor}\label{asdimprodnoteq}
There are coarse spaces $X$ and $Y$ such that $$\asdim(X\times Y) < \asdim(X) + \asdim(Y).$$
\end{cor}
\textbf{Proof.}
There are metric compacta $K$ and $L$ such that $\dim(K\times L)<\dim(K)+\dim(L)$.
Compare the section about Pontryagin's compacta, i.e Capter 5, Section 4.2 of \cite{Fedorchuk}.
The same reasoning as in the proof of Remark~\ref{asdimprodthmsc} translates this to asymptotic dimensions.
\qua

\begin{rem}
In \cite{buyalo2} the authors construct metric spaces $X$ and $Y$
such that $\asdim(X\times Y) < \asdim(X) + \asdim(Y)$.
\end{rem}


\chapter{Asymptotic dimension of $\text{CAT}(\kappa)$-spaces for $\kappa < 0$}

Let $\kappa < 0$.
We write $M_\kappa^n$ for the $n$-dimensional hyperbolic space with constant sectional curvature $\kappa$.
We need some basic facts about $\CAT(\kappa)$-spaces and some knowledge in hyperbolic geometry.
A good introduction on these subjects can be found in \cite{BH}.

We want to calculate the asymptotic dimension of $\CAT(\kappa)$-spaces.
As the following example shows, we should not expect
that the asymptotic dimension and the covering dimension of a $\CAT(\kappa)$-space always coincide.
\begin{example}\label{asdimnotdim}
Consider an $n$-simplex $\sigma$ in $M^n_{\kappa}$.
Take a $k$-dimensional face $f$ of $\sigma$.
There is an isometric embedding $i\colon M^k_{\kappa} \hookrightarrow M^n_\kappa$
such that $f\subseteq \im i$.
Now $\sigma \cup \im i$ is a $\text{CAT}(0)$-space with covering dimension $n$
and asymptotic dimension $k$.
\end{example}

Hence, for every pair of natural numbers $(k,n)$ with $k\leq n$,
there is a $\CAT(\kappa)$-space with asymptotic dimension $k$ and covering dimension $n$.

\section{Spaces with nicely covered spheres}

Let $(X,d)$ be a metric space and $x_0\in X$ a basepoint.
For $r\geq 0$ we define $S_r(x_0) = \{x\in X\mid d(x,x_0) = r\}$ and $D_r(x_0) = \{x\in X\mid d(x,x_0) < r\}$.

We consider $S_r(x_0)\subseteq X$ and $D_r(x_0)\subseteq X$
with the metric obtained by restricting the metric $d$ of $X$.

\begin{defn} \label{ncoloredspheres}
We say that $X$ has \emph{nicely $n$-colored spheres} (with respect to the basepoint $x_0\in X$)
if there is $\rho > 0$ and a cover $\mathcal{U}_k$ of $S_{k\cdot\rho}(x_0)$ for all $k\in\NNN$ such that
\begin{itemize}
\item there is a lower bound $\lambda>0$ on the Lebesgue numbers of $\mathcal{U}_k$ for $k\in\NNN$,
\item there is an upper bound on the mesh of the covers $\mathcal{U}_k$ and
\item $\mathcal{U}_k=\mathcal{U}_{k,1}\cup\cdots\cup\mathcal{U}_{k,n}$
with multiplicity $\mu(\mathcal{U}_{k,i})\leq 1$ for $i\in\{1,\ldots,n\}$.
\end{itemize}
\end{defn}

\begin{defn} \label{ncoveredspheres}
We say that $X$ has \emph{nicely $n$-covered spheres} (with respect to the basepoint $x_0\in X$)
if there is $\rho > 0$ and a cover $\mathcal{U}_k$ of $S_{k\cdot\rho}(x_0)$ for all $k\in\NNN$ such that
\begin{itemize}
\item there is a lower bound $\lambda>0$ on the Lebesgue numbers of $\mathcal{U}_k$ for $k\in\NNN$,
\item there is an upper bound on the mesh of the covers $\mathcal{U}_k$ and
\item $\mathcal{U}_k$ has multiplicity at most $n$.
\end{itemize}
\end{defn}

\begin{lemma}\label{coloredeqcovered}
A metric space $X$ has nicely $n$-colored spheres if and only if it has nicely $n$-covered spheres.
\end{lemma}
\textbf{Proof.}
Compare Remark~\ref{LebesguezahlBuchfuehrung}.
\qua

\begin{thm}\label{CATthm}
Let $\kappa < 0$.
For a $\CAT(\kappa)$-space $X$ with nicely $n$-covered spheres,
we have $\asdim (X) \leq n$.
\end{thm}

\begin{rem}
In some cases the upper bound on the asymptotic dimension of $X$ given by Theorem~\ref{CATthm}
is lower than the covering dimension of $X$.
Compare Example~\ref{asdimnotdim} and
observe that in this case the upper bound provided by Theorem~\ref{CATthm} is sharp.
\end{rem}

\section{Proof of Theorem~\ref{CATthm}}

Because of Lemma~\ref{coloredeqcovered} we may assume that $X$ has nicely $n$-colored spheres.
We will prove Theorem~\ref{CATthm} by constructing
uniformly bounded covers $\mathcal{M}_{\rho,N}$ of $X$ with multiplicity at most $n+1$.
In a second step we prove that for every $L > 0$ we can adjust the parameters $\rho$ and $N$
so that $\mathcal{M}_{\rho,N}$ has Lebesgue number at least $L$.

\subsection*{Step 1: Constructing covers}

Let $x_0$ be a basepoint of $X$ and $\rho>0$ such that for $k\in\NNN$ the sphere $S_k:=S_{k\cdot\rho}(x_0)$
is covered by $\mathcal{U}_k=\mathcal{U}_{k,1}\cup\cdots\cup\mathcal{U}_{k,n}$
with Lebesgue number $L(\mathcal{U}_k)\geq\lambda$, $\mesh(\mathcal{U}_k)\leq D$
and $\mu(\mathcal{U}_{k,i})=1$ for all $i\in\{1,\ldots,n\}$.
Observe that $\mathcal{U}_0 = \{\{x_0\}\}$. 
For $k\in\NNN$ we define $D_k:=D_{k\cdot\rho}(x_0)$.

Since $X$ is a $\CAT(\kappa)$-space, the space $X$ is uniquely geodesic.
We denote the geodesic segment from $x$ to $y$ by $[x,y]$.

Two parameters will appear in the following construction:
the distance $\rho$ between spheres and
the magnitude of shift $N\in\NNN$ towards the basepoint.
Note that whenever we have a space with nicely $n$-colored spheres of distance $\rho$,
we can choose any integral multiple of $\rho$ as a new distance of spheres,
i.e. we may choose $\rho$ as large as we want.

For $k\in\NNN$ we define a map $\theta_k\colon X\backslash D_k \to S_k$
by assigning to $x\in X\backslash D_k$
the intersecting point of the geodesic segment $[x_0,x]$ with $S_k$.
\begin{lemma}
The map $\theta_k$ is contracting\index{contracting} for any $k\in\NNN$.
\end{lemma}
\textbf{Proof.}
For $x\in X$ set $r(x) = d(x_0,x)$.
Let $x,y\in X\backslash D_k$. Without loss of generality assume $r(x)\leq r(y)$.
Take $z\in [x_0,y]$ with $r(z)=r(x)$.
Consider the triangle with vertices $x_0$, $x$ and $y$.
Take a comparison triangle in $M_\kappa^2$ and denote its vertices by $\overline{x_0}$, $\overline{x}$ and $\overline{y}$.
Then $$d(\theta_k(x),\theta_k(y)) \leq d\left(\overline{\theta_k(x)},\overline{\theta_k(y)}\right)
\leq d\left(\overline{x},\overline{z}\right) \leq d\left(\overline{x},\overline{y}\right) = d(x,y)\ .$$
Here $\overline{\theta_k(x)}$ and $\overline{\theta_k(y)}$ are the points on the chosen comparison triangle
corresponding to $\theta_k(x)$ and $\theta_k(y)$ respectively.
\qua

\bigskip

For $k\in\NNN$ we choose a map $\Theta_k\colon \mathcal{U}_{k+n}\to \mathcal{U}_{k}$.
We will see in Corollary~\ref{corollarthetachoice} that for $\rho$ large enough
we may assume that $\Theta_k$ has the following property:
\begin{equation}\label{thetaprop}
\theta_k(U) \subseteq \Theta_k(U) \quad \text{ for } k\in\NNN\text{ and }\ U\in\mathcal{U}_{k+n}
\end{equation}

Let $N\in\NNN$, $k\in n\cdot\NNN$ and $i\in\{1,\ldots,n\}$.
For $U\in\mathcal{U}_{k,i}$ we define as follows:
\begin{eqnarray*}
A(U) & := & \left\{
   \begin{array}{ll}
      D_{N+n} & \text{ if } k=0 \\
      \theta_k^{-1}(U) \ \cap \ D_{k+N+n}\backslash D_{k+N+i-1} & \text{ if } k>0
   \end{array} \right. \\   
B(U) & := & \theta_k^{-1}(U) \ \cap \ D_{k+N+i}\backslash D_{k+N} \\
U^\# & := & A(U)\ \cup \bigcup_{V\in \Theta_k^{-1}(U)} B(V)
\end{eqnarray*}
We get the following cover of $X$.
$$ \mathcal{M}_{\rho,N} := \left\{ U^\# \mid k\in n\cdot\NNN, \ U\in\mathcal{U}_k \right\} $$
Besides $\rho$ and $N$ this cover does depend on the choice of the maps $\Theta_k$.

\begin{lemma}
The multiplicity of $\mathcal{M}_{\rho,N}$ is at most $n+1$.
Moreover, if the maps $\Theta_k$ have property~(\ref{thetaprop}),
then $\mathcal{M}_{\rho,N}$ is uniformly bounded.
\end{lemma}
\textbf{Proof.}
Let $x\in X$. Choose $k\in n\cdot\NNN$ and $j\in\{1,\ldots,n\}$
such that $x\in D_{k+N+j}\backslash D_{k+N+j-1}$.
Observe that for $i\in\{1,\ldots,n\}$ the point $\theta_k(x)$ belongs to at most one set from $\mathcal{U}_{k,i}$.
Hence $x\in A(U)$ for at most $j$ sets $U\in\mathcal{U}_k$.
Similarly $x\in B(U)$ for at most $n+1-j$ sets $U\in\mathcal{U}_k$.

Assume that the maps $\Theta_k$ have property~(\ref{thetaprop}).
Let $U\in\mathcal{U}_k$ Observe that for $x,y\in U^\sharp$ we have
$$ d(x,y) \leq d(x,U)+ \diam(U) + d(U,y) \leq (N+2n)\cdot\rho + D + (N+2n)\cdot\rho \ . $$
Hence $\mathcal{M}_{\rho,N}$ is uniformly bounded.
\qua

\subsection*{Step 2: Adjusting parameters}\label{adjparameters}

In order to conclude the proof of Theorem \ref{CATthm},
it remains to prove that for every $L>0$ we can choose the parameters $\rho$ and $N$
such that $\mathcal{M}_{\rho,N}$ has Lebesgue number at least $L$.

\begin{lemma} \label{adjustrho}
Let $\delta>0$ and $B\subseteq S_{\rho\cdot k+a}(x_0)$ with $\diam(B)<a$.
The map $\theta_k|_B \colon B \to S_k$ is $\delta$-Lipschitz
if $a > \frac{2}{\sqrt{-\kappa}}\cdot\max\left\{1,\log\frac{2}{\delta}\right\}$.
\end{lemma}
\textbf{Proof.}
Set $r=\rho\cdot k+\frac{a}{2}$. Let $x,y\in B$.
Consider the geodesic triangle with vertices $x_0$, $x$ and $y$.
Denote the vertices of a comparison triangle in $M_\kappa^2$
by $\overline{x_0}$, $\overline{x}$ and $\overline{y}$ respectively.
We write $\overline{\theta_k(x)}$ and $\overline{\theta_k(y)}$
for the points on the chosen comparison triangle
corresponding to $\theta_k(x)$ and $\theta_k(y)$ respectively.
By $\overline{x}_m$ we denote the point on $[\overline{x_0},\overline{x}]$
with distance $\frac{a}{2}$ to $\overline{x}$.
Similarly, we define a point $\overline{y}_m$ on the segment $[\overline{x_0},\overline{y}]$.
Let $\alpha$ be the angle at $\overline{x_0}$ between
$[\overline{x_0},\overline{x}]$ and $[\overline{x_0},\overline{y}]$.
By $\gamma_r$ we denote the circular arc parametrized by arc length
from $\overline{x}_m$ to $\overline{y}_m$ in $S_r(\overline{x_0})\subset M_\kappa^2$
and by $\gamma_{\rho\cdot k}$ the corresponding arc from $\overline{\theta_k(x)}$ to $\overline{\theta_k(y)}$
in $S_{\rho\cdot k}(\overline{x_0})$.

Observe
\begin{eqnarray*}
d(\theta_k(x),\theta_k(y)) & \overset{(1)}{\leq} & d\left(\overline{\theta_k(x)},\overline{\theta_k(y)}\right)
\ \overset{(2)}{\leq} \ \text{ length of } \gamma_{\rho\cdot k} \\
& \overset{(3)}{=} & \frac{\alpha}{\sqrt{-\kappa}} \cdot \sinh\left(\sqrt{-\kappa}\cdot\rho\cdot k\right) \\
& < & \frac{\alpha}{2\cdot\sqrt{-\kappa}}\exp\left(\sqrt{-\kappa}\cdot\rho\cdot k\right) \\
& = & \frac{\alpha}{2\cdot\sqrt{-\kappa}}\exp\left(\sqrt{-\kappa}\cdot r\right)
 \cdot \exp\left(-\sqrt{-\kappa}\cdot \frac{a}{2}\right) \\
& \overset{(4)}{<} & \frac{\alpha}{2\cdot\sqrt{-\kappa}}\exp\left(\sqrt{-\kappa}\cdot r\right)\cdot\frac{\delta}{2} \\
& \overset{(5)}{<} & \frac{\alpha}{2\cdot\sqrt{-\kappa}}\Big(\exp\left(\sqrt{-\kappa}\cdot r\right)
 - \exp\left(-\sqrt{-\kappa}\cdot r\right)\Big)\cdot\delta \\
& = & \frac{\alpha}{\sqrt{-\kappa}} \cdot \sinh\left(\sqrt{-\kappa}\cdot r\right)\cdot\delta
\ \overset{(3)}{=} \ \delta\cdot\text{ length of } \gamma_r \\
& \overset{(6)}{\leq} & \delta\cdot d\left(\overline{x},\overline{y}\right) \ = \ \delta\cdot d(x,y)\ .
\end{eqnarray*}

Since $X$ is a $\CAT(\kappa)$-space, we can conclude $(1)$.
The fact that $\gamma_{\rho\cdot k}$ is a path from $\overline{\theta(x)}$ to $\overline{\theta(y)}$
implies inequality $(2)$.
For $(3)$ compare Proposition 6.17 in Part~I of \cite{BH}
about the Riemannian metric of $M_\kappa^2$ in polar coordinates.

The assumption $a>\frac{2}{\sqrt{-\kappa}}\cdot\log\frac{2}{\delta}$
implies $\exp\left(-\sqrt{-\kappa}\cdot \frac{a}{2}\right)<\frac{\delta}{2}$
and hence $(4)$.

The assumption $a>\frac{2}{\sqrt{-\kappa}}$ yields $\sqrt{-\kappa}\cdot r>1$.
This implies
$\exp(\sqrt{-\kappa}\cdot r) < 2\cdot\left(\exp\left(\sqrt{-\kappa}\cdot r\right) - \exp\left(-\sqrt{-\kappa}\cdot r\right)\right)$
which proves inequality $(5)$.

Since $d(\overline{x},\overline{y})\leq a$,
the geodesic segment $[\overline{x},\overline{y}]$ lies in $M_\kappa^2 \backslash D_r(\overline{x_0})$.
For a moment consider $M_\kappa^2 \backslash D_r(\overline{x_0})$ and $S_r(\overline{x_0})$
as metric spaces with their length metric.
    \footnote{Compare \cite{BH} for details on length metrics.}
Note that the radial projection $M_\kappa^2\backslash D_r(\overline{x_0})\to S_r(\overline{x_0})$ with center $\overline{x_0}$
is a contracting map which maps the geodesic segment $[\overline{x},\overline{y}]$ onto $\gamma_r$. This implies $(6)$.
\qua

\bigskip

Remember that $\lambda$ is a lower bound on the Lebesgue number
and  $D$ an upper bound on the mesh of the covers $\mathcal{U}_i$ for $i\in\NNN$.

\begin{cor}\label{corollarthetachoice}
If $\rho\cdot n > D$ and $\rho\cdot n > \frac{2}{\sqrt{-\kappa}}\cdot\max\left\{1,\log\frac{2\cdot D}{\lambda}\right\}$,
then we can choose $\Theta_k$ such that
$\theta_k(U) \subseteq \Theta_k(U)$ for all $k\in\NNN$ and $U\in\mathcal{U}_{k+n}$.
\end{cor}
\textbf{Proof.}
Set $\delta=\frac{\lambda}{D}$.
Lemma~\ref{adjustrho} implies that $\theta_k$ maps any $U\in\mathcal{U}_{k+n}$
into a ball $B\subseteq S_k$ of radius $\lambda$. Since $\mathcal{U}_k$ has Lebesgue number $\lambda$,
we may choose $\Theta_k(U)$ such that $B\subseteq\Theta_k(U)$.
\qua

\begin{cor} Let $L>0$. 
Choose $\rho > 2\cdot L$ and such that Corollary~\ref{corollarthetachoice} applies.
Choose $N\in\NNN$ such that $(N-1)\cdot\rho > 2L$
and $(N-1)\cdot\rho > \frac{2}{\sqrt{-\kappa}}\cdot\max\left\{1,\log\frac{2\cdot L}{\lambda}\right\}$.
Under these assumptions $\mathcal{M}_{\rho,N}$ has Lebesgue number at least $L$.
\end{cor}
\textbf{Proof.}
Let $x\in X$ and let $j\in\NNN$ be the smallest natural number such that 
the ball $D_L(x)$ of radius $L$ around $x$ is contained in $D_{j+N}$.
We claim that $D_L(x)\subseteq X\backslash D_{j+N-2}$.
Assume $D_L(x)\cap D_{j+N-2}\neq\emptyset$.
This implies $D_L(x)\subseteq\{x\in X\mid d(x,D_{j+N-2})<2L\}\subseteq D_{j+N-1}$.
Here the last inclusion follows from $\rho>2L$.
But this is a contradiction to the definition of $j$.

Take $k\in n\cdot\NNN$ and $l\in\{1,\ldots,n\}$ such that $j=k+l$.
Since $\theta_{k+N-1}$ is contracting,
the image of $D_L(x)$ is contained in a ball $B\subseteq S_{k+N-1}$ of radius $L$.
Set $\delta=\frac{\lambda}{L}$. Lemma~\ref{adjustrho} implies that $\theta_k(B)$
is contained in a ball $B'\subseteq S_k$ of radius $\lambda$.
Let $U\in\mathcal{U}_{k,i}$ be a covering set containing $B'$.
Then $D_L(x)\subseteq\theta_k^{-1}(U)$.

Hence $D_L(x)\subseteq A(U)$ if $i<l\leq n$.
On the other hand, we get $D_L(x)\subseteq B(U)$ if $1<l\leq i$.
Suppose $l=1$. Set $V:=\Theta_{k-n}(U)$. Note that Corollary~\ref{corollarthetachoice} yields
$\theta_k^{-1}(U) \cap D_{k+N}\backslash D_{k+N-1} \subseteq A(V)$.
Hence $D_L(x)\subseteq V^\#$.
\qua

\bigskip

This completes the proof of Theorem~\ref{CATthm}.

\section{Triangulations of bounded distortion}

\subsection*{Simplicial complexes}

A (combinatorial) \emph{simplicial complex} $C$ consists of a set $\vertices(C)$
and a collection $\Sigma(C)$ of finite, non-empty subsets of $\vertices(C)$
such that $\Sigma(C)$ is closed under taking subsets, i.e.
$\sigma\in\Sigma(C)$ and $\tau\subseteq\sigma$ imply $\tau\in\Sigma(C)$.
The collection $\Sigma(C)$ must contain $\{v\}$ for every $v\in\vertices(C)$. 
An element $v\in\vertices(C)$ is called a \emph{vertex} of $C$ and
$\sigma\in\Sigma(C)$ is called a \emph{simplex} of $C$.
The dimension of $\sigma$ is defined to be the cardinality of $\sigma$ minus one.
The dimension of $C$ is defined to be the supremum of the dimensions of its simplices, i.e.
$\dim(C) = \sup\{\dim(\sigma)\mid\sigma\in\Sigma(C)\}$.
By $C^{(k)}$ we denote the $k$-skeleton of $C$, i.e. the simplicial complex obtained from $C$
by removing all simplices of dimension greater than $k$.
A simplicial complex is said to be \emph{locally finite}
if each vertex is contained in only finitely many simplices.

By $\delta_v\colon\vertices(C)\to\CCC$ we denote the function which maps $v$ to $1$
and all the other vertices to $0$.
We identify $\vertices(C)$ as follows with a subset of $\ell^2(\vertices(C))$:
$$\vertices(C)\hookrightarrow\ell^2(\vertices(C))\ , \quad v \mapsto \frac{1}{\sqrt{2}} \cdot \delta_v$$
Observe that the distance of any two vertices in $\ell^2(\vertices(C))$ is one.
We write $\abs{\sigma}$ for the convex hull of $\sigma$ in $\ell^2(\vertices(C))$
and call $\abs{\sigma}$ a geometric simplex.
The euclidean realization $\abs{C}$ of $C$ is defined to be the union
of the geometric simplices.

There are two natural metrics on $\abs{C}$.
The \emph{affine metric}\index{affine metric} is just the restriction
of the standard metric on $\ell^2(\vertices(C))$ to $\abs{C}$.
The \emph{geodesic metric}\index{geodesic metric} is the length metric
associated to the affine metric. We will consider $\abs{C}$ with the geodesic metric.
Note that restricted to a single simplex both metrics coincide.

We write $\Star(v)$ for the open \emph{star}\index{star} of $v\in\vertices(C)$, i.e.
the interior of the union of all geometric simplices $\abs{\sigma}\subseteq\abs{C}$ with $v\in\sigma$.

The \emph{stability}\index{stability} of a finite dimensional simplicial complex $C$
is defined to be the biggest natural number $k$ such that there exist two different $k$-simplices
$\sigma,\tau\in\Sigma(C)$ such that $\sigma\cap\tau$ is a $(k-1)$-simplex of $C$.

\begin{lemma} \label{starcover}
Let $C$ be a simplicial complex of dimension $n$ and stability $k$.
The cover $\mathcal{V}= \{ \Star(v) \mid v\in\vertices(C) \}$ of $\abs{C}$
is uniformly bounded and has multiplicity $n+1$.
Moreover $\mathcal{V}$ has Lebesgue number $$\lambda_k = \frac{1}{\sqrt{2\cdot k\cdot (k+1)}}\ .$$
\end{lemma}
\textbf{Proof.}
The mesh of the cover $\mathcal{V}$ is at most $2$, i.e. $\mathcal{V}$ is uniformly bounded.
The multiplicity of $\mathcal{V}$ is $n+1$, since a simplex of $C$ has at most $n+1$ vertices.

The Lebesgue number of $\mathcal{V}$ is exactly the distance
of the center of a $k$-simplex to a face of dimension $k-1$.
Hence, an elementary computation yields the claim.
\qua

\subsection*{Spaces with triangulations of bounded distortion}

\begin{defn}
A \emph{triangulation of bounded distortion}\index{bounded distortion}\index{triangulation of bounded distortion}
of a metric space $X$ is a tuple $(C,\psi,l_1,l_2)$
where $C$ is a simplicial complex, $l_1,l_2 > 0$ and $\psi\colon \abs{C}\to X$ is a homeomorphism such that
for every $\sigma\in\Sigma(C)$ the map $\psi|_{\abs{\sigma}}$ is $(l_1,l_2)$-bi-Lipschitz, i.e.
$$l_1\cdot d(x,y) \leq d(\psi(x),\psi(y)) \leq l_2\cdot d(x,y)$$
for all $x,y\in\abs{\sigma}$.
We say that the triangulation is of dimension $\dim(C)$.
\end{defn}

Note that according to our definition the simplicial complex
of a triangulation of bounded distortion does not have to be locally finite.

\begin{prop} \label{starcoverBD}
Let $X$ be a metric space and $(C,\psi,l_1,l_2)$ a triangulation of bounded distortion
such that $C$ is of dimension $n$ and stability $k$.
Under these assumptions $\mathcal{U} = \{ \psi(\Star(v)) \mid v\in\vertices(C) \}$
is a cover of $X$ with the following properties:
\begin{itemize}
\item The cover $\mathcal{U}$ is uniformly bounded.
      More precisely, the diameter of any $U\in\mathcal{U}$ is at most $2\cdot l_2$.
\item When restricted to a single connected component,
      the cover $\mathcal{U}$ has Lebesgue number $L(\mathcal{U})\geq l_1\cdot\lambda_k$
      where $\lambda_k$ is defined as in Lemma~\ref{starcover}.
\item The cover $\mathcal{U}$ has multiplicity $n+1$.
\end{itemize}
\end{prop}

For the rest of this section let $X$ be a $\text{CAT}(\kappa)$-space
which admits a triangulation of bounded distortion $(C,\psi,l_1,l_2)$ of dimension $n$.
Choose a basepoint $x_0\in X$.
As above we set $S_r := \{ x\in X \mid d(x,x_0) = r \}$.

\begin{rem}
By restricting the cover $\mathcal{U}$ defined in Proposition~\ref{starcoverBD}
to the spheres $S_r$ with $r>0$
we see that $X$ has nicely $(n+1)$-covered spheres.
Hence, Theorem~\ref{CATthm} yields $\asdim(X)\leq n+1$.
\end{rem}

The original goal is to prove $\asdim(X)\leq n$,
but this seems to be much more difficult - at least in full generality. 
However, if there is a radial structure for $X$ as defined below,
we do obtain $\asdim(X)\leq n$.

\begin{defn}
Let $\rho>0$. For $k\in\NNN$ define $C_k$ to be the smallest simplicial complex
containing as simplices all $\sigma\in\Sigma(C)$
such that $\psi(\abs{\sigma})$ has non-empty intersection with $S_{k\cdot\rho}$.
Let $D>0$. Assume that for all $k\in\NNN$
there is an equivalence relation on $\vertices(C_k)$ such that
\begin{itemize}
\item in each $n$-simplex of $C_k$ there are two equivalent vertices and
\item the number $D$ is an upper bound for the diameter of the image under $\psi$ of any equivalence class.
\end{itemize}
A family of equivalence relations as described above will be called
a \emph{radial structure}\index{radial structure} for $X$.
\end{defn}

One might try to get a radial structure by defining the endpoints of ``radial'' edges to be equivalent.

\begin{prop}
If $X$ has a radial structure, then $\asdim(X)\leq n$.
\end{prop}
\textbf{Proof.}
Consider the collection $\mathcal{U}_k = \{ \psi(\Star(v)) \mid v\in\vertices(C_k) \}$.
We say $U_1, U_2\in\mathcal{U}_k$ are equivalent if they contain equivalent vertices.
For $U\in\mathcal{U}_k$ define $\widetilde{U}$ to be the union
of all $V\in\mathcal{U}_k$ which are equivalent to $U$.
By restricting $\widetilde{\mathcal{U}}_k := \{ \widetilde{U} \mid U\in\mathcal{U}_k \}$ to $S_{k\cdot\rho}$
we get a uniformly bounded cover of $S_{k\cdot\rho}$
with multiplicity at most $n$ and Lebesgue number at least $L(\mathcal{U})$.
Applying Theorem~\ref{CATthm} completes this proof.
\qua

\begin{question}
Is there a $\text{CAT}(\kappa)$-space with a triangulation of bounded distortion
which does not admit a radial structure?
\end{question}

\section{Complete, simply connected manifolds with bounded, strictly negative sectional curvature}

Theorem~\ref{thmHadamard} has been conjectured by Mikhael Gromov.
Compare Example~$1.\text{E}'_1$ of \cite{Gromov}.

\begin{defn}
Let $X$ be a Riemannian manifold.
For $x\in X$ and a two-dimensional subspace $p\subseteq T_xX$,
we denote the sectional curvature of the plane $p$ by $K(p)\in\RRR$.

We say that the sectional curvature of $X$ is bounded between $\kappa'$ and $\kappa$
or simply that $X$ has \emph{bounded sectional curvature}\index{sectional curvature}\index{bounded sectional curvature}
if $\kappa' \leq K(p) \leq \kappa$
for all $x\in X$ and all two-dimensional subspaces $p\subseteq T_xX$.

We say that $X$ has \emph{strictly negative sectional curvature}\index{strictly negative sectional curvature}
if there is $\kappa < 0$ such that $K(p) \leq \kappa$
for all $x\in X$ and all two-dimensional subspaces $p\subseteq T_xX$.
\end{defn}

\begin{thm} \label{thmHadamard}
For a complete, simply connected Riemannian manifold $X$ with bounded, strictly negative sectional curvature
we have $$\asdim(X)=\dim(X).$$
\end{thm}

We will prove this theorem by applying Theorem~\ref{CATthm}.
A different approach (using the boundary at infinity and a new invariant called capacity dimension)
has been taken in \cite{buyalo1} and \cite{buyalo2} to prove similar results.
These results include $\asdim(X)\leq\dim(X)$ if in addition to our assumptions $X$ is cobounded, i.e.
if there is a bounded subset $B\subset X$ such that the orbit of $B$ under the isometry group of $X$ covers $X$.
Compare Theorem~1.1 of \cite{buyalo1} and Corollary~1.2 and Proposition~6.2 of \cite{buyalo2}.

Before proceeding with the proof of Theorem~\ref{thmHadamard},
we will look at the differential geometry of the spheres in $X$.

\begin{lemma}\label{diffgeomsph} Let $X$ be a complete, simply connected Riemannian manifold of dimension $n$
with bounded, strictly negative sectional curvature.
Then there are $\kappa_1,\kappa_2\in\RRR$ and $\rho>0$ such that for all $x_0\in X$ and $r \geq 1$
\begin{enumerate}
\item[(a)] the sectional curvature of the sphere $S_r(x_0) = \{ x\in X\mid d(x,x_0) = r \}$
   is bounded between $\kappa_1$ and $\kappa_2$ and
\item[(b)] the injectivity radius of $S_r(x_0)$ is at least $\rho$.
\end{enumerate}
\end{lemma}
\textbf{Proof.} Let $x_0\in X$ and $r\geq 1$.

First we calculate bounds for the sectional curvature of $S_r(x_0)$.

Suppose the sectional curvature of $X$ is bounded between $-\kappa'$ and $-\kappa$.
Define $f\colon X \to \RRR, x\mapsto d(x_0,x)$.
Note that for each $x\in X$ the second fundamental tensor
$S_{\operatorname{grad} f}\colon T_xS_r(x_0) \to T_xS_r(x_0)$
is just the restriction of the Hessian\footnote{
   If $f\colon X\to \RRR$ is a differentiable function on a Riemannian manifold $X$,
   $\nabla$ the Levi-Civita connection and $Y$ a tangent vector field on $X$,
   then $\operatorname{Hess}(f)(Y) = \nabla_Y\operatorname{grad} f$.}
$\operatorname{Hess}(f)$ to the directions tangential to the spheres.

Let $x\in X$ with $f(x)=r\geq 1$.
Denote the eigenvalues of $\operatorname{Hess}(f)$ at $x$ by $\tau_1(x),\cdots,\tau_{n-1}(x)$.
Using Riccati comparison arguments, we get
$$ \sqrt{\kappa} \leq \sqrt{\kappa}\cdot \coth\left(\sqrt{\kappa}\cdot r\right) \overset{(\flat)}{\leq}
\tau_i(x)\overset{(\sharp)}{\leq} \sqrt{\kappa'}\cdot \coth\left(\sqrt{\kappa'}\cdot r\right)
\leq \sqrt{\kappa'}\cdot \coth\left(\sqrt{\kappa'}\right) . $$
For $(\flat)$ compare 1.7.3 of \cite{karcher}.
For $(\sharp)$ compare 1.7.1 of \cite{karcher} or Sections 1.4 and 1.6 of \cite{meyer1}.
Hence, the norm of the second fundamental form is bounded in terms of $\kappa'$ and $\kappa$.
Using the Gau\ss{} equation (see \cite{doCarmo} or \cite{Jost}) in order to relate the curvature of $X$ and the spheres,
we get universal bounds
$$\kappa_1 := \kappa -\kappa' - \kappa'\cdot \coth^2\left(\sqrt{\kappa'}\right) \quad\text{and}\quad
\kappa_2 := -2\cdot\kappa + \kappa'\cdot\coth^2\left(\sqrt{\kappa'}\right)$$
for the sectional curvature of $S_r(x_0)$.

Applying Corollary 4 of \cite{injRad} yields that the injectivity radius of $S_r(x_0)$ is at least
$\rho := \frac{\pi}{\sqrt{\kappa'}\cdot \coth\left(\sqrt{\kappa'}\right)}$.
\qua

\bigskip

\textbf{Proof of Theorem~\ref{thmHadamard}.}
First note that $X$ is a $\CAT(\kappa)$-space.
Compare Theorem II.1A.6 and the Cartan-Hadamard Theorem (II.4.1) of \cite{BH}.
We need completeness of $X$ to apply the Cartan-Hadamard Theorem.

Lemma~\ref{diffgeomsph} yields that $M := \bigcup_{n\in\NNN}S_n(x_0)\subset X$
is a Riemannian manifold of bounded geometry,
i.e. $M$ has bounded sectional curvature and positive injectivity radius.
Hence, we may apply the following theorem.
\begin{thm}\label{thmattie}
If $M$ is an $n$-dimensional Riemannian manifold of bounded geometry,
then $M$ admits a triangulation of bounded distortion
and the dimension of this triangulation is $n$.
\end{thm}
Theorem~\ref{thmattie} is an immediate consequence of Theorem 1.14 in \cite{attie}.
Moreover, Theorem~\ref{thmattie} and Proposition~\ref{starcoverBD} imply
that $X$ has nicely covered spheres.
Applying Theorem~\ref{CATthm} now yields $\asdim(X)\leq\dim(X)$.

Suppose $\asdim(X) =: k < n := \dim(X)$.
This implies that there is a uniformly bounded, open cover $\mathcal{U}$ of $X$
with multiplicity at most $k+1$.
Set $D:=\mesh(\mathcal{U})$. Choose a basepoint $x_0\in X$.
Let $\varepsilon>0$. We consider the following diffeomorphism $f_\varepsilon \colon X\to X$.
Define $f_\varepsilon(x)$ to be the point on the geodesic segment from $x_0$ to $x$
with distance $\frac{\varepsilon}{D}\cdot d(x_0,x)$ to $x_0$.
Observe that $f_\varepsilon(x)$ is $\frac{\varepsilon}{D}$-Lipschitz,
since the corresponding map on $\RRR^n$ is $\frac{\varepsilon}{D}$-Lipschitz and $X$ is a $\text{CAT}(0)$-space.
Now the image of the cover $\mathcal{U}$ under $f_\varepsilon$ gives an open cover of $X$ with sets of diameter
at most $\varepsilon$ and with multiplicity at most $k+1$.
Applying Theorem~1.6.12 of \cite{Engelking} yields $\dim(X)\leq k < n$.
This is a contradiction. Hence $\asdim(X) \geq \dim(X)$.
\qua

\begin{rem}
We actually proved $\asdim(X)\geq\dim(X)$ for any complete, simply connected
Riemannian manifold with strictly negative sectional curvature.
\end{rem}

\begin{question}
Is there a complete, simply connected Riemannian manifold
with strictly negative, but unbounded sectional curvature
and $\asdim(X) > \dim(X)$?
\end{question}


\nocite{BD2, bartels, Yufin, MacL, petersen}
\bibliographystyle{amsalpha}
\bibliography{asdim}

\clearpage
\thispagestyle{empty}\ \newpage 

\printindex

\section*{Lebenslauf}

\vspace{0.5cm}

\subsection*{Pers\"onliche Daten}
\begin{tabular}{ll}
Name & Bernd Grave \\
Geburtsdatum & 6.\! Oktober 1975 \\
Geburtsort & Damme, Niedersachsen \\
Staatsangeh\"origkeit & deutsch
\end{tabular}

\subsection*{Schulausbildung}
\begin{tabular}{ll}
1982 -- 1988 & Grundschule und Orientierungsstufe in Holdorf \\
1988 -- 1995 & Gymnasiums in Damme \\
1995 & Abitur \\
\end{tabular}

\subsection*{Zivildienst}
\begin{tabular}{ll}
1995 -- 1996 & Zivildienst in einem Altenpflegeheim in Holdorf \\
\end{tabular}

\subsection*{universit\"are Ausbildung}
\begin{tabular}{ll}
1996 -- 1999 & Studium der Mathematik und der Physik an der \\
             & Westf\"alischen-Wilhelms-Universit\"at M\"unster \\
1998         & Vordiplom in Mathematik und in Physik \\
1999 -- 2000 & Studium der Mathematik an der Universidad \\
             & Complutense de Madrid \\
2000 -- 2001 & Studium der Mathematik an der Universit\"at M\"unster \\
2002         & Diplom in Mathematik \\
seit 2002    & Promotionsstudent im Graduiertenkolleg Gruppen und \\
             & Geometrie an der Georg-August-Universit\"at G\"ottingen \\
2004         & Forschungsaufenthalt an der Vanderbilt University \\
             & in Nashville, Tennessee \\
\end{tabular}

\vspace{1cm}

G\"ottingen, den 19.\! Dezember 2005

\vspace{0.5cm}

Bernd Grave

\end{document}